\documentclass[AMA,STIX1COL]{WileyNJD-v2}

\articletype{RESEARCH ARTICLE}

\usepackage{amsmath}
\usepackage{amssymb}
\usepackage{amsfonts}
\usepackage{amsbsy}

\usepackage{newlfont}

\usepackage{bm}
\usepackage{siunitx}

\usepackage{float}
\usepackage{placeins}
\usepackage{afterpage}

\DeclareMathOperator{\dive}{\nabla\cdot}

\begin{document}

\title{An efficient and accurate implicit DG solver for the incompressible Navier-Stokes equations}

\author[1]{Giuseppe Orlando*}

\author[2]{Alessandro Della Rocca}

\author[1]{Paolo Francesco Barbante}

\author[1]{Luca Bonaventura}

\author[1]{Nicola Parolini}

\authormark{G. ORLANDO \textsc{et al}}

\address[1]{\orgdiv{MOX - Dipartimento di Matematica}, \orgname{Politecnico di Milano}, \orgaddress{\state{Piazza Leonardo da Vinci 32, 20133 Milano}, \country{Italy}}}

\address[2]{\orgdiv{Global R\&D}, \orgname{Tenova S.p.A.}, \orgaddress{\state{Via Albareto 31, 16153 Genova}, \country{Italy}}}

\corres{*Giuseppe Orlando, MOX - Dipartimento di Matematica, Politecnico di Milano \\
		Piazza Leonardo da Vinci 32, 20133 Milano, Italy. \email{giuseppe.orlando@polimi.it}}


\abstract[Summary]{We propose an efficient, accurate and robust implicit solver for the incompressible Navier-Stokes equations, based on a DG spatial discretization and on the TR-BDF2 method for time discretization. The effectiveness of the method is demonstrated in a number of classical benchmarks, which highlight its superior efficiency with respect to other widely used implicit approaches. The parallel implementation of the proposed method in the framework of the \textit{deal.II} software package allows for accurate and efficient adaptive simulations in complex geometries, which makes the proposed solver attractive for large scale industrial applications.}


\maketitle

\pagebreak

\section{Introduction}
\label{sec:intro} \indent
 
The efficient numerical solution of the incompressible Navier-Stokes equations is one of the most relevant goals of computational fluid dynamics.
A great number of methods have been proposed in the literature, see for example, among many others, the reviews in \cite{quartapelle:2013, quarteroni:2008}. Since the seminal proposals \cite{chorin:1968, temam:1969}, projection methods \cite{guermond:2006} have become very popular for the time discretization of this problem. 
Several spatial discretization approaches have been proposed and finite volume techniques using unstructured meshes \cite{fletcher:1997} have become the state of the art for industrial applications, in particular when implemented in parallel software packages like OpenFoam \cite{chen:2014,jasak:2007,weller:1998}. Indeed, in previous work by one of the authors \cite{dellarocca:2018}, a wide range of projection methods was implemented in OpenFoam and their performance was compared, as a preliminary step towards the development of a computational fluid dynamics tool for combustion simulations of industrial interest.
On the other hand, high order finite elements, both in their continuous and discontinuous versions \cite{giraldo:2020, karniadakis:2005}, have gained increasing popularity in the academic community and also in many applications, but are still far from being the reference tool for industrial use. More specifically, Discontinuous Galerkin methods for the Navier-Stokes equations have been proposed by many authors, we refer for example to  \cite{bassi:2005, fehn:2019, fehn:2017, fehn:2018, giorgiani:2014, schotzau:2003}.
 
In this work, we seek to combine, on the one hand, accurate and flexible discontinuous finite element spatial discretizations, and on the other hand, efficient and unconditionally stable time discretizations, following an approach that has been shown to be quite successful for applications to numerical weather prediction in \cite{tumolo:2015, tumolo:2013}. Building on the experience of \cite{dellarocca:2018}, we propose an accurate, efficient  and robust projection method, based on the second order TR-BDF2 method \cite{bank:1985, hosea:1996, tumolo:2015}. This solver is implemented using discontinuous finite elements, in the framework of the numerical library \textit{deal.II} \cite{bangerth:2007}, in order to provide a reliable and easily accessible tool for large scale industrial applications. It is important to remark that time discretizations of the Navier-Stokes equations based on accurate implicit solvers have been proposed in a number of papers, see among many others \cite{bassi:2007, bassi:2015, rhebergen:2013, tavelli:2014, tavelli:2016}. The specific combination of techniques presented in this work does not entail major conceptual novelties with respect to any of the above references, but we claim that
it constitutes an optimal combination for the development of a second order $h-$adaptive flow solver that can be competitive for industrial applications with more conventional finite volume techniques. Furthermore, while the TR-BDF2 method is only second order in time, the wide range of simulations presented in \cite{tumolo:2015} show that this method still allows to achieve quite accurate results even when coupled to higher order discretizations in space.
The paper is organized as follows: the time discretization approach is outlined and discussed in Section \ref{sec:modeleq}. The spatial discretization is presented in Section \ref{sec:numeth}. Some implementation issues, the validation of the proposed method
and its application to a number of significant CFD benchmarks are reported in Section
\ref{sec:tests}. Some conclusions and perspectives for future work are described in Section \ref{sec:conclu}.

\section{The Navier-Stokes equations and the time discretization strategy}
\label{sec:modeleq} \indent
 
Let \(\Omega \subset \mathbb{R}^{d}, 2 \le d \le 3\) be a connected open bounded set with a sufficiently smooth boundary \(\partial\Omega\) and denote by \(\mathbf{x}\) the spatial coordinates and by \(t\) the temporal coordinate. We consider the classical unsteady incompressible Navier-Stokes equations, written in non-dimensional form as:

\begin{align}
\label{eq:ns_incomp}
\frac{\partial \mathbf{u}}{\partial t} + \dive\left(\mathbf{u} \otimes\mathbf{u}\right) + \nabla p &= \frac{1}{Re}\Delta\mathbf{u} + \mathbf{f} \nonumber \\
\dive\mathbf{u} &= 0   
\end{align}

for $ \mathbf{x} \in \Omega, t \in (0, T], $ supplied with suitable initial and boundary conditions. Here \(T\) is the final time, \(\mathbf{u}\) is the fluid velocity, \(p\) is the pressure divided by density and \(Re\) is the Reynolds number, which is usually defined as $ Re = UL/\nu, $ where $U$ denotes a reference value of the velocity magnitude, $L$ a reference length scale and $\nu $ the fluid kinematic viscosity. The velocity \(\mathbf{u}\) and the pressure \(p\) are coupled together by the incompressibility constraint in \eqref{eq:ns_incomp}, which leads, after space discretization, to a system of differential and algebraic equations whose numerical solution presents several difficulties widely discussed in the literature. Furthermore, in the specific case of projection methods, difficulties arise in choosing the boundary conditions to be imposed for the Poisson equation which is to be solved at each time step to compute the pressure, see e.g. the discussion in \cite{guermond:2006}.

An alternative that allows to avoid or reduce some of these problems is the so-called artificial compressibility formulation, originally introduced in \cite{chorin:1967}. In this formulation, the incompressibility constraint is relaxed and a time evolution equation for the pressure is introduced, which is characterized by an artificial sound speed $ c, $ so as to obtain

\begin{align}
\label{eq:ns_artcomp}
\frac{\partial\mathbf{u}}{\partial t} + \dive\left(\mathbf{u}\otimes\mathbf{u}\right) + \nabla p &= \frac{1}{Re}\Delta\mathbf{u} + \mathbf{f} 
\nonumber\\
\frac{1}{c^2}\frac{\partial p}{\partial t} + \dive\mathbf{u} &= 0.   
\end{align}

For the sake of simplicity, we shall only consider \(\mathbf{f} = \mathbf{0}\) and Dirichlet boundary conditions for the velocity, i.e., \(\mathbf{u}\rvert_{\partial\Omega} = \mathbf{u}_D(t)\), while we consider homogeneous Neumann boundary conditions for the pressure. While most commonly discretized by explicit methods, see e.g.  \cite{nithiarasu:2003, rogers:1990} among many others, implicit methods have also been applied to this formulation, see e.g. \cite{ekaterinaris:2004, merkle:1987, rahman:2008}.  

Our goal here is to extend the projection method based on the TR-BDF2 scheme introduced in  \cite{dellarocca:2018} for the formulation \eqref{eq:ns_incomp} to the time discretization of system \eqref{eq:ns_artcomp}. This allows to avoid the introduction of stabilization parameters and to exploit the special properties of the TR-BDF2 method, which will be reviewed here briefly. Our development is also inspired by the first order semi-implicit methods \cite{casulli:1984, dumbser:2016}, which were proposed originally for the compressible Navier-Stokes equations but which could also be applied in the pseudo-incompressible case.
Introducing a discrete time step $ \Delta t = T/N $ and discrete time levels $ t^n = n\Delta t,$ $n= 0,\dots,N, $ for a generic time dependent problem \( \bm{u}^{\prime}=\mathcal{N}(\bm{u})\) the incremental form of the TR-BDF2 method can be described in terms of two stages, the first from \(t^n\) to \(t^{n+\gamma} = t^n + \gamma \Delta t\) and the second from \(t^{n+\gamma}\) to \(t^{n+1}\), which can be written as:

\begin{align}
\label{eq:trbdf2}
\frac{\bm{u}^{n+\gamma} - \bm{u}^{n}}{\gamma\Delta t} &= \frac{1}{2}\mathcal{N}\left(\bm{u}^{n+\gamma}\right) + \frac{1}{2}\mathcal{N}\left(\bm{u}^{n}\right) \\
\frac{\bm{u}^{n+1} - \bm{u}^{n + \gamma}}{\left(1 - \gamma\right)\Delta t} &= \frac{1}{2 - \gamma}\mathcal{N}\left(\bm{u}^{n+1}\right) + \frac{1 - \gamma}{2\left(2 - \gamma\right)}\mathcal{N}\left(\bm{u}^{n+\gamma}\right) + \frac{1 - \gamma}{2\left(2 - \gamma\right)}\mathcal{N}\left(\bm{u}^{n}\right). \nonumber
\end{align}

Here, \(\bm{u}^{n}\) denotes the approximation at time  \( n = 0,...,N\).
Notice that, in order to guarantee L-stability, one has to choose \(\gamma = 2 - \sqrt{2}\).
This second order implicit method, originally introduced in \cite{bank:1985} as a combination of the Trapezoidal Rule (or Crank-Nicolson) method and of the Backward Differentiation Formula method of order 2, has been fully analyzed in \cite{hosea:1996}.
While we will use here its original formulation, the method was shown in \cite{hosea:1996} to be an L-stable Explicit first step, Diagonally Implicit Runge Kutta method (ESDIRK). Explicit methods that complement TR-BDF2 as second order IMEX pairs have been introduced in \cite{giraldo:2013} and successfully employed in \cite{bonaventura:2018b}, \cite{garres:2021}.
Unconditionally strong stability preserving extensions of TR-BDF2 have been derived in \cite{bonaventura:2017}.
While the third order method that constitutes an embedded pair with TR-BDF2 is only conditionally stable, see the discussion in \cite{hosea:1996}, a first order embedded method is derived in \cite{kennedy:2016}, thus allowing for efficient time adaptation strategies. Finally, the analysis presented in \cite{bonaventura:2021} shows that
the method is optimal among second order methods for typical structural mechanics equations, thus making it an excellent candidate also for applications to fluid-structure interaction problems. While we do not pursue these developments in the present work, we would like to highlight these features as strong motivations for our specific choice of the time discretization method.

Following then the projection approach described in \cite{dellarocca:2018} and applying method \eqref{eq:trbdf2} to system \eqref{eq:ns_artcomp}, the momentum predictor equation for the first stage reads:
 
\begin{eqnarray}
\label{eq:first_stage}
&&\frac{\mathbf{u}^{n+\gamma,*} - \mathbf{u}^{n}}{\gamma\Delta t} - \frac{1}{2Re}\Delta\mathbf{u}^{n+\gamma,*} + \frac{1}{2}\dive\left(\mathbf{u}^{n+\gamma,*}\otimes\mathbf{u}^{n+\frac{\gamma}{2}}\right) = \nonumber \\
&&\frac{1}{2Re}\Delta\mathbf{u}^{n} - \frac{1}{2}\dive\left(\mathbf{u}^{n}\otimes\mathbf{u}^{n+\frac{\gamma}{2}}\right) - \nabla p^n\\
&&\mathbf{u}^{n+\gamma,*}\rvert_{\partial\Omega} = \mathbf{u}_D^{n+\gamma}. \nonumber
\end{eqnarray}

Notice that, in order to avoid solving a nonlinear system at each time step, an approximation is introduced in the nonlinear momentum advection term, so that \(\mathbf{u}^{n + \frac{\gamma}{2}}\) is defined by extrapolation as

$$\mathbf{u}^{n + \frac{\gamma}{2}} = \left(1 + \frac{\gamma}{2\left(1-\gamma\right)}\right)\mathbf{u}^{n} - \frac{\gamma}{2\left(1-\gamma\right)}\mathbf{u}^{n-1}.$$

Alternatively, \(\mathbf{u}^{n + \frac{\gamma}{2}}\) can be replaced by $\mathbf{u}^{n+\gamma,*} $ in the left hand side and by $ \mathbf{u}^{n} $ in the right hand side of \eqref{eq:first_stage}, respectively, and $ \mathbf{u}^{n+\gamma,*} $ can be determined by fixed point iteration. Numerical experiments show that this fully nonlinear formulation is necessary to achieve accurate results for larger Courant number values,
see the discussion in Section \ref{sec:tests}.
Following \cite{bell:1989}, we set then $ \delta p^{n+\gamma} = p^{n+\gamma} - p^{n} $ and impose

\begin{eqnarray}
\label{eq:discpress}
&&\frac{\mathbf{u}^{n+\gamma} - \mathbf{u}^{n+\gamma,*}}{\gamma\Delta t} =-\nabla \delta p^{n+\gamma}  \nonumber \\
&&\frac{1}{c^2}\frac{\delta p^{n+\gamma} }{\gamma\Delta t} + \dive\mathbf{u}^{n+\gamma}  =0.
\end{eqnarray}

Substituting the first equation into the second in \eqref{eq:discpress}, one obtains
the Helmholtz equation

\begin{equation}
\label{eq:helmholtz2} 
\frac{1}{c^2 \gamma^2\Delta t^2 }\delta p^{n+\gamma} - \Delta \delta p^{n+\gamma} =
 - \frac{1}{  \gamma\Delta t }\dive\mathbf{u}^{n+\gamma,*},
\end{equation}

which is solved with the boundary condition $ \nabla \delta p^{n+\gamma}\cdot\mathbf{n}\rvert_{\partial\Omega} = 0. $ Once this equation is solved, the final velocity update for the first stage \(\mathbf{u}^{n+\gamma} = \mathbf{u}^{n+\gamma,*} -  \gamma\Delta t\nabla \delta p^{n+\gamma}\) can be computed. Notice that the previous procedure is equivalent to introducing the intermediate update \(\mathbf{u}^{n+\gamma,**} = \mathbf{u}^{n+\gamma,*} + \gamma\Delta t\nabla  p^{n}\), solving

\begin{equation}
\label{eq:helmholtz2a}
\frac{1}{c^2}\frac{p^{n+\gamma}}{\gamma^2\Delta t^2} -\Delta p^{n+\gamma} = 
- \frac{1}{\gamma\Delta t} \dive\mathbf{u}^{n+\gamma,**}  + \frac{1}{c^2}\frac{p^{n }}{\gamma^2\Delta t^2}
\end{equation}

and then setting \(\mathbf{u}^{n+\gamma} = \mathbf{u}^{n+\gamma,**} - \gamma\Delta t\nabla  p^{n+\gamma}\). The second TR-BDF2 stage is performed in a similar manner. We first define the second momentum predictor:
 
\begin{eqnarray}
\label{eq:second_stage}
&&\frac{\mathbf{u}^{n+1,*} - \mathbf{u}^{n + \gamma}}{\left(1 - \gamma\right)\Delta t} - \frac{a_{33}}{Re}\Delta\mathbf{u}^{n+1,*} + a_{33}\dive\left(\mathbf{u}^{n+1,*}\otimes\mathbf{u}^{n+\frac{3}{2}\gamma}\right) =  \\
&&\frac{a_{32}}{Re}\Delta\mathbf{u}^{n + \gamma}  - a_{32}\dive\left(\mathbf{u}^{n + \gamma}\otimes\mathbf{u}^{n + \gamma}\right) + \frac{a_{31}}{Re}\Delta\mathbf{u}^{n} - a_{31}\dive\left(\mathbf{u}^{n}\otimes\mathbf{u}^{n}\right) - \nabla p^{n+\gamma} \nonumber \\
&&\mathbf{u}^{n+1,*}\rvert_{\partial\Omega^D} = \mathbf{u}_D^{n+1}, \nonumber
\end{eqnarray}

where one has

$$
a_{31} = \frac{1-\gamma}{2\left(2-\gamma\right)} \ \ \ 
a_{32} = \frac{1-\gamma}{2\left(2-\gamma\right)} \ \ \
a_{33} = \frac{1}{2-\gamma}.
$$

Again, in order to avoid solving a nonlinear system at each time step, an approximation is introduced in the nonlinear momentum advection term, so that $ \mathbf{u}^{n+\frac{3}{2}\gamma} $ is defined by extrapolation as 

$$\mathbf{u}^{n+\frac{3}{2}\gamma} = \left(1 + \frac{1+\gamma}{\gamma}\right)\mathbf{u}^{n+\gamma} - \frac{1-\gamma}{\gamma}\mathbf{u}^{n}.$$ 

Alternatively, $\mathbf{u}^{n+\frac{3}{2}\gamma}$ can be replaced by $\mathbf{u}^{n+1,*}, $ which can then be determined by fixed point iteration.
We set then $\delta p^{n+1} =p^{n+1} - p^{n+\gamma} $ and impose

\begin{eqnarray}
\label{eq:discpress2}
&&\frac{\mathbf{u}^{n+1} - \mathbf{u}^{n+1,*}}{(1-\gamma)\Delta t} = -\nabla \delta p^{n+1}  \nonumber \\
&&\frac{1}{c^2}\frac{\delta p^{n+1}}{(1-\gamma)\Delta t} + \dive\mathbf{u}^{n+1}  =0.
\end{eqnarray}

Substituting the first equation into the second in \eqref{eq:discpress2}, one obtains
the Helmholtz equation

\begin{equation}
\label{eq:helmholtz3} 
\frac{1}{c^2 (1-\gamma)^2\Delta t^2}\delta p^{n+1} - \Delta \delta p^{n+1} = - \frac{1}{(1-\gamma)\Delta t}\dive\mathbf{u}^{n+1,*},
\end{equation}

which is solved with the boundary condition $ \nabla \delta p^{n+1}\cdot\mathbf{n}\rvert_{\partial\Omega} = 0. $ Once this equation is solved, the final velocity update

$$ \mathbf{u}^{n+1} = \mathbf{u}^{n+1,*}-(1-\gamma)\Delta t\nabla \delta p^{n+1} $$ 

can be computed. Also for this second stage, notice that the procedure is equivalent to setting \(\mathbf{u}^{n+1,**} = \mathbf{u}^{n+1,*} + (1-\gamma)\Delta t\nabla  p^{n+\gamma}\), solving

\begin{equation}
\label{eq:helmholtz3a}
\frac{1}{c^2}\frac{ p^{n+1} }{(1-\gamma)^2\Delta t^2} -\Delta p^{n+1}= 
- \frac{1 }{(1-\gamma)\Delta t} \dive\mathbf{u}^{n+1,**}  + \frac{1}{c^2}\frac{p^{n+\gamma }}{(1-\gamma)^2\Delta t^2}
\end{equation}

and then setting \(\mathbf{u}^{n+1} = \mathbf{u}^{n+1,**} -  (1-\gamma)\Delta t\nabla  p^{n+1}\).

For the purposes of the comparisons that will be reported in Section \ref{sec:tests}, we also present two alternative and very popular second order projection methods, proposed respectively in \cite{bell:1989} and in \cite{guermond:1998}, which are based on the parent methods of TR-BDF2, i.e. the Crank-Nicolson (or Trapezoidal Rule) method and the BDF2 method, respectively. We briefly recall the formulation of these  schemes in the framework of the artificial compressibility formulation. 
The momentum predictor for the Bell-Colella-Glaz  \cite{bell:1989}  projection method reads as follows

\begin{eqnarray}
\label{eq:momentum_BCG}
&&\frac{\mathbf{u}^{n+1,*} - \mathbf{u}^{n}}{\Delta t} - \frac{1}{2Re}\Delta\mathbf{u}^{n+1,*} +
\left[\left(\mathbf{u}\cdot\nabla\right)\mathbf{u}\right]^{n + \frac{1}{2},*} = \nonumber \\
&&\frac{1}{2Re}\Delta\mathbf{u}^{n} - \nabla p^n\\
&&\mathbf{u}^{n+1,*}\rvert_{\partial\Omega} = \mathbf{u}_D^{n+1}. \nonumber
\end{eqnarray}

Notice that here we have set \(\mathbf{u}^{n + \frac{1}{2},*} = \frac{1}{2}\left(\mathbf{u}^{n+1,*} +  \mathbf{u}^{n }\right),\) so that the scheme
is fully nonlinear. On the other hand, setting \(\delta p^{n + 1} = p^{n + 1} - p^{n}\), we obtain the following Helmholtz equation for the projection stage

\begin{eqnarray}
\label{eq:helmholtz_BCG} 
&&\frac{1}{c^2 \Delta t^2 }\delta p^{n+1} - \Delta \delta p^{n+1} =
- \frac{1}{  \Delta t }\dive\mathbf{u}^{n+1,*} \\
&&\nabla \delta p^{n+1}\cdot\mathbf{n}\rvert_{\partial\Omega} = 0. \nonumber
\end{eqnarray}

Eventually, the velocity has to be updated with the gradient of the pressure increment:

\begin{equation}
\mathbf{u}^{n+1} = \mathbf{u}^{n+1,*} - \Delta t\nabla\delta p^{n+1}.
\end{equation}

It is apparent that this method is essentially based on the Crank-Nicolson time discretization approach. A method based on the BDF2 scheme has been presented instead by Guermond and Quartapelle in \cite{guermond:1998}. The momentum predictor reads as follows

\begin{eqnarray}
\label{eq:momentum_HI}
&&\frac{3\mathbf{u}^{n+1,*} - 4\mathbf{u}^{n} + \mathbf{u}^{n-1}}{\Delta t} - \frac{1}{Re}\Delta\mathbf{u}^{n+1,*} + \left(\mathbf{u}^{n,*}\cdot\nabla\right)\mathbf{u}^{n+1,*} \nonumber \\
&&+ \frac{1}{2}\left(\nabla\cdot\mathbf{u}^{n,*}\right)\mathbf{u}^{n+1,*} = -\nabla p^n\\
&&\mathbf{u}^{n+1,*}\rvert_{\partial\Omega} = \mathbf{u}_D^{n+1}. \nonumber
\end{eqnarray}

The Helmholtz equation for the projection stage is

\begin{eqnarray}
\label{eq:helmholtz_HI} 
&&\frac{1}{c^2 \Delta t^2 }p^{n+1} - \Delta p^{n+1} =
- \frac{1}{\Delta t}\dive\mathbf{u}^{n+1,**} + \frac{1}{c^2 \Delta t^2 }p^{n} \\
&&\nabla \delta p^{n+1}\cdot\mathbf{n}\rvert_{\partial\Omega} = 0, \nonumber
\end{eqnarray}

where \(\mathbf{u}^{n+1,**} = \mathbf{u}^{n+1,*} + \frac{2}{3}\Delta t\nabla p^{n}\). Eventually, the velocity is updated with the gradient of the computed pressure:

\begin{equation}
\mathbf{u}^{n+1} = \mathbf{u}^{n+1,**} - \frac{2}{3}\Delta t\nabla p^{n+1}.
\end{equation}

\section{The spatial discretization}
\label{sec:numeth} \indent

For the spatial discretization, we consider discontinuous finite element approximations, due to their great flexibility in performing mesh adaptation.
We consider a decomposition of the domain \(\Omega\) into a family of hexahedra \(\mathcal{T}_h\) (quadrilaterals in the two-dimensional case) and denote each element by \(K\). The skeleton \(\mathcal{E}\) denotes the set of all element faces and \(\mathcal{E} = \mathcal{E}^{I} \cup \mathcal{E}^{B}\), where \(\mathcal{E}^{I}\) is the subset of interior faces and \(\mathcal{E}^{B}\) is the subset of boundary faces. We also introduce the following finite element spaces

\[Q_k = \left\{v \in L^2(\Omega) : v\rvert_K \in \mathbb{Q}_k \quad \forall K \in \mathcal{T}_h\right\}\] 

and

\[\mathbf{Q}_k = \left[Q_k\right]^d,\]

where \(\mathbb{Q}_k\) is the space of polynomials of degree \(k\) in each coordinate direction. Considering the well-posedness analyses in \cite{schotzau:2003, toselli:2002}, the finite element spaces that will be used for the discretization of velocity and pressure are \(\mathbf{V}_h = \mathbf{Q}_k\) and \(Q_h = Q_{k-1} \cap L^2_{0}(\Omega)\), respectively,  where \(k \ge 2\) and \(L^2_{0}(\Omega) = \left\{v \in L^2(\Omega) : \int_{\Omega} vd\Omega = 0 \right\}\). Notice that, while for the sake of coherence with the time discretization and of comparison with second order finite volume methods we will mostly consider the case $ k = 2 $ in the following, the formulation we present is completely general and also the implementation validated in Section \ref{sec:tests} supports arbitrary values of $ k. $ Furthermore, notice that
the above choice for the finite element spaces corresponds to that implemented in the \textit{deal.II} library, which will be employed for the numerical computation. The proposed approach can in principle also be applied to tetrahedral meshes and \(P\)-spaces.
Suitable jump and average operators can then be defined as customary for finite element discretizations, see e.g. \cite{brezzi:2002}. 
A face \(\Gamma \in \mathcal{E}^{I}\) shares two elements that we denote by \(K^{+}\) with outward unit normal \(\mathbf{n}^{+}\) and \(K^{-}\) with outward unit normal \(\mathbf{n}^{-}\), whereas for a face \(\Gamma \in \mathcal{E}^{B}\) we denote by \(\mathbf{n}\) the outward unit normal.
For a scalar function \(\varphi\) the jump is defined as

\[\left[\left[\varphi\right]\right] = \varphi^{+}\mathbf{n}^{+} + \varphi^{-}\mathbf{n}^{-} \quad \text{if }\Gamma \in \mathcal{E}^{I} \qquad \left[\left[\varphi\right]\right] = \varphi\mathbf{n} \quad \text{if }\Gamma \in \mathcal{E}^{B}.\]

The average is defined as

\[\left\{\left\{\varphi\right\}\right\} = \frac{1}{2}\left(\varphi^{+} + \varphi^{-}\right) \quad \text{if }\Gamma \in \mathcal{E}^{I} \qquad \left\{\left\{\varphi\right\}\right\} = \varphi \quad \text{if }\Gamma \in \mathcal{E}^{B}.\]

Similar definitions apply for a vector function \(\boldsymbol{\varphi}\):

\begin{align*}
\left[\left[\boldsymbol{\varphi}\right]\right] &= \boldsymbol{\varphi}^{+}\cdot\mathbf{n}^{+} 
+\boldsymbol{\varphi}^{-}\cdot\mathbf{n}^{-} \quad \text{if }\Gamma \in \mathcal{E}^{I} \qquad 
\left[\left[\boldsymbol{\varphi}\right]\right] = \boldsymbol{\varphi}\cdot\mathbf{n} \quad \text{if }\Gamma \in \mathcal{E}^{B} \\
\left\{\left\{\boldsymbol{\varphi}\right\}\right\} &= \frac{1}{2}\left(\boldsymbol{\varphi}^{+} + \boldsymbol{\varphi}^{-}\right) \quad \text{if }\Gamma \in \mathcal{E}^{I} \quad\qquad \left\{\left\{\boldsymbol{\varphi}\right\}\right\} = \boldsymbol{\varphi} \quad \text{if }\Gamma \in \mathcal{E}^{B}.
\end{align*}

For vector  functions, it is also useful  to define a tensor jump as:

\[\left<\left<\boldsymbol{\varphi}\right>\right> = \boldsymbol{\varphi}^{+}\otimes\mathbf{n}^{+} 
+ \boldsymbol{\varphi}^{-}\otimes\mathbf{n}^{-} \quad \text{if }\Gamma \in \mathcal{E}^{I} 
\qquad \left<\left<\boldsymbol{\varphi}\right>\right> = \boldsymbol{\varphi}\otimes\mathbf{n} \quad \text{if }\Gamma \in \mathcal{E}^{B}.\]

Given these definitions, the weak formulation of the momentum predictor equation for the first stage is obtained multiplying equation \eqref{eq:first_stage} by a test function \(\mathbf{v} \in \mathbf{V}_h\), integrating over \(K \in \mathcal{T}_h\) and applying Green's theorem. 
To impose the boundary conditions, we set \(\left(\mathbf{u}^{n+\gamma,*}\right)^{-} = -\left(\mathbf{u}^{n+\gamma,*}\right)^{+} + 2\mathbf{u}_D^{n+\gamma}\) with \(\left[\nabla\left(\mathbf{u}^{n+\gamma,*}\right)^{+}\right]\cdot\mathbf{n} = \left[\nabla\left(\mathbf{u}^{n+\gamma,*}\right)^{-}\right]\cdot\mathbf{n}\). 
 
We now treat separately the discretization of the diffusion and advection contributions, respectively. The approximation of the diffusion term is based on the symmetric interior penalty method (SIP) \cite{arnold:1982}. We denote the scalar product between two second-order tensors by

\[\mathbf{A}:\mathbf{B} = \sum_{i,j}A_{ij}B_{ij}. \]

Following \cite{fehn:2019}, we set for each face \(\Gamma\) of a cell \(K\) 

\begin{equation}
\label{eq:penconst1}
    \sigma^{\mathbf{u}}_{\Gamma,K} = \left(k + 1 \right)^2\frac{\text{diam}(\Gamma)}{\text{diam}(K)}
\end{equation}

and we define the penalization constant for the SIP method as 

$$C_u = \frac{1}{2}\left(\sigma^{\mathbf{u}}_{\Gamma,K^{+}} + \sigma^{\mathbf{u}}_{\Gamma,K^{-}}\right)$$

if \(\Gamma \in \mathcal{E}^{I}\) and \(C_u = \sigma^{\mathbf{u}}_{\Gamma,K}\) otherwise. 
Taking into account boundary conditions as previously discussed and summing over all \(K \in \mathcal{T}_h\), we can define the following bilinear form:
 
\begin{eqnarray}
\label{eq:bilin_diff}
a_{\mathbf{u}}^{(1)}(\mathbf{u},\mathbf{v}) &=& \frac{1}{2Re}\sum_{K\in\mathcal{T}_h}\int_K\nabla\mathbf{u}:\nabla\mathbf{v}d\Omega 
- \frac{1}{2Re}\sum_{\Gamma \in \mathcal{E}^{I}}\int_\Gamma\left\{\left\{\nabla\mathbf{u}\right\}\right\}:\left<\left<\mathbf{v}\right>\right>d\Sigma \nonumber \\
&-& \frac{1}{2Re}\sum_{\Gamma \in \mathcal{E}^{B}}\int_\Gamma\left(\nabla\mathbf{u}\right)\mathbf{n}\cdot\mathbf{v}d\Sigma -\frac{1}{2Re}\sum_{\Gamma \in \mathcal{E}^{I}}\int_\Gamma \left<\left<\mathbf{u}\right>\right>:\left\{\left\{\nabla\mathbf{v}\right\}\right\}d\Sigma\nonumber\\
&-& \frac{1}{2Re}\sum_{\Gamma \in \mathcal{E}^{B}}\int_\Gamma \left(\mathbf{u}\otimes\mathbf{n}\right):\nabla\mathbf{v}d\Sigma
+\frac{1}{2Re}\sum_{\Gamma \in \mathcal{E}^I}\int_\Gamma C_u\left<\left<\mathbf{u}\right>\right>:\left<\left<\mathbf{v}\right>\right>d\Sigma \nonumber\\
&+& \frac{1}{2Re}\sum_{\Gamma \in \mathcal{E}^B}\int_\Gamma 2C_u\left(\mathbf{u}\cdot\mathbf{v}\right)d\Sigma. 
\end{eqnarray}

The approximation of the advection term employs the widely used local Lax-Friedrichs (LF) flux, see e.g. \cite{giraldo:2020}. Setting 
$$\lambda = \max\left(\left|\left(\mathbf{u}^{n+\frac{\gamma}{2}} \right)^{+}\cdot\mathbf{n}\right|,\left|\left(\mathbf{u}^{n+\frac{\gamma}{2}} \right)^{-}\cdot\mathbf{n}\right|\right) $$

with \(\mathbf{n} = \mathbf{n}^{\pm}\) and taking into account boundary conditions,
we define the trilinear form

\begin{eqnarray}
\label{eq:trilinear}
c^{(1)}(\mathbf{u}^{n+\frac{\gamma}{2}}, \mathbf{u},\mathbf{v}) &= &-\frac{1}{2}\sum_{K\in\mathcal{T}_h}\int_K\left(\mathbf{u}\otimes\mathbf{u}^{n+\frac{\gamma}{2}}\right):\nabla\mathbf{v}d\Omega + \frac{1}{2}\sum_{\Gamma \in \mathcal{E}^{I}}\int_{\Gamma}\left(\left\{\left\{\mathbf{u}\otimes\mathbf{u}^{n+\frac{\gamma}{2}}\right\}\right\}\right):\left<\left<\mathbf{v}\right>\right>d\Sigma \nonumber \\
&+& \frac{1}{2}\sum_{\Gamma \in \mathcal{E}^{I}}\int_{\Gamma}\frac{\lambda}{2}\left<\left<\mathbf{u}\right>\right>:\left<\left<\mathbf{v}\right>\right>d\Sigma + \frac{1}{2}\sum_{\Gamma \in \mathcal{E}^{B}}\int_\Gamma\lambda\left(\mathbf{u}\cdot\mathbf{v}\right)d\Sigma.
\end{eqnarray}

Finally, we also define the functional  
 
\begin{eqnarray}
\label{eq:rhs1}
F_{\mathbf{u}}^{(1)}(\mathbf{v})^{n+\gamma} = &-&\frac{1}{2Re}\sum_{K\in\mathcal{T}_h}\int_K\nabla\mathbf{u}^{n}:\nabla\mathbf{v}d\Omega  
+ \frac{1}{2Re}\sum_{\Gamma \in \mathcal{E}}\int_\Gamma\left\{\left\{\nabla\mathbf{u}^{n}\right\}\right\}:\left<\left<\mathbf{v}\right>\right>d\Sigma  \nonumber \\
&+& \frac{1}{2}\sum_{K\in\mathcal{T}_h}\int_K\left(\mathbf{u}^{n}\otimes\mathbf{u}^{n+\frac{\gamma}{2}}\right):\nabla\mathbf{v}d\Omega -\frac{1}{2}\sum_{\Gamma \in \mathcal{E}}\int_{\Gamma}\left(\left\{\left\{\mathbf{u}^{n}\otimes\mathbf{u}^{n+\frac{\gamma}{2}}\right\}\right\}\right):\left<\left<\mathbf{v}\right>\right>d\Sigma \nonumber \\
&+&\sum_{K\in\mathcal{T}_h}\int_K p^n\dive\mathbf{v}d\Omega - \sum_{\Gamma \in \mathcal{E}}\int_{\Gamma}\left\{\left\{p^n\right\}\right\}\left[\left[\mathbf{v}\right]\right]d\Sigma - \frac{1}{2Re}\sum_{\Gamma \in \mathcal{E}^{B}}\int_\Gamma\left(\mathbf{u}_D^{n+\gamma}\otimes\mathbf{n}\right):\nabla\mathbf{v}d\Sigma \nonumber \\ 
&+& \frac{1}{2Re}\sum_{\Gamma \in \mathcal{E}^{B}}\int_\Gamma2C_u\left(\mathbf{u}_D^{n+\gamma}\cdot\mathbf{v}\right)d\Sigma - \frac{1}{2}\sum_{\Gamma \in \mathcal{E}^{B}}\int_\Gamma \left(\mathbf{u}_D^{n+\gamma}\otimes\mathbf{u}^{n+\frac{\gamma}{2}}\right)\mathbf{n}\cdot\mathbf{v}d\Sigma \nonumber \\
&+& \frac{1}{2}\sum_{\Gamma \in \mathcal{E}^{B}}\int_\Gamma\lambda\left(\mathbf{u}_D^{n+\gamma}\cdot\mathbf{v}\right)d\Sigma,
\end{eqnarray} 

which also includes the terms representing the weak form of Dirichlet boundary conditions. It is worth to point out that in the right-hand side no penalization terms have been introduced for the variables computed at previous time-steps. Moreover, for the sake of clarity, the face integrals related to the quantities at previous time-steps are reported on the whole skeleton \(\mathcal{E}\), without distinguishing between interior and boundary faces.

The complete weak formulation of the first stage velocity update reads then as follows:
\textit{ given \(\mathbf{u}^{n+\frac{\gamma}{2}}, \mathbf{u}^{n} \in \mathbf{V}_h\) and \(p^n \in Q_h\), find \(\mathbf{u}^{n+\gamma,*} \in \mathbf{V}_h\) such that}: 

\begin{eqnarray}
\label{eq:weakform_u_1}
&&\sum_{K \in \mathcal{T}_h}\int_K \frac{1}{\gamma\Delta t}\mathbf{u}^{n+\gamma,*}\cdot\mathbf{v}d\Omega + a^{(1)}_{\mathbf{u}}(\mathbf{u}^{n+\gamma,*},\mathbf{v})  + c^{(1)}(\mathbf{u}^{n+\frac{\gamma}{2}}, \mathbf{u}^{n+\gamma,*},\mathbf{v})  \nonumber \\
&&=\sum_{K \in \mathcal{T}_h}\int_K \frac{1}{\gamma\Delta t}\mathbf{u}^{n}\cdot\mathbf{v}d\Omega + F^{(1)}_{\mathbf{u}}(\mathbf{v})^{n+\gamma} \quad \forall\mathbf{v}\in\mathbf{V}_h.
\end{eqnarray}
\pagebreak

For the projection steps defined by equation \eqref{eq:helmholtz2a} we apply again the SIP method. In order to impose homogeneous Neumann boundary conditions we prescribe \(\left[\nabla\left(p^{n+\gamma}\right)^{-}\right]\mathbf{n} = -\left[\nabla\left(p^{n+\gamma}\right)^{+}\right]\mathbf{n}\): for this reason, no contribution from boundary faces arises. We then multiply by a test function \(q \in Q_h\), we apply Green's theorem and we define:
 
\begin{eqnarray}
\label{eq:weakform_p_1}
a_{p}(p,q)  &=& \sum_{K \in \mathcal{T}_h}\int_K\nabla p\cdot\nabla qd\Omega  
- \sum_{\Gamma \in \mathcal{E}^{I}}\int_\Gamma \left\{\left\{\nabla p^{n+\gamma}\right\}\right\}\cdot\left[\left[q\right]\right]d\Sigma \nonumber \\
&-&\sum_{\Gamma \in \mathcal{E}^{I}}\int_\Gamma \left[\left[p\right]\right]\cdot\left\{\left\{\nabla q\right\}\right\}d\Sigma + \sum_{\Gamma \in \mathcal{E}^{I}}\int_\Gamma C_p\left[\left[p\right]\right]\cdot\left[\left[q\right]\right]d\Sigma \\
F^{(1)}_{p}(q)^{n+\gamma} &=& \sum_{K \in \mathcal{T}_h}\int_K\frac{1}{\gamma\Delta t}\mathbf{u}^{n+\gamma,**}\cdot\nabla qd\Omega - \sum_{\Gamma \in \mathcal{E}}\int_\Gamma\frac{1}{\gamma\Delta t}\left\{\left\{\mathbf{u}^{n+\gamma,**}\right\}\right\}\cdot\left[\left[q\right]\right]d\Sigma
\end{eqnarray} 

and again we set

\begin{equation}
\sigma^{p}_{\Gamma,K} = k^2\frac{\text{diam}(\Gamma)}{\text{diam}(K)},  
\end{equation}

while, if \(\Gamma \in \mathcal{E}^{I}\), we set \(C_p = \frac{1}{2}\left(\sigma^{p}_{\Gamma,K^{+}} + \sigma^{p}_{\Gamma,K^{-}}\right)\), otherwise \(C_p = \sigma^{p}_{\Gamma,K}\). The weak formulation of equation \eqref{eq:helmholtz2a} reads then: \textit{given \(p^n \in Q_h\), find \(p^{n+\gamma} \in Q_h\) such that}

\begin{equation}
\label{helmholtz_1}
\sum_{K \in \mathcal{T}_h}\int_K \frac{1}{c^2\gamma^2\Delta t^2}p^{n+\gamma}q d\Omega + a_p(p^{n+\gamma},q) = \sum_{K \in \mathcal{T}_h}\int_K \frac{1}{c^2\gamma^2\Delta t^2}p^{n} q d\Omega + F^{(1)}_p(q)^{n+\gamma} \qquad \forall q \in Q_h.
\end{equation}

The second stage can be described in a similar manner. We start defining the bilinear forms for the second momentum predictor as

\begin{align}
a_{\mathbf{u}}^{(2)}(\mathbf{u},\mathbf{v}) &= \frac{a_{33}}{Re}\sum_{K\in\mathcal{T}_h}\int_K\nabla\mathbf{u}:\nabla\mathbf{v}d\Omega - \frac{a_{33}}{Re}\sum_{\Gamma \in \mathcal{E}^{I}}\int_\Gamma\left\{\left\{\nabla\mathbf{u}\right\}\right\}:\left<\left<\mathbf{v}\right>\right>d\Sigma \nonumber \\
&- \frac{a_{33}}{Re}\sum_{\Gamma \in \mathcal{E}^{B}}\int_\Gamma\left(\nabla\mathbf{u}\right)\mathbf{n}\cdot\mathbf{v}d\Sigma 
- \frac{a_{33}}{Re}\sum_{\Gamma \in \mathcal{E}^{I}}\int_\Gamma \left<\left<\mathbf{u}\right>\right>:\left\{\left\{\nabla\mathbf{v}\right\}\right\}d\Sigma \nonumber \\
&- \frac{a_{33}}{Re}\sum_{\Gamma \in \mathcal{E}^{B}}\int_\Gamma \left(\mathbf{u}\otimes\mathbf{n}\right):\nabla\mathbf{v}d\Sigma  
+ \frac{a_{33}}{Re}\sum_{\Gamma \in \mathcal{E}^I}\int_\Gamma C_u\left<\left<\mathbf{u}\right>\right>:\left<\left<\mathbf{v}\right>\right>d\Sigma \nonumber \\
&+ \frac{a_{33}}{Re}\sum_{\Gamma \in \mathcal{E}^B}\int_\Gamma 2C_u\left(\mathbf{u}\cdot\mathbf{v}\right)d\Sigma 
\end{align}

\begin{align}
c^{(2)}(\mathbf{u}^{n+\frac{3}{2}\gamma}, \mathbf{u},\mathbf{v}) = &-a_{33}\sum_{K\in\mathcal{T}_h}\int_K\left(\mathbf{u}\otimes\mathbf{u}^{n+\frac{3}{2}\gamma}\right):\nabla\mathbf{v}d\Omega + a_{33}\sum_{\Gamma \in \mathcal{E}^{I}}\int_{\Gamma}\left(\left\{\left\{\mathbf{u}\otimes\mathbf{u}^{n+\frac{3}{2}\gamma}\right\}\right\}\right):\left<\left<\mathbf{v}\right>\right>d\Sigma \nonumber \\
&+ a_{33}\sum_{\Gamma \in \mathcal{E}^{I}}\int_{\Gamma}\frac{\lambda}{2}\left<\left<\mathbf{u}\right>\right>:\left<\left<\mathbf{v}\right>\right>d\Sigma 
+ a_{33}\sum_{\Gamma \in \mathcal{E}^{B}}\int_\Gamma\lambda\left(\mathbf{u}\cdot\mathbf{v}\right)d\Sigma,
\end{align}

where \(\lambda = \max\left(\left|\left(\mathbf{u}^{n+\frac{3}{2}\gamma} \right)^{+}\cdot\mathbf{n}\right|,\left|\left(\mathbf{u}^{n+\frac{3}{2}\gamma} \right)^{-}\cdot\mathbf{n}\right|\right)\) with \(\mathbf{n} = \mathbf{n}^{\pm}\).

We also define the linear functional:

\begin{align}
F^{(2)}_{\mathbf{u}}(\mathbf{v})^{n+1} = &-\frac{a_{32}}{Re}\sum_{K\in\mathcal{T}_h}\int_K\nabla\mathbf{u}^{n + \gamma}:\nabla\mathbf{v}d\Omega + \frac{a_{32}}{Re}\sum_{\Gamma \in \mathcal{E}}\int_\Gamma\left\{\left\{\nabla\mathbf{u}^{n + \gamma}\right\}\right\}:\left<\left<\mathbf{v}\right>\right>d\Sigma \nonumber \\
&-\frac{a_{31}}{Re}\sum_{K\in\mathcal{T}_h}\int_K\nabla\mathbf{u}^{n}:\nabla\mathbf{v}d\Omega + \frac{a_{31}}{Re}\sum_{\Gamma \in \mathcal{E}}\int_\Gamma\left\{\left\{\nabla\mathbf{u}^{n}\right\}\right\}:\left<\left<\mathbf{v}\right>\right>d\Sigma \nonumber \\
&+ a_{32}\sum_{K\in\mathcal{T}_h}\int_K\left(\mathbf{u}^{n+\gamma}\otimes\mathbf{u}^{n+\gamma}\right):\nabla\mathbf{v}d\Omega - a_{32}\sum_{\Gamma \in \mathcal{E}}\int_{\Gamma}\left(\left\{\left\{\mathbf{u}^{n+\gamma}\otimes\mathbf{u}^{n+\gamma}\right\}\right\}\right):\left<\left<\mathbf{v}\right>\right>d\Sigma \nonumber \\
&+ a_{31}\sum_{K\in\mathcal{T}_h}\int_K\left(\mathbf{u}^{n}\otimes\mathbf{u}^{n}\right):\nabla\mathbf{v}d\Omega - a_{31}\sum_{\Gamma \in \mathcal{E}}\int_{\Gamma}\left(\left\{\left\{\mathbf{u}^{n}\otimes\mathbf{u}^{n}\right\}\right\}\right):\left<\left<\mathbf{v}\right>\right>d\Sigma \nonumber \\
&+\sum_{K\in\mathcal{T}_h}\int_K p^{n+\gamma}\dive\mathbf{v}d\Omega - \sum_{\Gamma \in \mathcal{E}}\int_{\Gamma}\left\{\left\{p^{n+\gamma}\right\}\right\}\left[\left[\mathbf{v}\right]\right]d\Sigma \nonumber \\
&-\frac{a_{33}}{Re}\sum_{\Gamma \in \mathcal{E}^{B}}\int_\Gamma\left(\mathbf{u}_D^{n+1}\otimes\mathbf{n}\right):\nabla\mathbf{v}d\Sigma + \frac{a_{33}}{Re}\sum_{\Gamma \in \mathcal{E}^{B}}\int_\Gamma2C_u\left(\mathbf{u}_D^{n+1}\cdot\mathbf{v}\right)d\Sigma
\nonumber \\
&- a_{33}\sum_{\Gamma \in \mathcal{E}^{B}}\int_\Gamma \left(\mathbf{u}_D^{n+1}\otimes\mathbf{u}^{n+\frac{3}{2}\gamma}\right)\mathbf{n}\cdot\mathbf{v}d\Sigma 
+ a_{33}\sum_{\Gamma \in \mathcal{E}^{B}}\int_\Gamma\lambda\left(\mathbf{u}_D^{n+1}\cdot\mathbf{v}\right)d\Sigma.
\end{align}

Finally, the weak formulation for the equation \eqref{eq:second_stage} reads as follows:
\textit{given \(\mathbf{u}^{n+\frac{3}{2}\gamma}, \mathbf{u}^{n + \gamma} \in \mathbf{V}_h\) and \(p^{n + \gamma} \in Q_h\), find \(\mathbf{u}^{n+\gamma,*} \in \mathbf{V}_h\) such that:}

\begin{eqnarray}
\label{eq:weakform_u_2}
&&\sum_{K \in \mathcal{T}_h}\int_K \frac{1}{\left(1 - \gamma\right)\Delta t}\mathbf{u}^{n+1,*}\cdot\mathbf{v}d\Omega + a^{(2)}_{\mathbf{u}}(\mathbf{u}^{n+1,*},\mathbf{v})  + c^{(2)}(\mathbf{u}^{n+\frac{3}{2}\gamma}, \mathbf{u}^{n+1,*},\mathbf{v})  \nonumber \\
&&=\sum_{K \in \mathcal{T}_h}\int_K \frac{1}{\left(1 - \gamma\right)\Delta t}\mathbf{u}^{n + \gamma}\cdot\mathbf{v}d\Omega + F^{(2)}_{\mathbf{u}}(\mathbf{v})^{n+1} \quad \forall\mathbf{v}\in\mathbf{V}_h.
\end{eqnarray}

We can then immediately define the functional associated to the second projection step as

\begin{align}
F^{(2)}_p(q)^{n+1} &= \sum_{K \in \mathcal{T}_h}\int_K\frac{1}{\left(1 - \gamma\right)\Delta t}\mathbf{u}^{n+1,**}\cdot\nabla qd\Omega 
- \sum_{\Gamma \in \mathcal{E}}\int_\Gamma\frac{1}{\left(1-\gamma\right)\Delta t}\left\{\left\{\mathbf{u}^{n+1,**}\right\}\right\}\cdot\left[\left[q\right]\right]d\Sigma.
\end{align}

Therefore, the weak formulation for \eqref{eq:helmholtz3a} reads as follows: 
\textit{given \(p^n \in Q_h\),  find \(p^{n+1} \in Q_h\) such that:}

\begin{align}
\label{helmholtz_2}
\sum_{K \in \mathcal{T}_h}\int_K \frac{1}{c^2(1-\gamma)^2\Delta t^2}p^{n+1}qd\Omega + a(p^{n+1},q)^{n+1} = \sum_{K \in \mathcal{T}_h}\int_K \frac{1}{c^2(1-\gamma)^2\Delta t^2}p^{n+\gamma}qd\Omega + F^{(2)}_p(q)^{n+1} \qquad \forall q \in Q_h.    
\end{align}

We now derive the fully discrete algebraic expressions corresponding to each of the two stages. We denote  by \(\boldsymbol{\varphi}_i(\mathbf{x})\) the basis functions for the space \(\mathbf{V}_h\) and by \(\psi_i(\mathbf{x})\) the basis functions for the space \(Q_h\), respectively, so that the discrete approximations of \(\mathbf{u}\) and \(p\) read as follows

\begin{equation*}
\mathbf{u}\approx \mathbf{u}_h = \sum_{j = 1}^{\text{dim}(\mathbf{V}_h)}u_j(t)\boldsymbol{\varphi}_j(\mathbf{x}) \qquad p \approx p_h = \sum_{j = 1}^{\text{dim}(Q_h)}p_j(t)\psi_j(\mathbf{x}).
\end{equation*}

For the first stage, we take \(\mathbf{v} = \boldsymbol{\varphi}_i\), \(i = 1,\dots,\text{dim}(\mathbf{V}_h)\) and we exploit the representation introduced above to obtain the matrices

\begin{eqnarray}
\label{eq:mass_matrix}
M_{ij} &= &\sum_{K\in\mathcal{T}_h}\int_K\boldsymbol{\varphi}_j\cdot\boldsymbol{\varphi}_id\Omega\\
A_{ij}^{n+\gamma} &=& a^{(1)}_{\mathbf{u}}\left(\boldsymbol{\varphi}_j, \boldsymbol{\varphi}_j\right) \\
C_{ij}\left(\mathbf{u}^{n + \frac{\gamma}{2}}\right) &=& c^{(1)}\left(\mathbf{u}^{n + \frac{\gamma}{2}}, \boldsymbol{\varphi}_j, \boldsymbol{\varphi}_i\right)
\end{eqnarray}

After computing the integrals in the previous formulae by appropriate quadrature rules,
one obtains the algebraic system

\begin{equation}
\left(\frac{1}{\gamma\Delta t}\mathbf{M} + \mathbf{A}^{n+\gamma} + \mathbf{C}\left(\mathbf{u}^{n+\frac{\gamma}{2}}\right)\right)\mathbf{U}_{h}^{n+\gamma,*} = \frac{1}{\gamma\Delta t}\mathbf{M}\mathbf{U}_{h}^{n} + \mathbf{F}_{\mathbf{u}}^{n+\gamma},
\end{equation}

where $\mathbf{U}_{h}$ denotes the vector of the discrete degrees of freedom associated to the velocity field and $\mathbf{F}_{\mathbf{u}}^{n+\gamma}$ is the vector obtained evaluating \(F^{(1)}_{\mathbf{u}}(\boldsymbol{\varphi}_i)^{n+\gamma}, i = 1,\dots,\text{dim}(\mathbf{V}_h).\)
The same procedure can be applied for the projection step, obtaining the matrices

\begin{eqnarray}
M^p_{{ij}}  &=& \sum_{K\in\mathcal{T}_h}\int_K\psi_j\psi_id\Omega \\ \label{eq:pressure_mass_matrix} 
K_{ij}  &=& a_p(\psi_j, \psi_i). \label{eq:weakform_p_1_algebraic_entries}
\end{eqnarray}
 
After computing the integrals in the previous formulae by appropriate quadrature rules,
one obtains the algebraic counterpart of \eqref{helmholtz_1} 
\begin{equation}
\left(\frac{1}{c^2\gamma^2\Delta t^2}\mathbf{M}^p + \mathbf{K} \right)\mathbf{P}_{h}^{n+\gamma} = \frac{1}{c^2\gamma^2\Delta t^2}\mathbf{M}^p\mathbf{P}_{h}^{n} + \mathbf{F}_{p}^{n+\gamma}
\end{equation}
where again $\mathbf{P}_{h} $ denotes the vector of the discrete degrees of freedom associated to pressure and $\mathbf{F}_{p}^{n+\gamma}$ is the vector obtained evaluating \(F^{(1)}_{p}(\boldsymbol{\psi}_i)^{n+\gamma}, i = 1,\dots,\text{dim}(Q_h).\)
For the second stage, we proceed in a similar manner; for the momentum predictor \eqref{eq:weakform_u_2} we obtain

\begin{equation}
\left(\frac{1}{\left(1 - \gamma\right)\Delta t}\mathbf{M} + \mathbf{A}^{n+1} + \mathbf{C}\left(\mathbf{u}^{n+\frac{3}{2}\gamma}\right)\right)\mathbf{U}_{h}^{n+1,*} = \frac{1}{\left(1 - \gamma\right)\Delta t}\mathbf{M}\mathbf{U}_{h}^{n + \gamma} + \mathbf{F}_{\mathbf{u}}^{n+1},
\end{equation}

where we set

\begin{equation}
A_{ij}^{n+1} = a^{(2)}_{\mathbf{u}}\left(\boldsymbol{\varphi}_j,\boldsymbol{\varphi}_i\right) \ \ \ \
C_{ij}\left(\mathbf{u}^{n+\frac{3}{2}\gamma}\right) = c^{(2)}\left(\mathbf{u}^{n+\frac{3}{2}\gamma}, \boldsymbol{\varphi}_j,\boldsymbol{\varphi}_i\right)
\end{equation}

and \(\mathbf{F}_{\mathbf{u}}^{n+1}\) is the vector obtained evaluating \(F^{(2)}_{\mathbf{u}}(\boldsymbol{\varphi}_i)^{n+1}, i = 1...\text{dim}(\mathbf{V}_h).\) Eventually, as algebraic counterpart of \eqref{helmholtz_2} we obtain

\begin{equation}
\left(\frac{1}{c^2\left(1-\gamma\right)^2\Delta t^2}\mathbf{M}^p + \mathbf{K} \right)\mathbf{P}_{h}^{n+1} = \frac{1}{c^2\left(1-\gamma\right)^2\Delta t^2}\mathbf{M}^p\mathbf{P}_{h}^{n+\gamma} + \mathbf{F}_{p}^{n+1},
\end{equation}

where again $\mathbf{F}_{p}^{n+1}$ is the vector obtained evaluating \(F^{(2)}_{p}(\boldsymbol{\psi}_i)^{n+1}, i = 1...\text{dim}(Q_h).\) Notice that, in the evaluation of \(\mathbf{F}_{p}^{n+\gamma}\) and \(\mathbf{F}_{p}^{n+1}\), there is also a preliminary stage which is the projection of \(\nabla p^{n}\) and \(\nabla p^{n+\gamma}\) into \(\mathbf{V}_h\) to compute \(\mathbf{u}^{n+\gamma,**}\) and \(\mathbf{u}^{n+1,**}\), respectively. In particular, we define the projection matrix \(\mathbf{P}\)  

\begin{equation}
P_{ij} = \sum_{K \in \mathcal{T}_h}\int_K\nabla\psi_j\cdot\boldsymbol{\varphi}_id\Omega
\end{equation} 

and we solve the linear systems \(\mathbf{M}\tilde{\mathbf{u}}^{n+\gamma,**} = \mathbf{P}\mathbf{p}^{n}\) for the first stage and \(\mathbf{M}\tilde{\mathbf{u}}^{n+1,**} = \mathbf{P}\mathbf{p}^{n + \gamma}\), where \(\tilde{\mathbf{u}}^{n+\gamma,**}\) and \(\tilde{\mathbf{u}}^{n+1,**}\) denote the two required projections. The same procedure has to applied also in the final update of the velocity; in particular, for the first stage we set 

\begin{equation}
\mathbf{u}^{n+\gamma} = \mathbf{u}^{n+\gamma,*} - \gamma\Delta t\left(\tilde{\mathbf{u}}^{n+1,**} - \tilde{\mathbf{u}}^{n+\gamma,**}\right),
\end{equation}

while for the second stage we solve \(\mathbf{M}\tilde{\mathbf{u}}^{n+1} = \mathbf{P}\mathbf{p}^{n + 1}\) and then we compute

\begin{equation}
\mathbf{u}^{n+1} = \mathbf{u}^{n+1,*} - \left(1 - \gamma\right)\Delta t\left(\tilde{\mathbf{u}}^{n+1} - \tilde{\mathbf{u}}^{n+1,**}\right).
\end{equation}

\section{Numerical experiments}
\label{sec:tests} \indent

The numerical method outlined in the previous Sections has been validated in a number of relevant benchmarks. Notice that, following e.g. \cite{tumolo:2015, tumolo:2013}, we set $\cal H =\min\{\mathrm{diam}(K) |  K \in \mathcal{T}_{h}\} $ and we define the stability parameters:

\begin{equation}
\label{eq:stabpar}
    C = kU \Delta t/\cal H, \ \ \ \ \    \mu = k^2 \Delta t/(Re\cal H^2),
\end{equation}

where \(U\) is the magnitude of a characteristic velocity and $ \mu $ defines the typical stability parameter in the discretization of parabolic terms. We also recall here that \(k\) is the polynomial degree of the finite element space chosen for the discretization of the velocity. As stated before, the proposed method has been implemented using the numerical library \textit{deal.II}, which is based on a matrix-free approach \cite{bangerth:2007}. The 9.2.0  version of \textit{deal.II} was employed  and most simulations have been run in parallel with MPI.
No global sparse matrix is built and only the action of the linear operators defined in Section \ref{sec:numeth} on a vector is actually implemented. Another feature of the library employed during the numerical simulations is the mesh adaptation capability, as we will see in the presentation of the results. In the following tests, unless differently stated, we take $ c=10^3\  m/s, $ which is the order of magnitude of the speed of sound in water. Moreover, the preconditioned conjugate gradient method implemented in the function \textit{SolverCG} of the \textit{deal.II} library was employed to solve the Helmholtz equations, while the GMRES solver for the momentum equations is implemented in the function \textit{SolverGMRES} of the same library. A Jacobi preconditioner is used for the two momentum predictors, whereas a Geometric Multigrid preconditioner is employed for the Helmholtz equations. 
 
\subsection{Case tests with analytical solution}
\label{ssec:convergence}

In order to verify the correctness of our implementation and to assess the convergence property of the scheme, we first perform numerical convergence studies in two and three dimensions, respectively. In two dimensions, we consider as a benchmark the classical Taylor-Green vortex \cite{green:1937} in the box \(\Omega = \left(0,2\pi\right)^2\), for which an analytical solution is available:

\begin{eqnarray}
\label{eq:greentaylor}
\mathbf{u}(\mathbf{x}, t) &=& \left(\begin{array}{c}
\cos(x_1)\sin(x_2)e^{-\frac{2t}{Re}} \\
-\sin(x_1)\cos(x_2)e^{-\frac{2t}{Re}}
\end{array} \right)  \\
p(\mathbf{x},t) &=& -\frac{1}{4}\left(\cos(2x_1) + \cos(2x_2)\right)e^{-\frac{4t}{Re}}.
\end{eqnarray}

In three dimensions, an analogous study has been carried out for the Arnoldi-Beltrami-Childress (ABC) flow, see e.g. \cite{galloway:1987}, whose exact solution is

\begin{eqnarray}
\label{eq:abc}
\mathbf{u}(\mathbf{x}, t) &=& \left(\begin{array}{c}
\left(\sin(x_3) + \cos(x_2)\right)e^{-\frac{t}{Re}} \\
\left(\sin(x_1) + \cos(x_3)\right)e^{-\frac{t}{Re}} \\
\left(\sin(x_2) + \cos(x_1)\right)e^{-\frac{t}{Re}}
\end{array} \right). \\
p(\mathbf{x},t) &=& -\sin(x_1)\cos(x_3) - \sin(x_2)\cos(x_1) - \sin(x_3)\cos(x_2).
\end{eqnarray}
 
\begin{figure}[!h]
	\centering
	\includegraphics[width=0.8\textwidth]{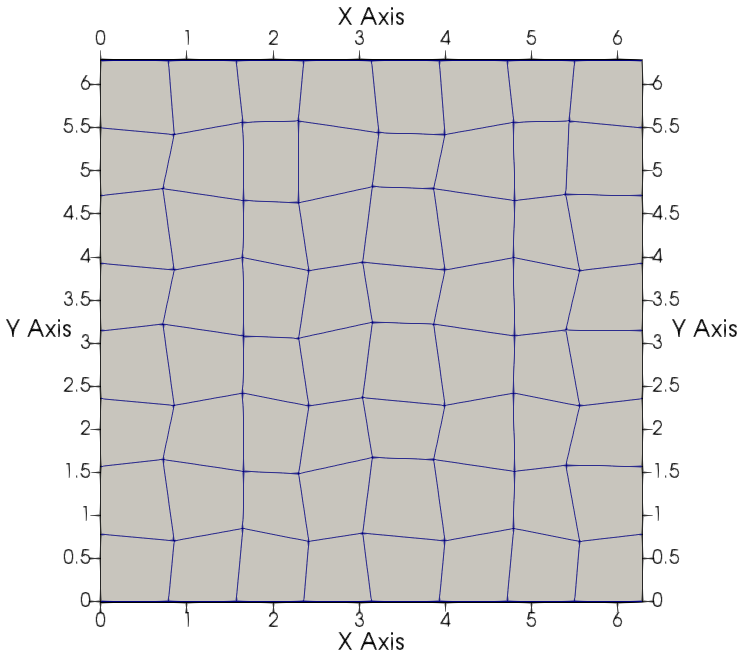}  
	\caption{Example of distorted mesh.}
	\label{fig:distorted}
\end{figure}

\begin{figure}[!h]
	\centering
	\includegraphics[width=0.45\textwidth]{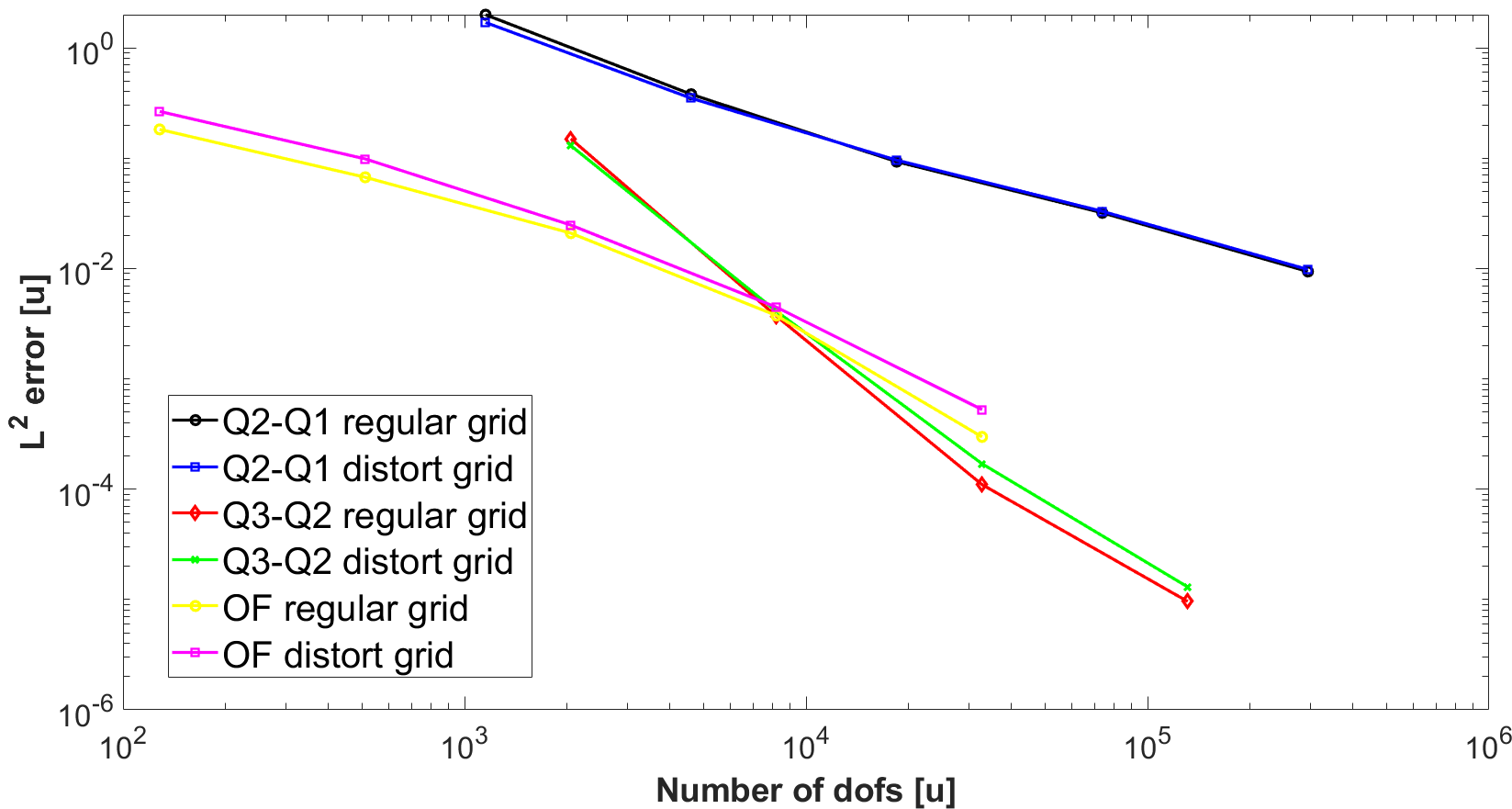} a)
	\includegraphics[width=0.45\textwidth]{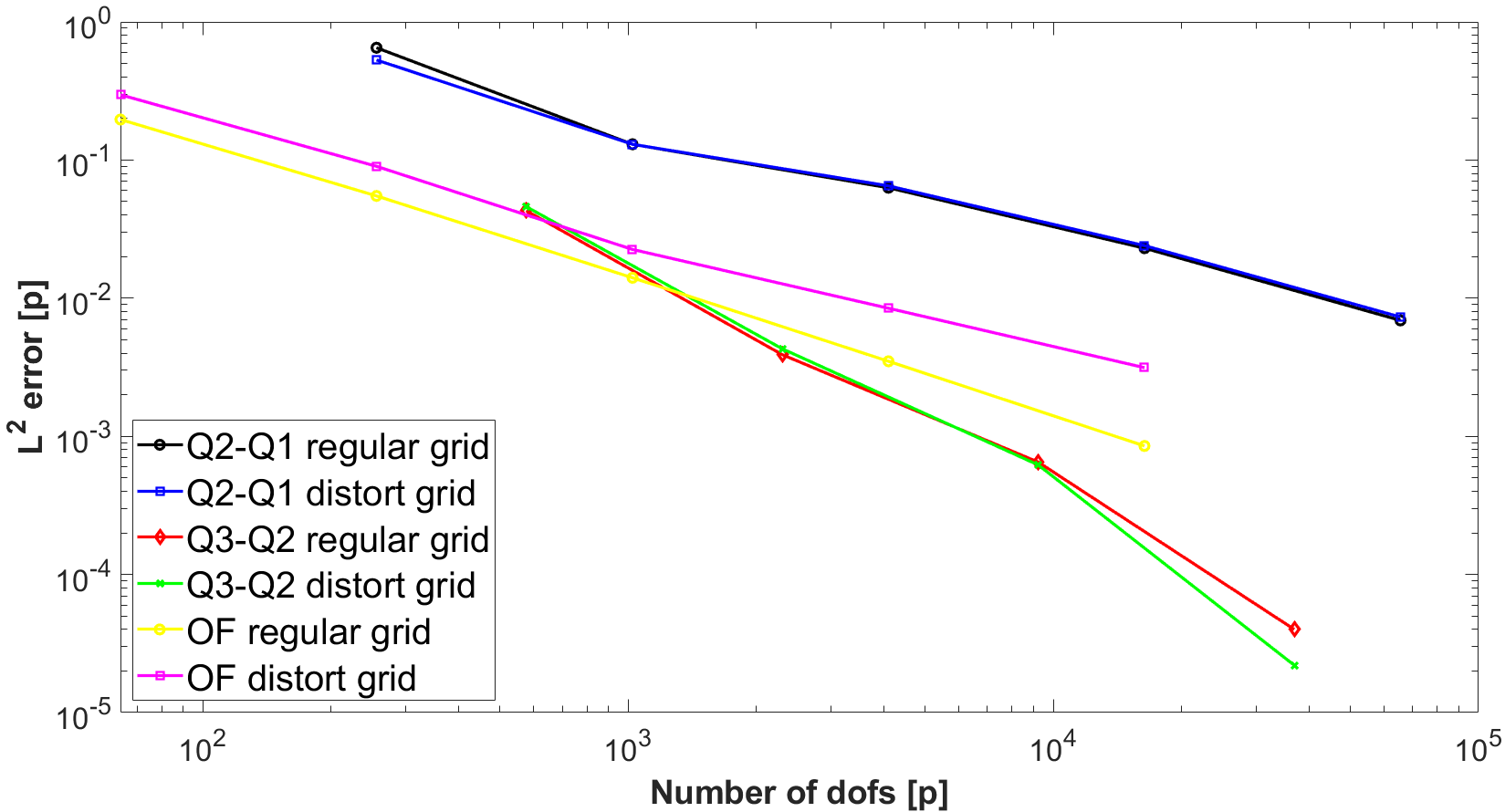} b)
	\includegraphics[width=0.45\textwidth]{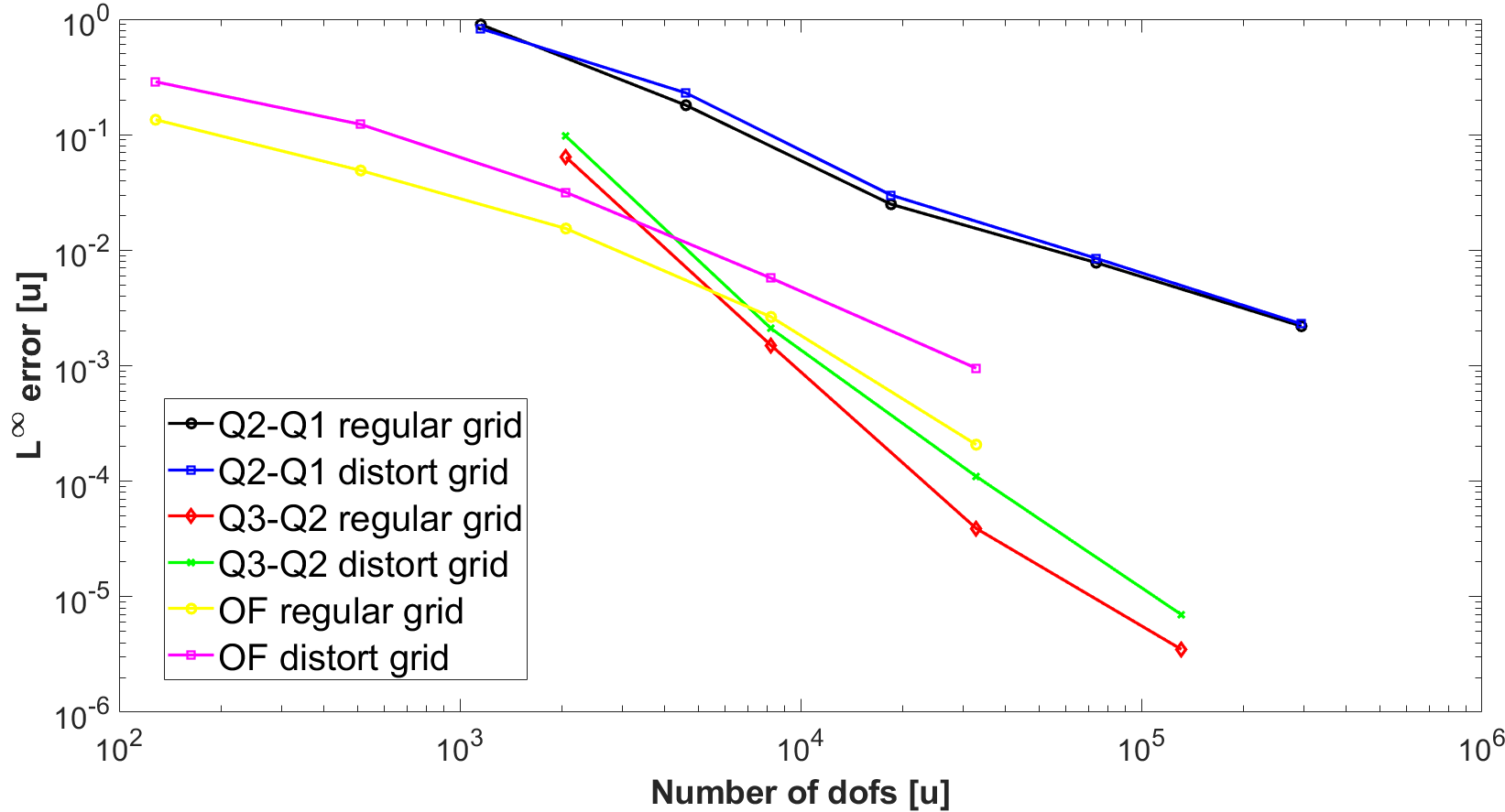} c)
	\includegraphics[width=0.45\textwidth]{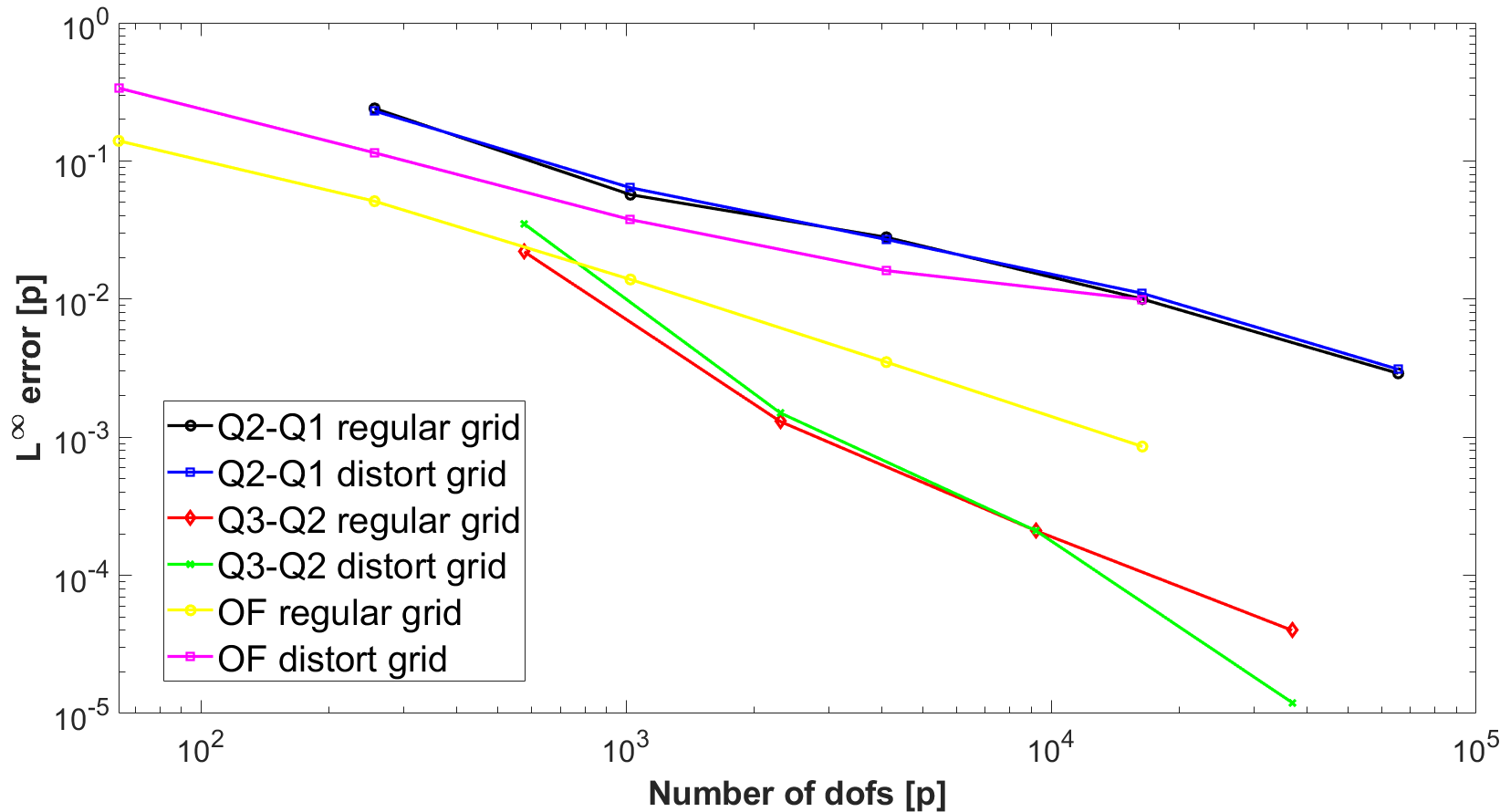} d)
 	\caption{Taylor-Green vortex at \(Re = 100\) and \(t = 3.2\), absolute errors as function of  the number of degrees of freedom, a) \(L^2\) errors for the velocity, b) \(L^2\) errors for the pressure, c) \(L^\infty\) errors for the velocity, d) \(L^\infty\) errors for the pressure. The black line denotes the solution with \(\mathbf{Q}_{2}-Q_{1}\) on regular grids, the blue line represents the results with \(\mathbf{Q}_{2}-Q_{1}\) on distort grids, the red line reports the results with \(\mathbf{Q}_{3}-Q_{2}\) on regular grids, the green line denotes the solution with \(\mathbf{Q}_{3}-Q_{2}\) on distort grids, the yellow line represents the results with OpenFoam on regular grids and the magenta line represents the results with OpenFoam on distort grids.}
	\label{fig:errors_TG}
\end{figure}

For the two dimensional case, we performed a convergence test at $ T= 3.2 $ for $ Re=100 $ 
starting with an initial Cartesian mesh of $ 8 \times 8 $ elements and doubling several  times the number of elements $N_{el} $ in each direction.
The time step was chosen so as to keep \(C = 1.63\) constant for all resolutions (hyperbolic scaling), so as to test the accuracy of the method for values of the time steps beyond the stability limit of explicit schemes but not large enough to affect the second order accuracy.
The results for the $\mathbf{Q}_2-Q_1$  and $\mathbf{Q}_3-Q_2$ cases are reported in Tables  \ref{tab:greentaylor_q2q1_u}, \hspace{0.05cm} \ref{tab:greentaylor_q2q1_p} and \ref{tab:greentaylor_q3q2_u}, \hspace{0.05cm} \ref{tab:greentaylor_q3q2_p}, respectively. It can be observed that the expected convergence rates are recovered, without the necessity of employing fixed point iterations to determine the velocity in the two stages. Analogous results are obtained, see Table \ref{tab:greentaylor_q2q1_dist_u}, \hspace{0.05cm} \ref{tab:greentaylor_q2q1_dist_p} if distorted meshes with analogous characteristics are employed.\\
The same test was repeated, for the case of $\mathbf{Q}_2-Q_1$ elements, using the  alternative methods \cite{bell:1989, guermond:1998} summarized in Section \ref{sec:modeleq}. 
It can be observed from the results reported in Tables \ref{tab:greentaylor_q2q1_bell_u}, \hspace{0.05cm} \ref{tab:greentaylor_q2q1_bell_p}, \hspace{0.05cm} \ref{tab:greentaylor_q2q1_guermond_u}, \hspace{0.05cm} \ref{tab:greentaylor_q2q1_guermond_p} \hspace{0.05cm} that, while the convergence rates are analogous, the relative errors in the $L^2$ norm are about $ 50\% $ smaller for the TR-BDF2 solver. \\
As mentioned in Section \ref{sec:modeleq}, when we increase the Courant number, also the TR-BDF2 scheme requires fixed point iterations in the momentum predictor stages in order to preserve its accuracy. As it can be noticed in Tables \ref{tab:greentaylor_q2q1_CFL_3_u}, \ref{tab:greentaylor_q2q1_CFL_3_p} that the second order convergence rate is still maintained.
\\~\\
For the three dimensional case, an analogous convergence test was performed again at $ T= 3.2 $ but using $ Re=1, $ due to the stability characteristics of the ABC flow, see e.g. the discussion in \cite{galloway:1987}. 
We have considered an initial Cartesian mesh of $8 \times 8 \times 8$ elements and we have refined the mesh by doubling each time the number of elements $ N_{el} $ in each direction, while keeping \(C = 1.63\) constant (hyperbolic scaling). The results for the $\mathbf{Q}_2-Q_1$  and $\mathbf{Q}_3-Q_2$ cases are reported in Tables  \ref{tab:abc_q2q1_u}, \hspace{0.05cm} \ref{tab:abc_q2q1_p} and \ref{tab:abc_q3q2_u}, \hspace{0.05cm} \ref{tab:abc_q3q2_p}, respectively. It can be observed that the expected convergence rates are recovered for the lower degree case, also in this case without the necessity of fixed point iterations, while less accurate results are obtained in the higher degree case. Since in this case the problem is diffusion dominated, rather than advection dominated, the loss of accuracy can be readily explained by the very large values obtained in this test for the parabolic stability parameter $\mu.$ Repeating the test at constant $\mu$ (parabolic scaling), one obtains the results displayed in Tables \ref{tab:abc_q2q1_muc_u}, \hspace{0.05cm} \ref{tab:abc_q2q1_muc_p} and \ref{tab:abc_q3q2_muc_u}, \hspace{0.05cm} \ref{tab:abc_q3q2_muc_p}, which show a clear improvement both in errors and convergence rates.
\\~\\
We have also used the two dimensional Taylor Green benchmark at \(Re = 100\) to compare our results  with analogous simulations carried out using using the classical PISO method \cite{issa:1986} as implemented in the OpenFoam package. In both cases, the computation was carried out at a very small value of the Courant number, so that the spatial discretization error is dominant. We are aware of the difficulties of comparing different discretizations schemes both in space and time implemented in different frameworks and, therefore, the following analysis has to be interpreted merely as first stress test to highlight the superior flexibility of the proposed DG implementation. We have performed a test using both $\mathbf{Q}_2-Q_1$ and $\mathbf{Q}_3-Q_2$ elements on regular and distorted meshes. An example of the coarsest distorted mesh is shown in Figure \ref{fig:distorted}, while the results of the convergence test for both \(L^{2}\) and \(L^{\infty}\) norms are reported in Figure \ref{fig:errors_TG}. While the OpenFoam discretization appears to outperform the $\mathbf{Q}_2-Q_1$  DG approximation at lower resolutions, it can be seen that it is much more sensitive to the mesh  distortion than DG approximations, especially with respect to \(L^\infty\) errors. Furthermore, as expected from polynomial approximation theory, the $\mathbf{Q}_3-Q_2$ DG approximation clearly shows its faster convergence properties, which are achieved within the same mathematical and implementation framework. Instead, higher order accuracy for finite volume formulations entails the use of complex and often \textit{ad hoc} reconstruction procedures with large stencils.

\FloatBarrier
 
\subsection{Two-dimensional lid driven cavity}
\label{ssec:2dcavity}

The lid driven cavity flow is a classical benchmark for the two-dimensional incompressible Navier-Stokes equations. Reference solutions obtained with high order techniques are reported, among many others, in \cite{auteri:2002, botella:1998, bruneau:2006}. For this two-dimensional problem, is it customary to represent the flow also in terms of the streamfunction \(\Psi\), which is defined as the solution of the Laplace problem

\begin{eqnarray}
\label{eq:streamline}
&&-\Delta\Psi = \nabla\times\mathbf{u} = \omega \qquad \text{in }\Omega\\
&&\Psi\rvert_{\partial\Omega} = 0 \nonumber
\end{eqnarray}

where the symbol \(\nabla\times\) denotes the curl operator and the vorticity is the scalar field defined as

\[\omega = \frac{\partial v}{\partial x_1} - \frac{\partial u}{\partial x_2}.\]

\begin{figure}[!h]
	\centering
 	\includegraphics[width=0.45\textwidth]{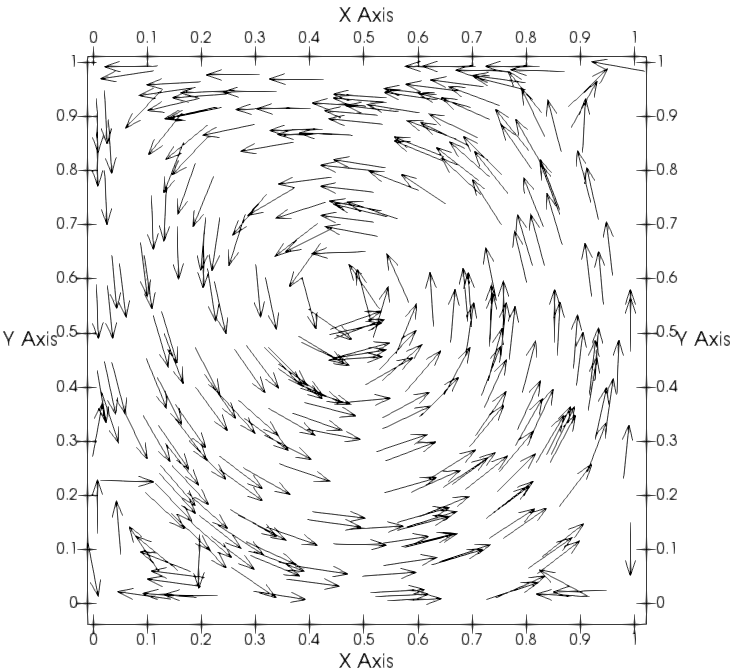} a)
 	\includegraphics[width = 0.45\textwidth]{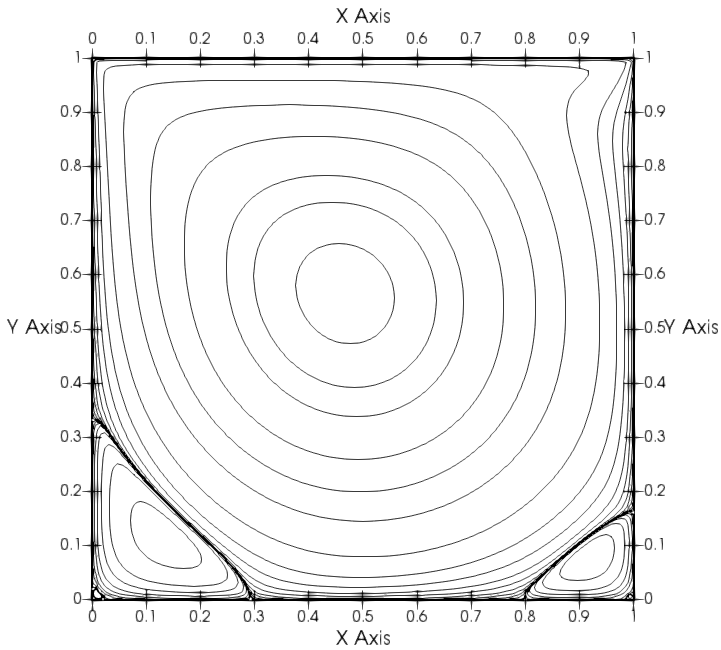} b)
 	\caption{Lid driven cavity benchmark at \(Re = 1000\): a) flow field, b) streamfunction contours. Contour values  are chosen as in \cite{bruneau:2006}.}
 	\label{fig:streamline_fixed}
\end{figure}

We consider the case \(Re = 1000\) computed with $\mathbf{Q}_2-Q_1$ elements on a Cartesian
mesh composed of $N_e=128$ square elements in each coordinate direction, with a  time step chosen so that the Courant number is approximately \(1.3\). The computation is performed until the steady state is reached up to a tolerance of \(10^{-7}\), which occurs around \(T = 70\). The streamfunction contours at steady state are shown in Figure \ref{fig:streamline_fixed} using the same isoline values as in \cite{bruneau:2006}. It can be observed that all the main flow structures are correctly reproduced.

\FloatBarrier
 
\begin{figure}[!h]
	\includegraphics[width=0.45\textwidth]{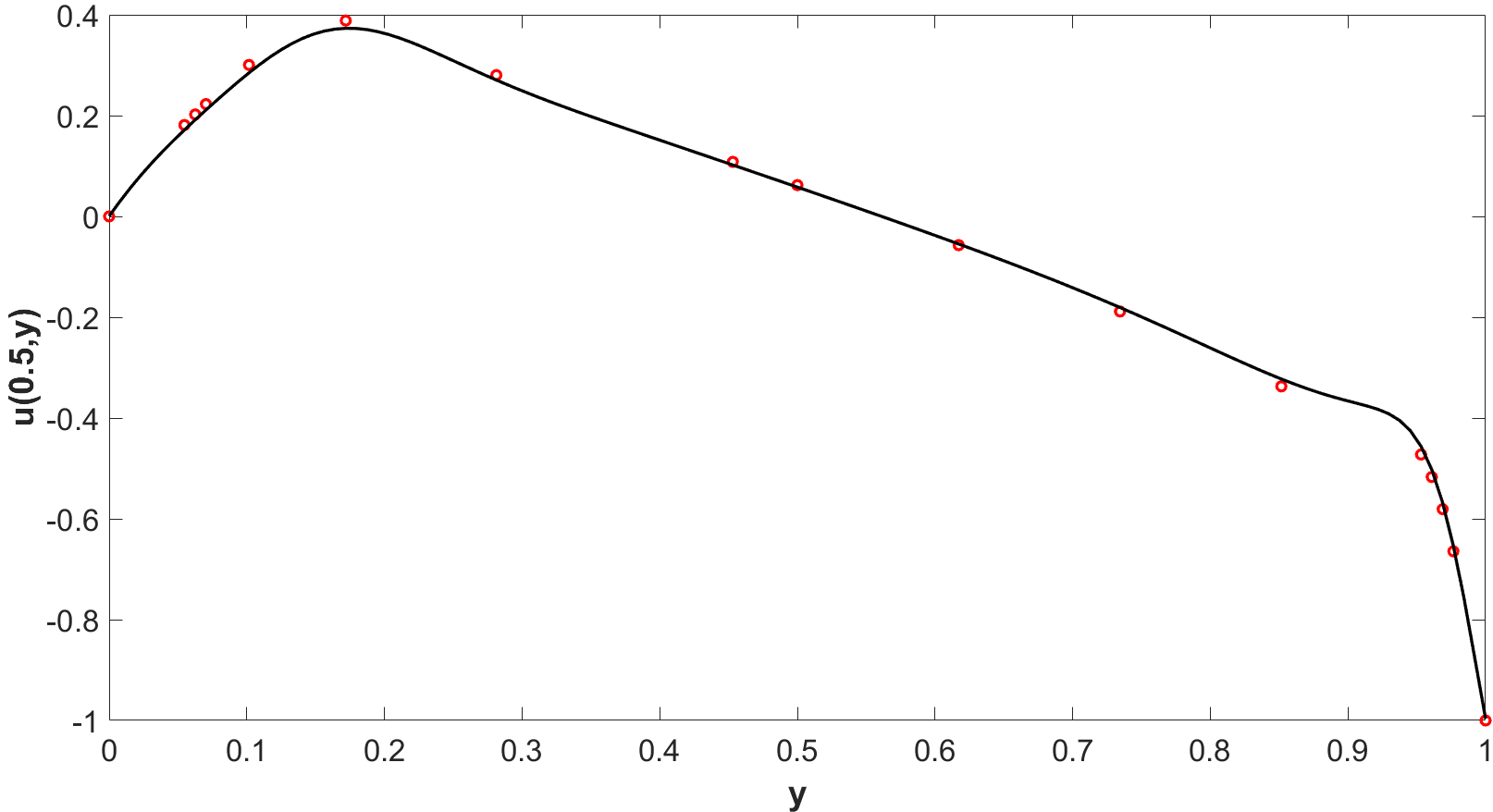} a)
 	\includegraphics[width=0.45\textwidth]{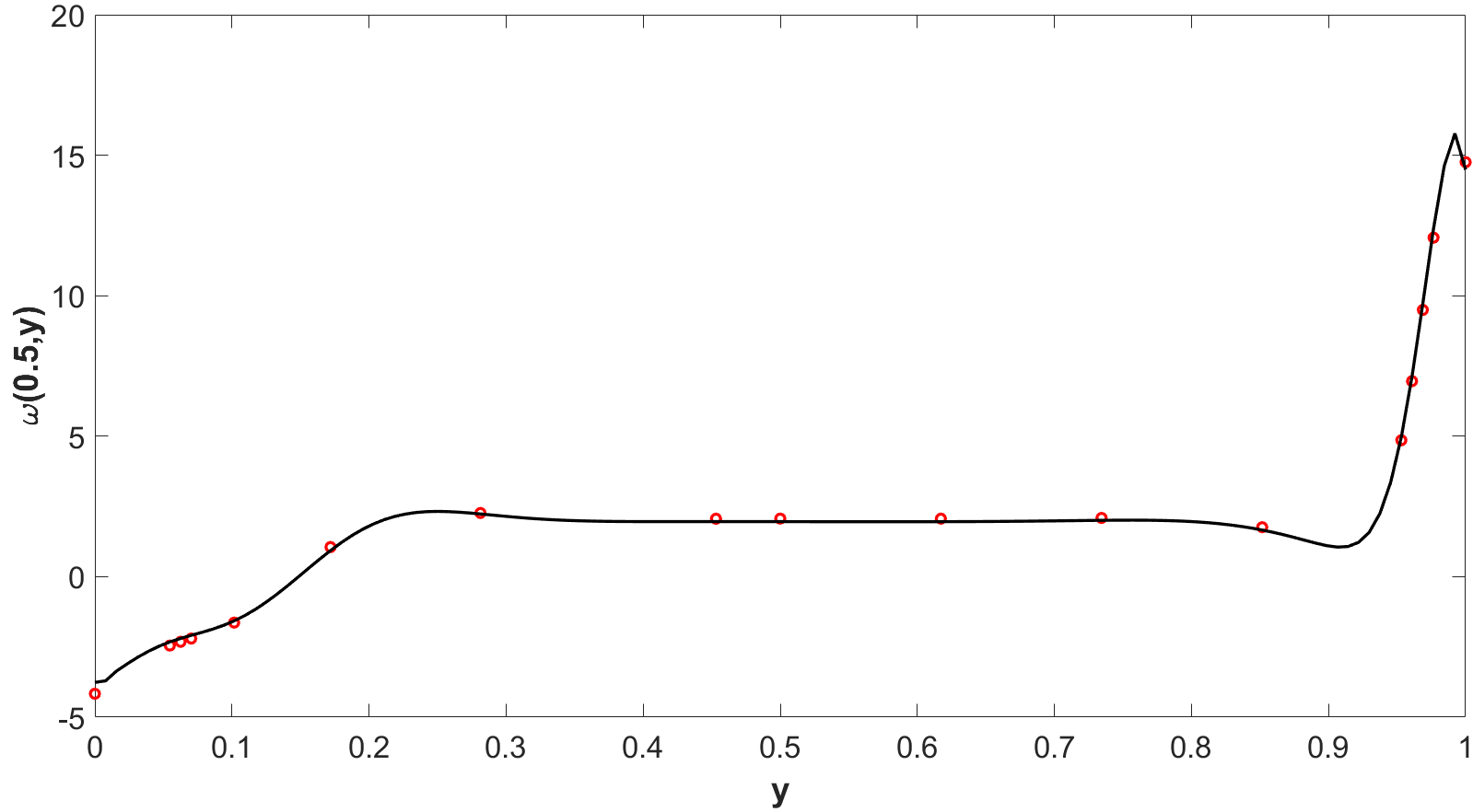} b)
 	\caption{Lid driven cavity benchmark at \(Re = 1000\): a) \(u\) velocity component values along the middle of the cavity, b) \(\omega\) values along the middle of the cavity. The continuous line denotes the numerical solution  and the dots the reference solution values from \cite{botella:1998}.}
 	\label{fig:u_om_middle_fixed}
\end{figure}

\begin{figure}[!h]
 	\includegraphics[width=0.45\textwidth]{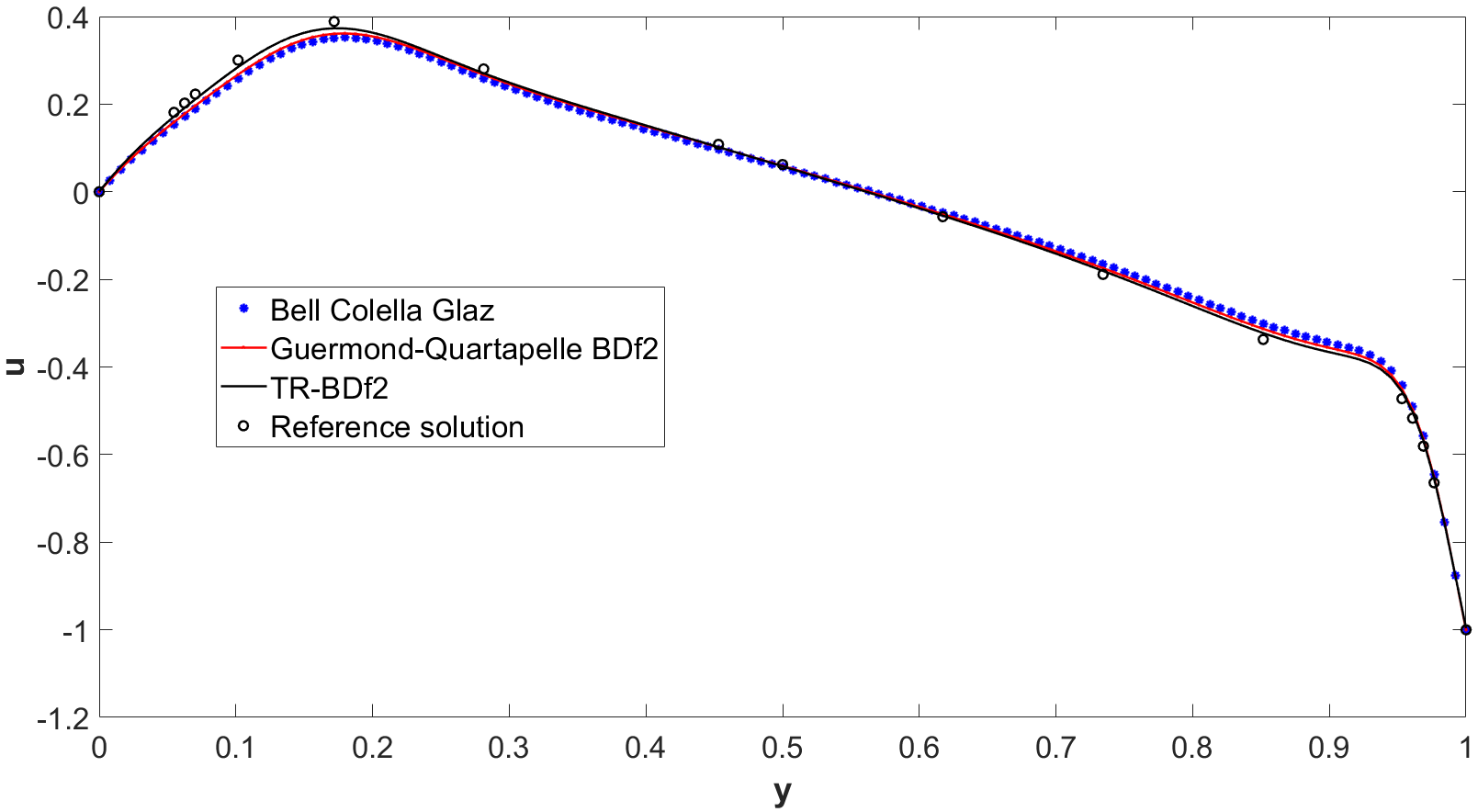} a)
 	\includegraphics[width=0.45\textwidth]{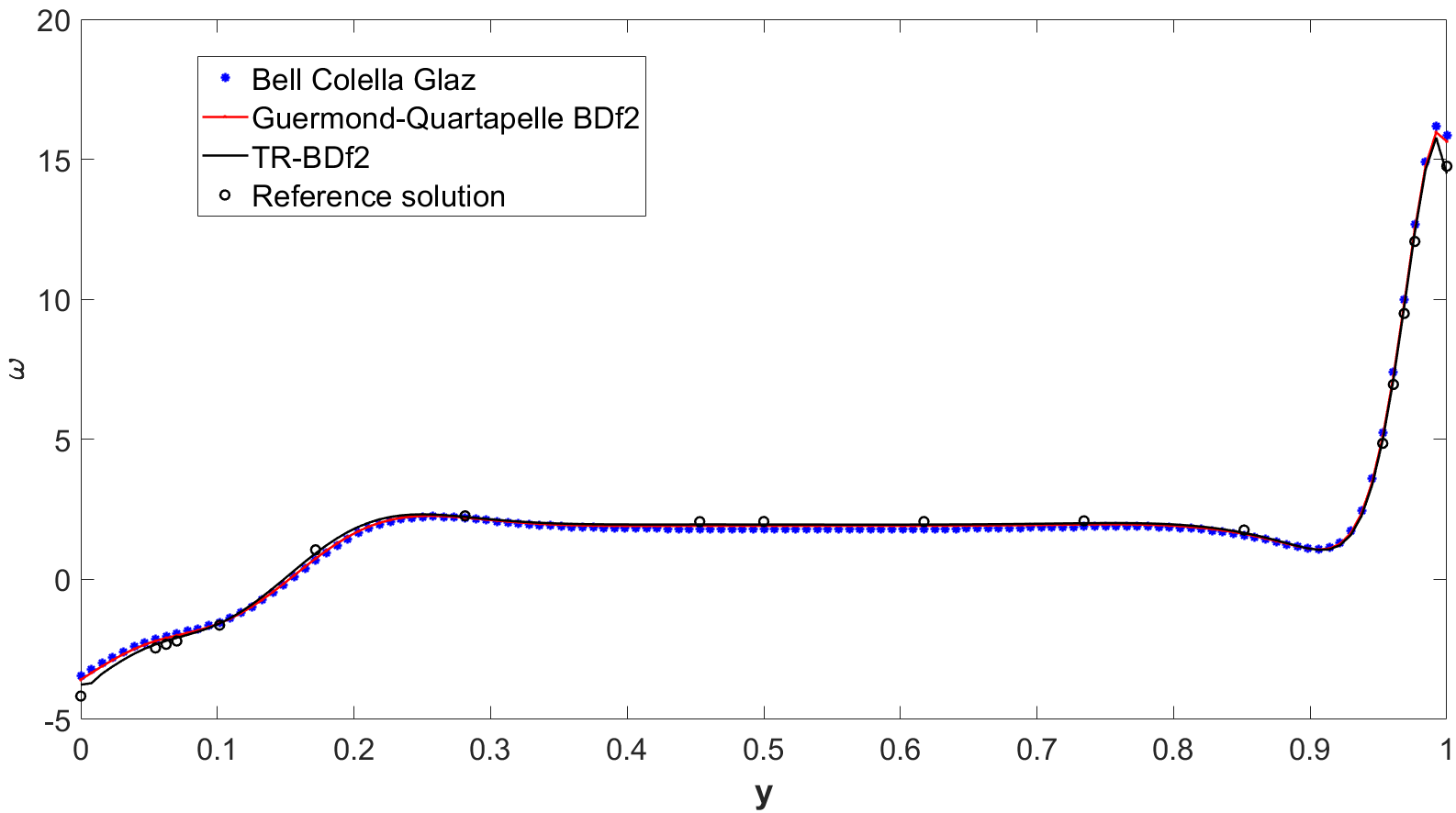} b)
 	\includegraphics[width=0.45\textwidth]{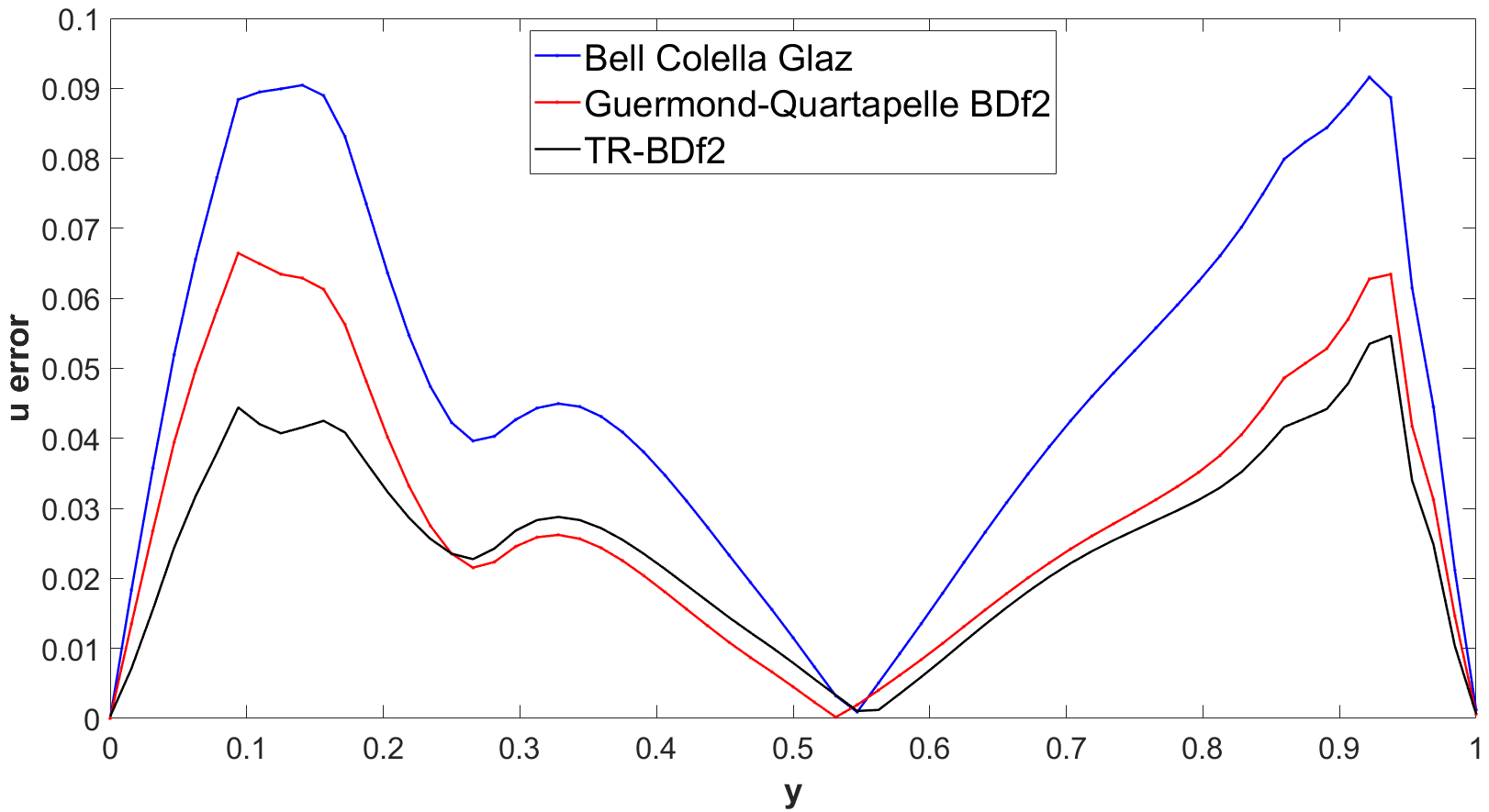} c)
 	\caption{Lid driven cavity benchmark at \(Re = 1000\): a) \(u\) velocity component values along the middle of the cavity, b) \(\omega\) values along the middle of the cavity. The continuous lines denote the numerical solutions	 with the methods \cite{bell:1989}, \cite{guermond:1998} and with the present method, the circles the reference solution values from \cite{botella:1998}, c) absolute error on \(u\) velocity component with respect to reference solution of \cite{botella:1998} interpolated along the middle of the cavity. The continuous black line denotes the result with the proposed method, the red one the results of \cite{guermond:1998} and the blue dots the results of \cite{bell:1989}.}
 	\label{fig:u_om_middle_fixed_comp}
\end{figure}

For a more quantitative comparison, we report in Figure \ref{fig:u_om_middle_fixed} the \(u\) component of the velocity and the vorticity \(\omega\) along the middle of the cavity, together with the reference results of \cite{botella:1998}. Good agreement with the reference solution is achieved. The maximum horizontal velocity along the centerline was computed as \(u_{max} = 0.3732\) which implies a relative error with respect to the reference solution of the order of \(10^{-2}\). The vorticity value at the center of the cavity was computed as \(\omega_{cen} = 1.9594\), which implies again a relative error with respect to the reference solution of the order of \(10^{-2}\). For comparison, the same test was repeated also using for the time discretization the parent methods described in  \cite{bell:1989}, \cite{guermond:1998}. The results are plotted in Figure \ref{fig:u_om_middle_fixed_comp}, highlighting the better performance of the proposed method based on TR-BDF2.
  
\FloatBarrier
 
\begin{figure}[!h]
	\centering
 	\includegraphics[scale = 0.35]{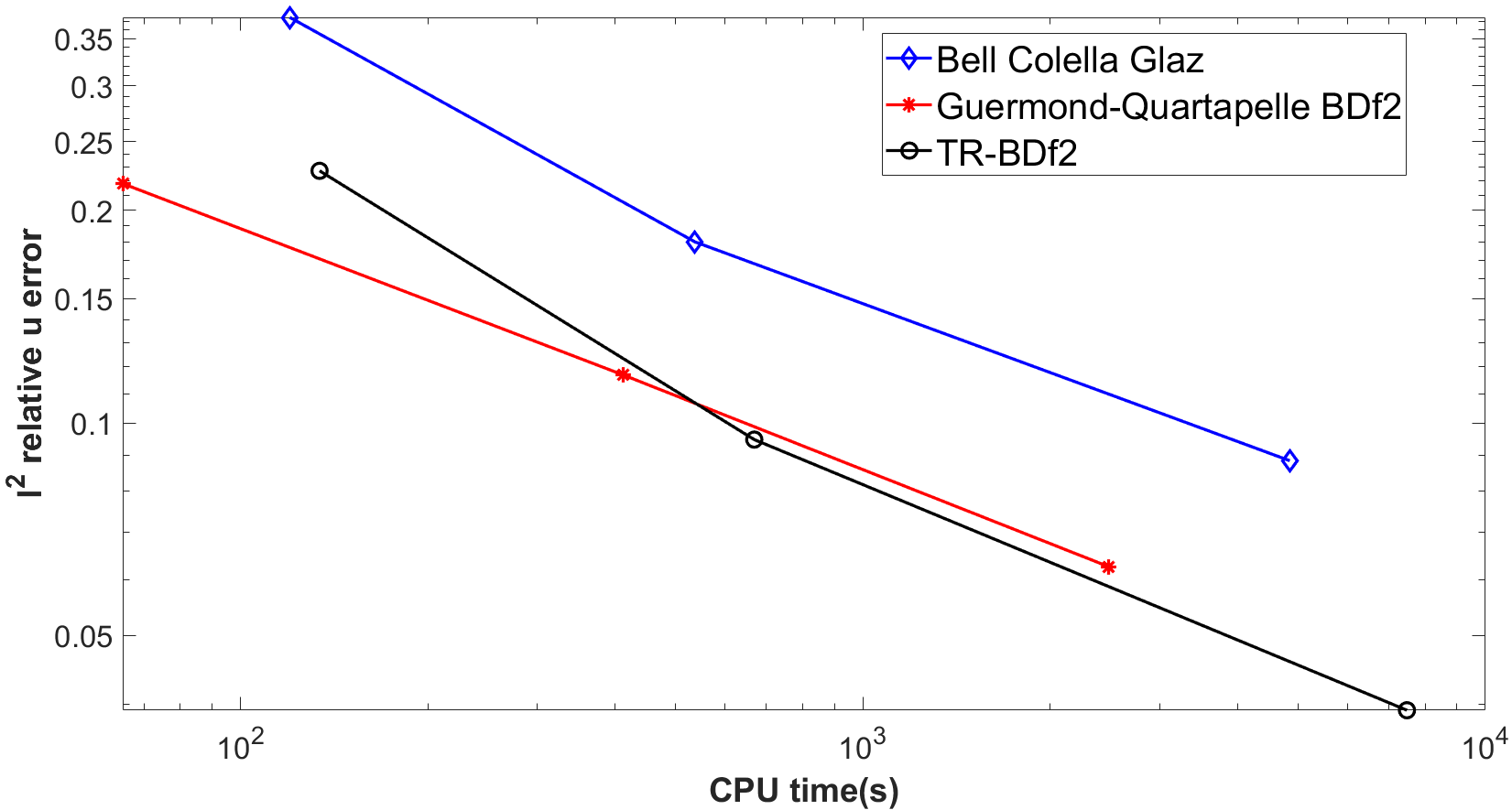}   	\caption{Lid driven cavity benchmark at \(Re = 1000\), \(l^2\) relative errors with respect to the CPU time required by the simulations with 8 MPI processes on Intel(R) Xeon(R) CPU Xeon E5-2640 v4 @ 2.4GHz. The continuous black line denotes the result with the proposed method, the red one the results of \cite{guermond:1998} and the blue one the results of \cite{bell:1989}.}
 	\label{fig:CPU_time_lid_driven}
\end{figure}

Moreover, we have compared the computational time required by the three methods for \(\mathcal{H} = \frac{1}{32},\frac{1}{64},\frac{1}{128}\), keeping the Courant number fixed. This assessment is important to show potential drawbacks of the two stage of the TR-BDF2 method with respect to the single stage methods employed in \cite{bell:1989, guermond:1998}. As shown in Figure \ref{fig:CPU_time_lid_driven}, the TR-BDF2 method shows superior efficiency with respect to the Bell-Colella-Glaz method, while it behaves similarly to the BDF2 method of \cite{guermond:1998}. Multistep methods, however, entail a memory overhead that is not appealing for large scale applications.

\FloatBarrier
 
\begin{figure}[!h]
 	\centering
 	\includegraphics[width=0.45\textwidth]{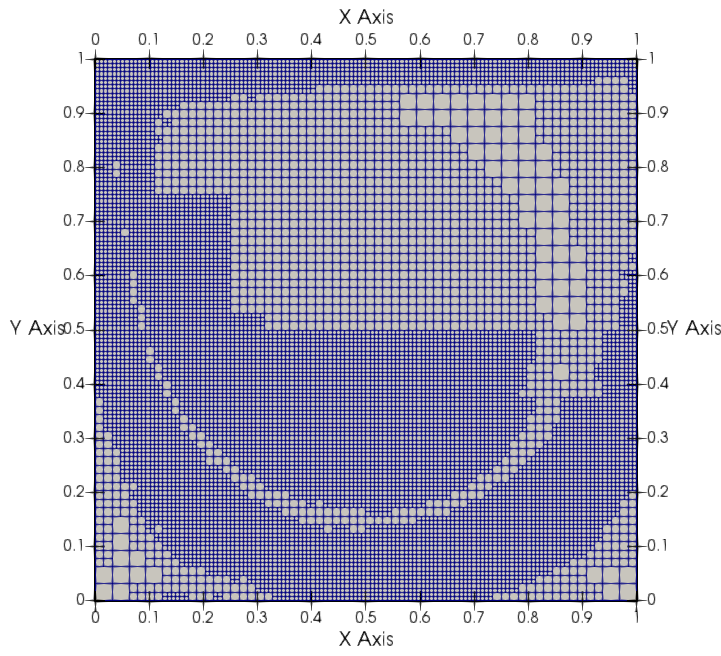} a)
 	\includegraphics[width=0.45\textwidth]{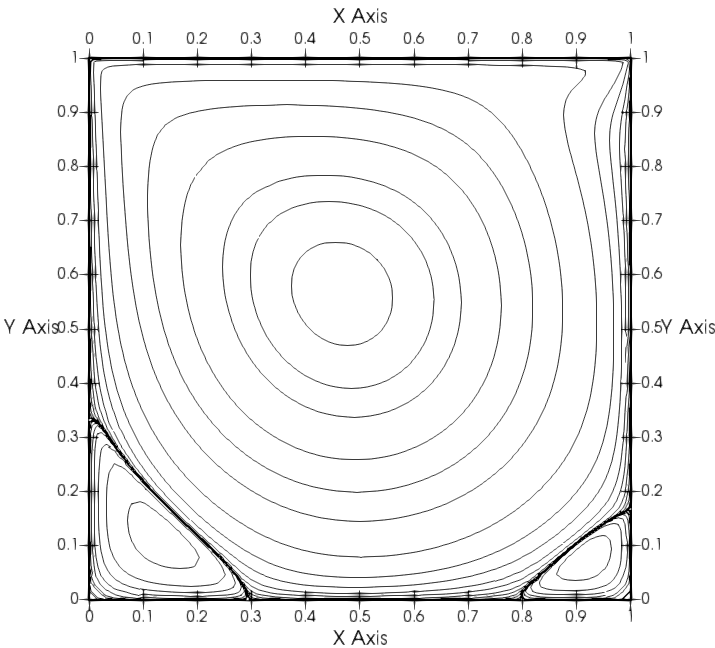} b)
 	\caption{Lid driven cavity benchmark at \(Re = 1000\), adaptive simulation: a) final mesh after adaptive refinement, b) streamfunction contours. Contour values as in \cite{bruneau:2006}.}
 	\label{fig:mesh_stream_adapted}
\end{figure}

We have also repeated this test using the adaptive tools present in the \textit{deal.II} library, as mentioned at the beginning of the Section. In each element \(K\) we define the quantity

\begin{equation}
	\eta_K = \text{diam}(K)^2\left\|\boldsymbol{\omega}\right\|^2_K
\end{equation}

that acts as local refinement indicator. We then started from a uniform Cartesian mesh with $ N_e= 8 $ in each coordinate direction and we allowed refinement or coarsening based on the distribution of the values of $ \eta_K, $ refining 10\% of the elements with largest indicator values and coarsening 30\% of the elements with the smallest indicator values. This remeshing procedure was carried out every 1000 time steps.
However, in order to avoid using a too coarse mesh for too long in the initial stages of the simulation, every 50 time steps the maximum difference between the velocities at two consecutive time steps was checked and the remeshing was performed whenever this quantity was greater then \(10^{-2}\). 
The minimum element diameter allowed was \(\mathcal{H} = \frac{1}{128}\), so as to obtain again \(C \approx 1.3\). A maximum element diameter equal to \(\frac{1}{32}\) was also required, in order to avoid an excessive reduction of the spatial resolution. 
The final adapted mesh and the streamline contours are reported in Figure \ref{fig:mesh_stream_adapted}. It can be observed that the refinement indicator allows to enhance automatically the resolution along the top boundary of the domain and in other regions of large vorticity values. 
 
\FloatBarrier
 
\begin{figure}[!h]
 	\includegraphics[width=0.45\textwidth]{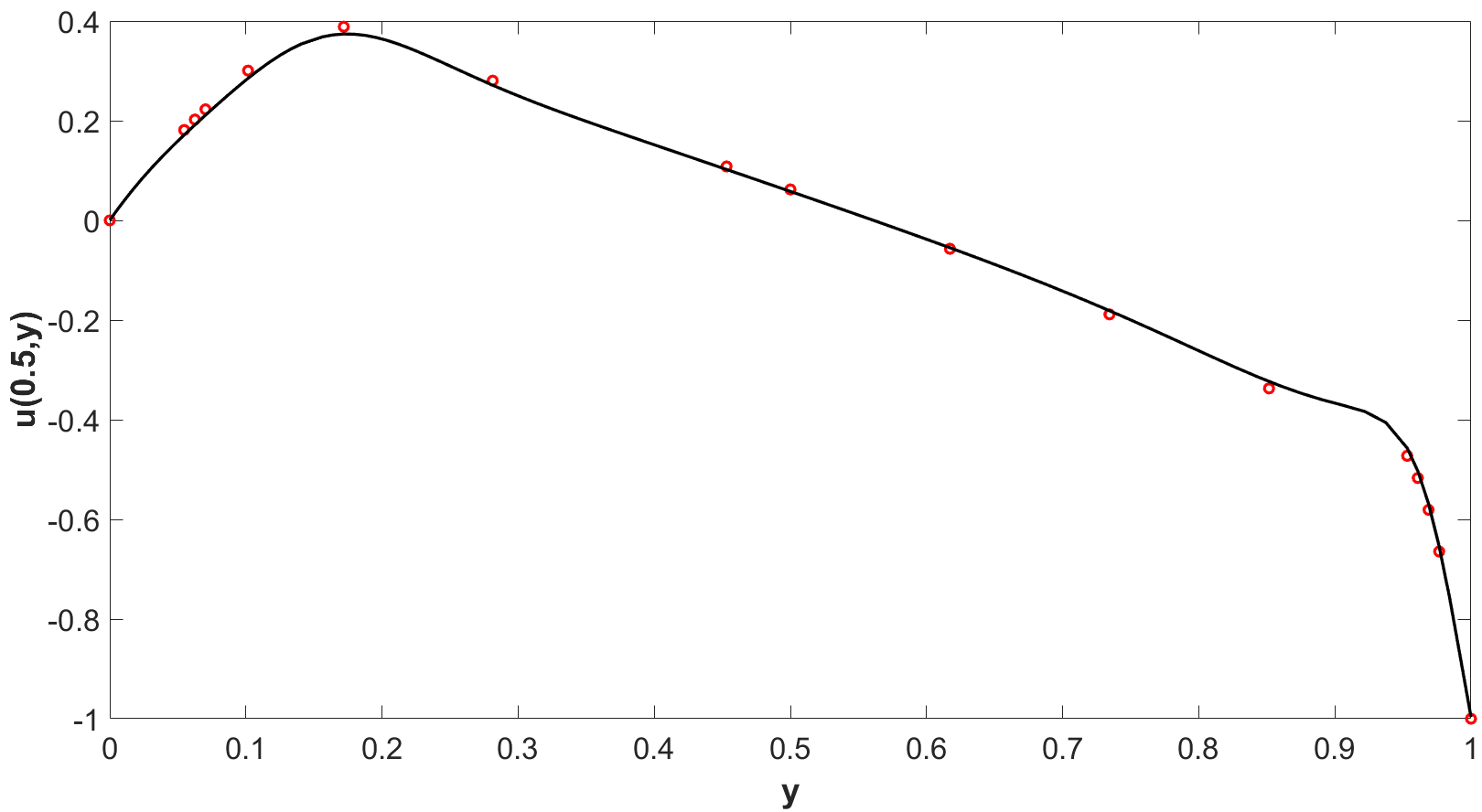} a)
 	\includegraphics[width=0.45\textwidth]{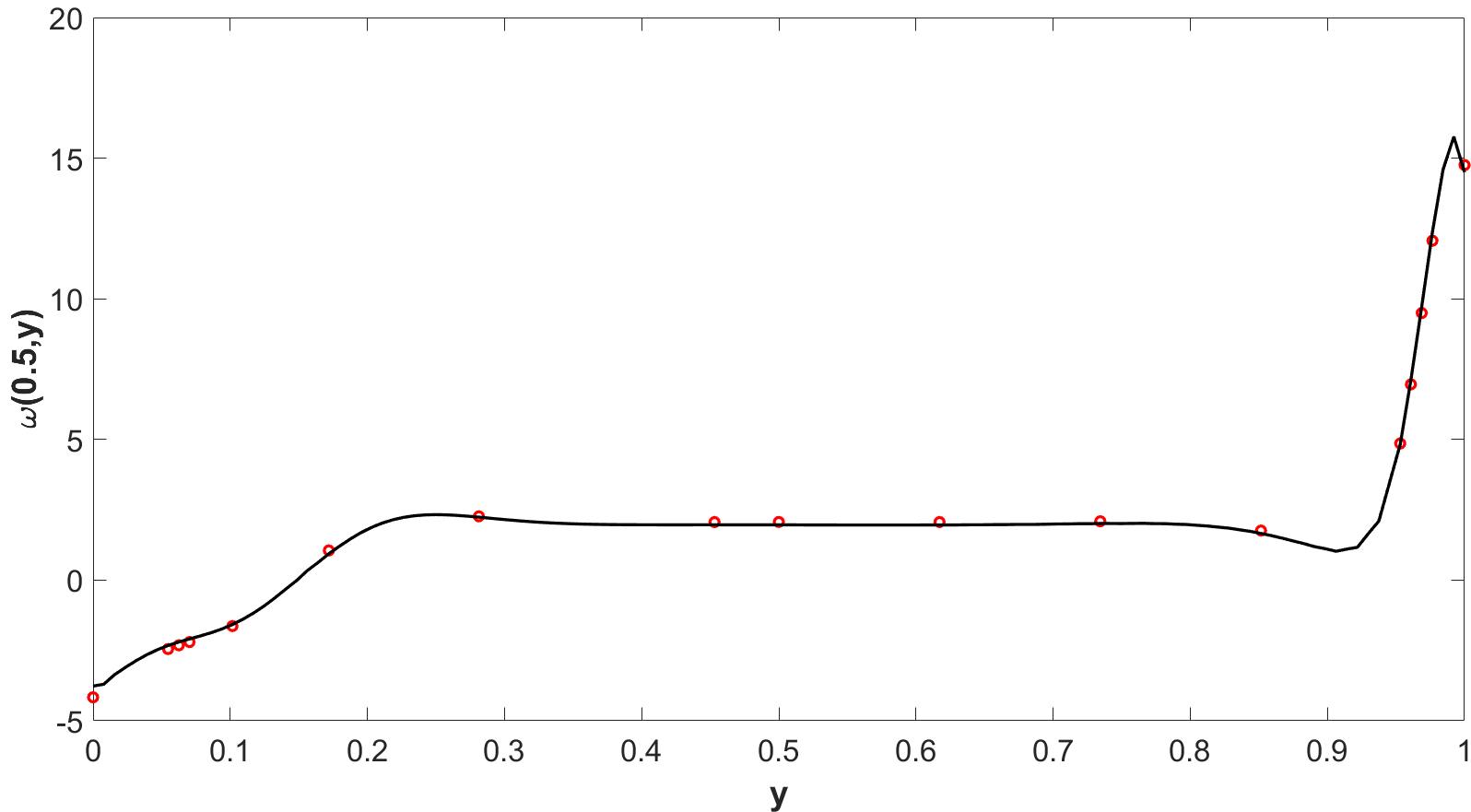} b)
 	\caption{Lid driven cavity benchmark at \(Re = 1000\), adaptive simulation: a) \(u\) velocity component values along the middle of the cavity,
 	b) \(\omega\) values along the middle of the cavity. The continuous line denotes the numerical solution  and the dots the reference solution values from \cite{botella:1998}.}
 	\label{fig:u_om_middle_adapted}
\end{figure}

\begin{figure}[!h]
	\begin{center}
	 	\includegraphics[width=0.7\textwidth]{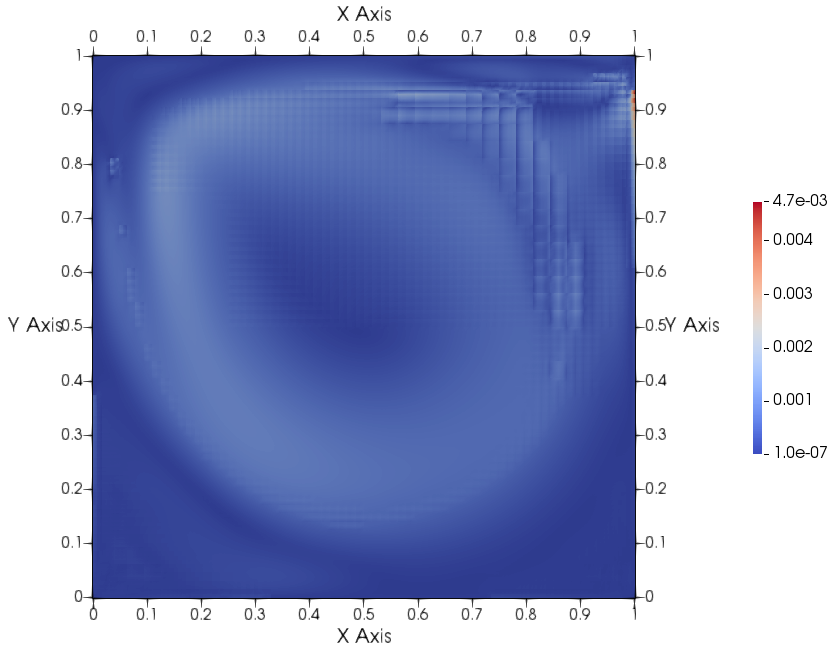}	
	\end{center}
	 \caption{Lid driven cavity benchmark at \(Re = 1000\), difference for velocity magnitude between the fixed grid simulation and the adaptive simulation (interpolated to the fixed grid).}
	 \label{fig:difference_fixed_adapted}
\end{figure}

For a more quantitative point of view, we compare again in Figure \ref{fig:u_om_middle_adapted} the \(u\) component of the velocity and the vorticity  \(\omega\) along the middle of the cavity with the reference results in \cite{botella:1998}. The maximum horizontal velocity along the centerline is now \(u_{max} = 0.3739\) which implies a relative error of the order of \(10^{-2}\), as in the corresponding non adaptive simulation. The vorticity value at the center of the cavity is now \(\omega_{cen} = 1.9652\), which also implies a relative error with respect to the reference solution of the order of \(10^{-2}\). In Figure \ref{fig:difference_fixed_adapted}, instead, the absolute difference between the velocities computed in the fixed mesh and adaptive simulations is plotted over the whole domain, showing that no substantial loss of accuracy has occurred. This result has been obtained with a reduction of about 25\% of the required computational time. While showing the potential of the adaptivity procedures available in the present implementation, this is still far from optimal. Experiments with more specific error indicators and less restrictive options for the refinement parameters will be carried out in future work.    
  
\subsection{Three-dimensional lid driven cavity}
\label{ssec:3dcavity}

We now consider the three-dimensional analog of the previously studied lid driven cavity benchmark.
Among several others, we consider the configuration and reference solutions provided in \cite{albensoeder:2005}, which we summarize here for convenience.
We consider a rectangular cavity of the size $d \times h \times l $ in the $ x, y $ and $ z $ direction, respectively. The flow is driven by the wall at $x=d/2,$ which moves tangentially in the $y$ direction with constant velocity $V.$ The length $d$ is used to introduce non dimensional space variables, so that the effective computational domain is given by 

$$\Omega =\left[-\frac{\Gamma}2,\frac{\Gamma}2\right ]\times\left [-\frac12,\frac12\right]\times\left[-\frac{\Lambda}2,\frac{\Lambda}2\right],$$

where the aspect ratios in the $x $ and $z $ directions are defined as 

\begin{equation}
 \label{eq:aspect}
 \Gamma=\frac{h}{d}, \ \ \ \ \ \ \Lambda= \frac{l}{d}.
\end{equation}

\begin{figure}[!h]
 	\includegraphics[width=0.45\textwidth]{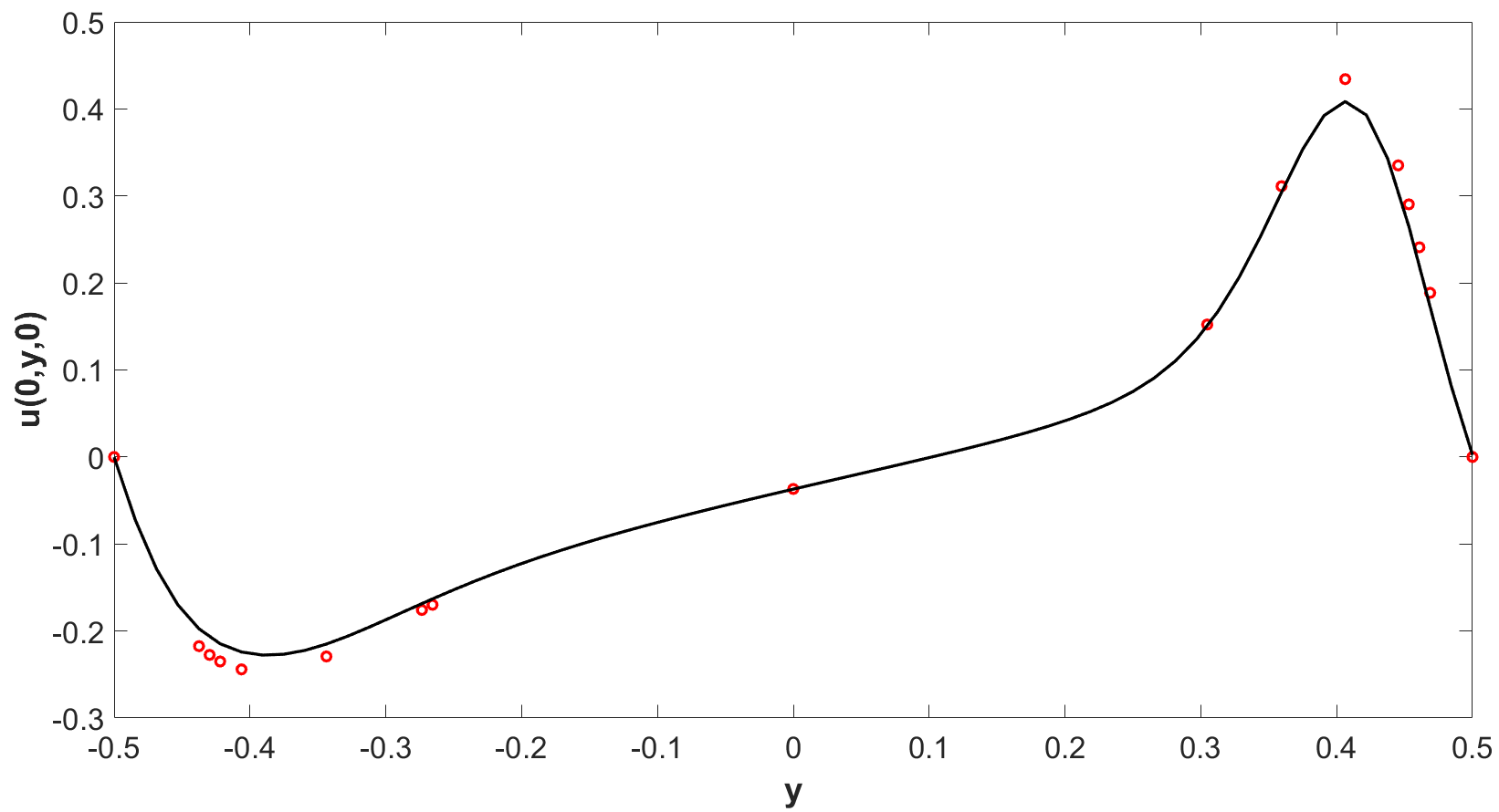} a)
	\includegraphics[width=0.45\textwidth]{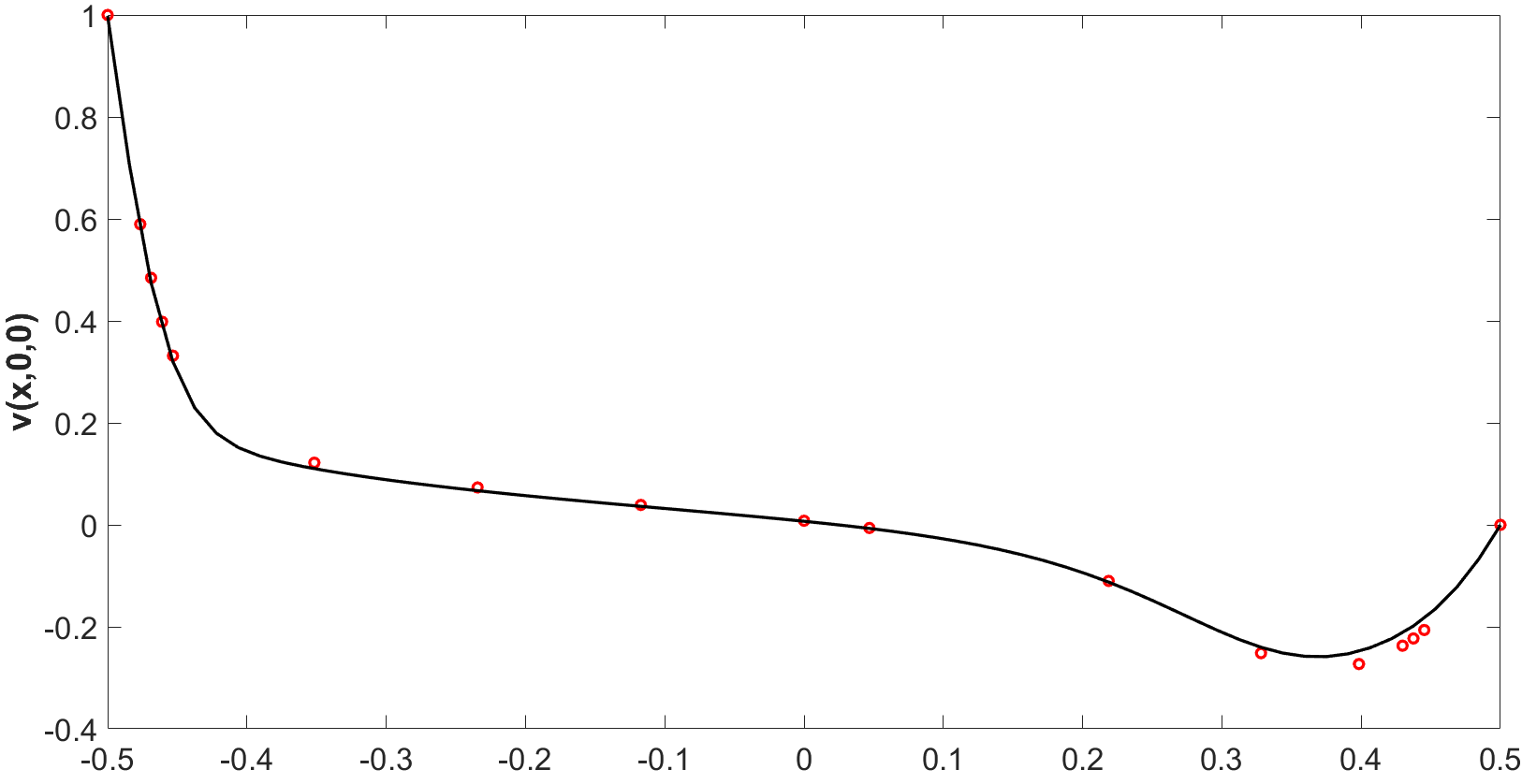} b) 
	\caption{3D lid driven cavity benchmark at \(Re = 1000\), fixed mesh simulation, a) \(v\) velocity component values along the $x $ axis b) \(u\) velocity component values along the $y $ axis. The continuous line denotes the numerical solution and the dots the reference solution values from \cite{albensoeder:2005}.}
 	\label{fig:3dcav_fixed}
\end{figure}

We have considered here the $ \Gamma=1,\Lambda=1 $ case, computed with $\mathbf{Q}_2-Q_1$ elements on a Cartesian mesh composed of \(64 \times 64 \times 48\) square elements, with a  time step chosen so that the Courant number is approximately \(1\). 
Notice that the same mesh was employed in \cite{albensoeder:2005}, which however employed a much more accurate spectral collocation method.
The computation is performed until the steady state is reached up to a tolerance of \(10^{-4}\), which  is achieved around \(T = 40\). 
We take as reference results those presented in Tables 5 and 6 in \cite{albensoeder:2005}. Notice that, in that paper, a different non dimensional scaling is employed, so that their results have been appropriately rescaled in order to compare them with those obtained here. In Figure \ref{fig:3dcav_fixed} we report the results for the \(v\) velocity component values along the $x$ axis and the \(u\) component of the velocity along the $y$ axis, respectively. We see that, in spite of the relatively coarse mesh, a reasonable accuracy is achieved.
   
\FloatBarrier
 
\begin{figure}[!h]
  	\includegraphics[width=0.29\textwidth]{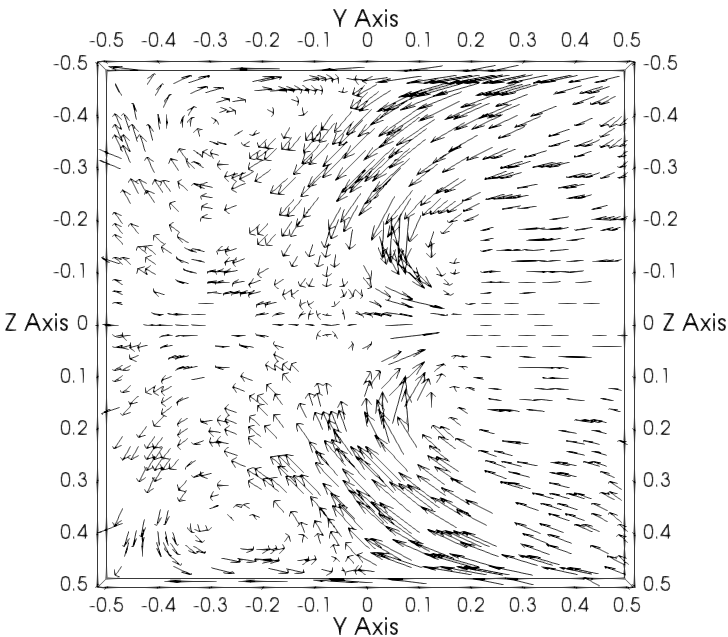} a)
  	\includegraphics[width=0.29\textwidth]{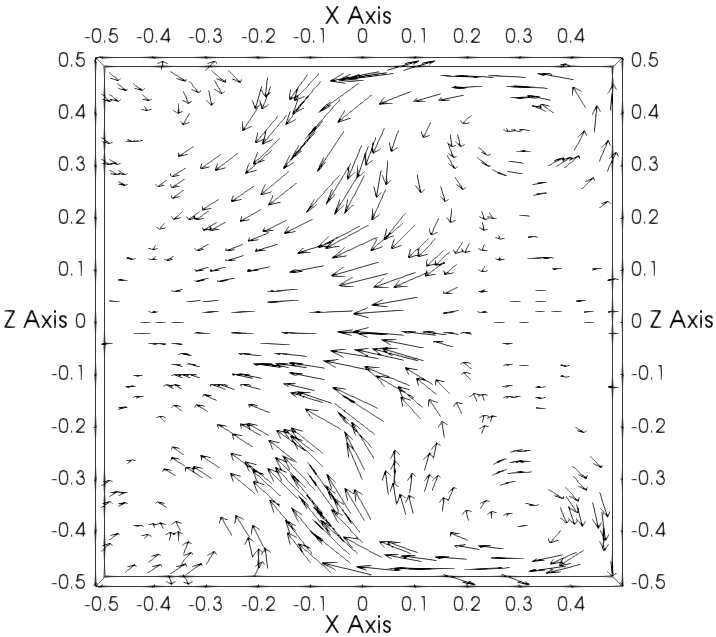} b)
  	\includegraphics[width=0.29\textwidth]{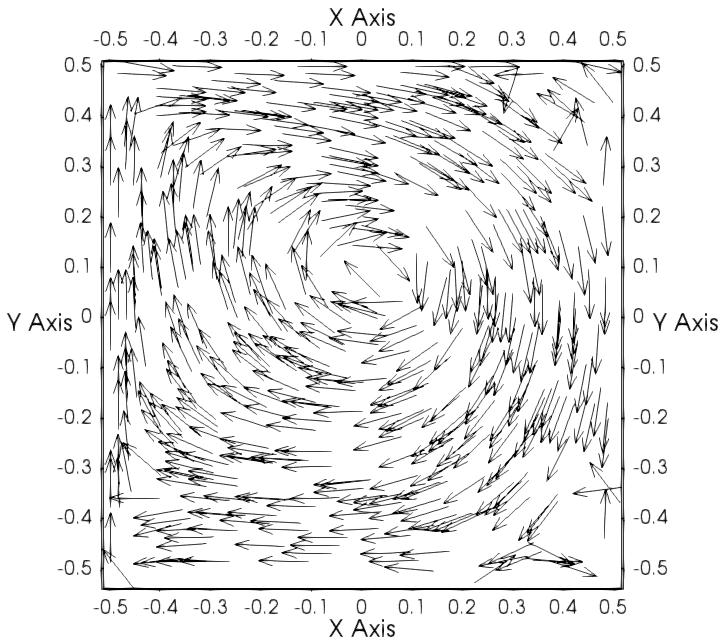} c)
  	\caption{3D lid driven cavity benchmark at \(Re = 1000\), a) Flow field vectors for the plane \(x = 0\), b) Flow field vectors for the plane \(y = 0\), c) Flow field vectors for the plane \(z = 0\).}
  	\label{fig:velocity_fields_3D_lid_driven}
\end{figure}

\begin{figure}[!h]
  	\includegraphics[width=0.29\textwidth]{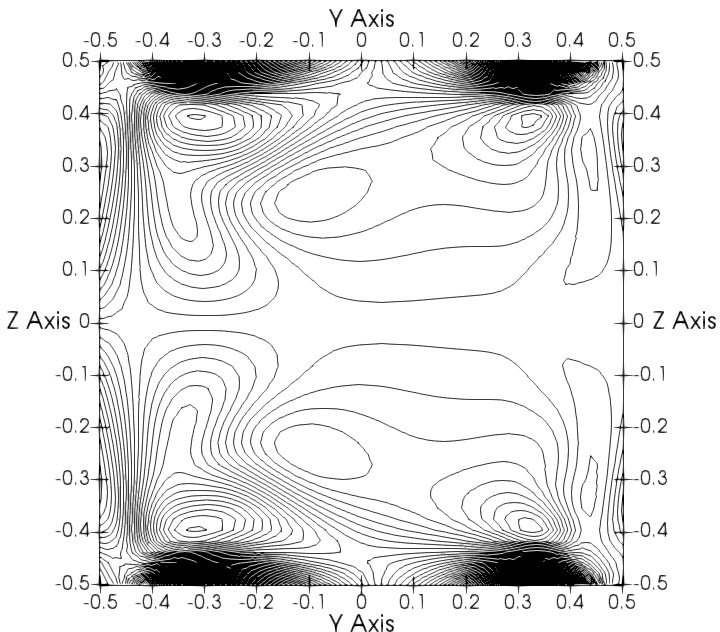} a)
  	\includegraphics[width=0.29\textwidth]{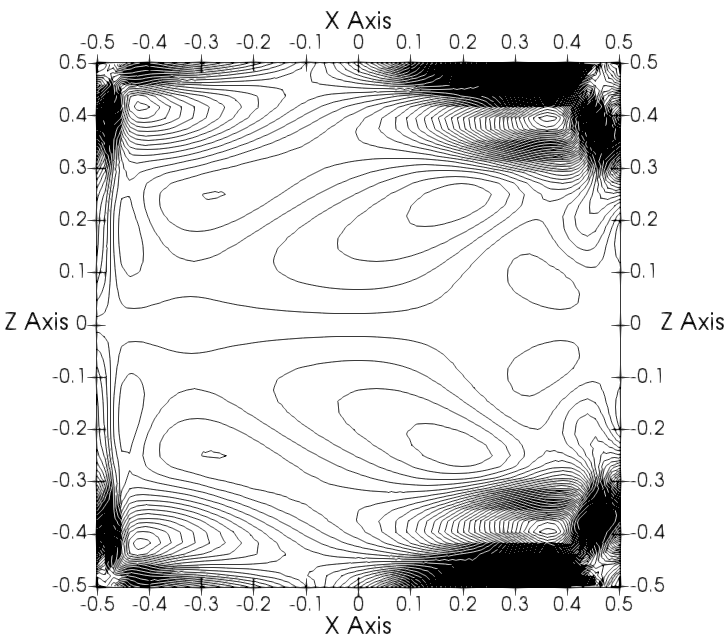} b)
  	\includegraphics[width=0.29\textwidth]{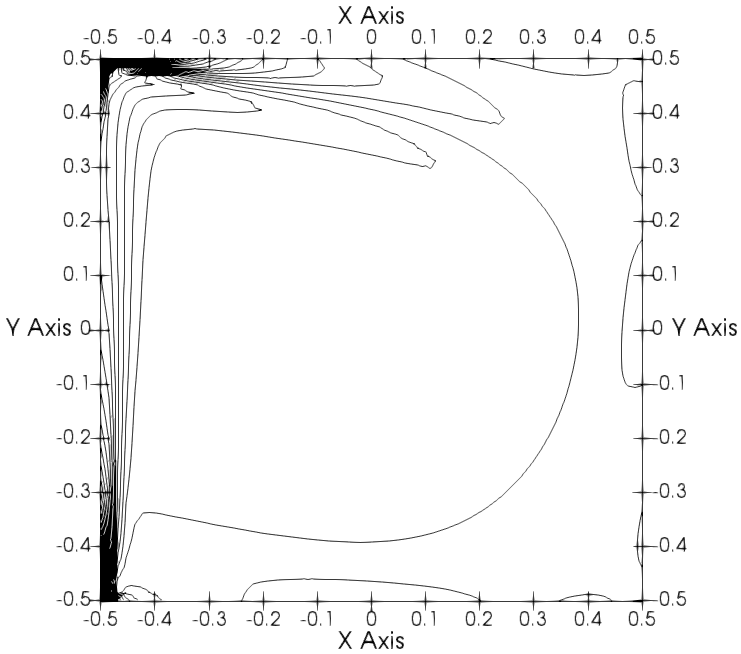} c)
  	\caption{3D lid driven cavity benchmark at \(Re = 1000\), a) Vorticity (\(\omega_x\)) contours at \(x = 0\), b) Vorticity (\(\omega_y\)) contours at \(y = 0\), Vorticity (\(\omega_z\)) contours at \(z = 0\)}
  	\label{fig:vorticity_contours_3D_lid_driven}
\end{figure}

In Figure \ref{fig:velocity_fields_3D_lid_driven} and \ref{fig:vorticity_contours_3D_lid_driven} we show instead the velocity field on the three median plane sections of the cavity, highlighting the presence of vortices near the centerline of the cavity. The results are in good qualitative agreement with those reported in \cite{jiang:1994}.
 
\begin{figure}[!h]
  	\includegraphics[width=0.45\textwidth]{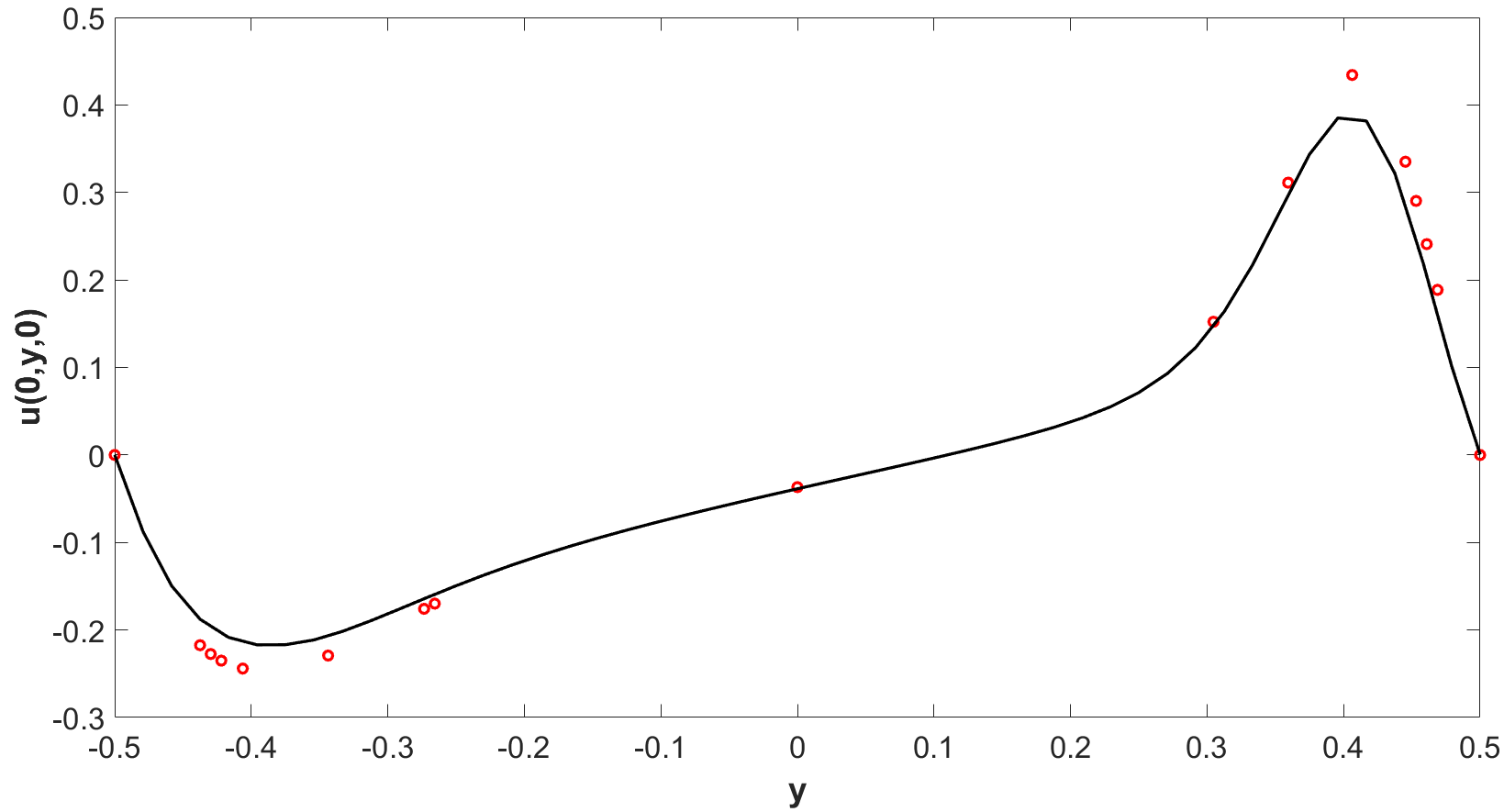} a)
  	\includegraphics[width=0.45\textwidth]{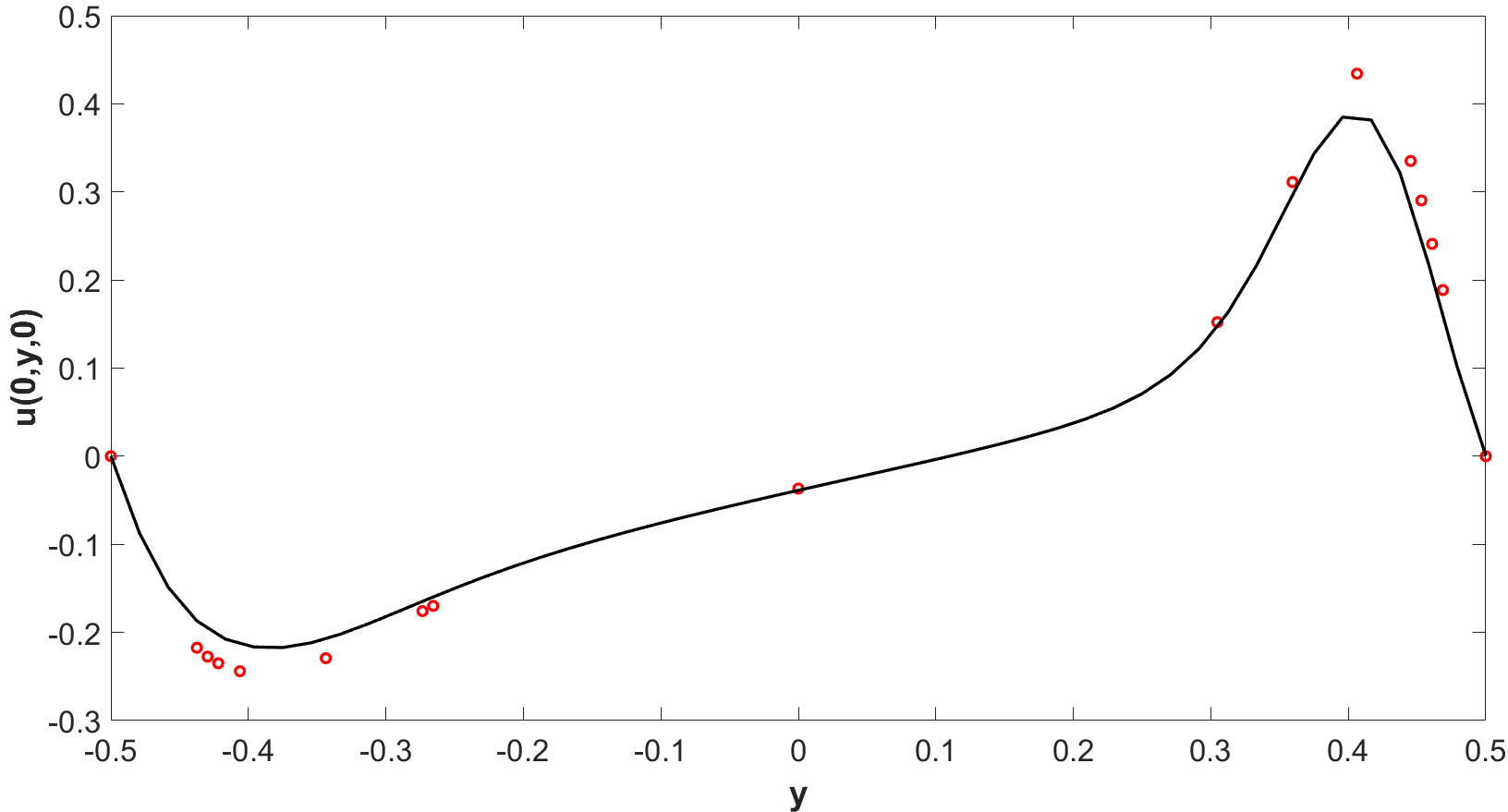} b) 
  	\caption{3D lid driven cavity benchmark at \(Re = 1000\), a) \(u\) velocity component values along the $y $ axis for adaptive mesh simulation, b) \(u\) velocity component values along the $y $ axis for fixed grid simulation. The continuous line denotes the numerical solution and the dots the reference solution values from \cite{albensoeder:2005}.}
  	\label{fig:3dcav_u_y_adapt}
\end{figure}

\begin{figure}[!h]
  	\includegraphics[width=0.45\textwidth]{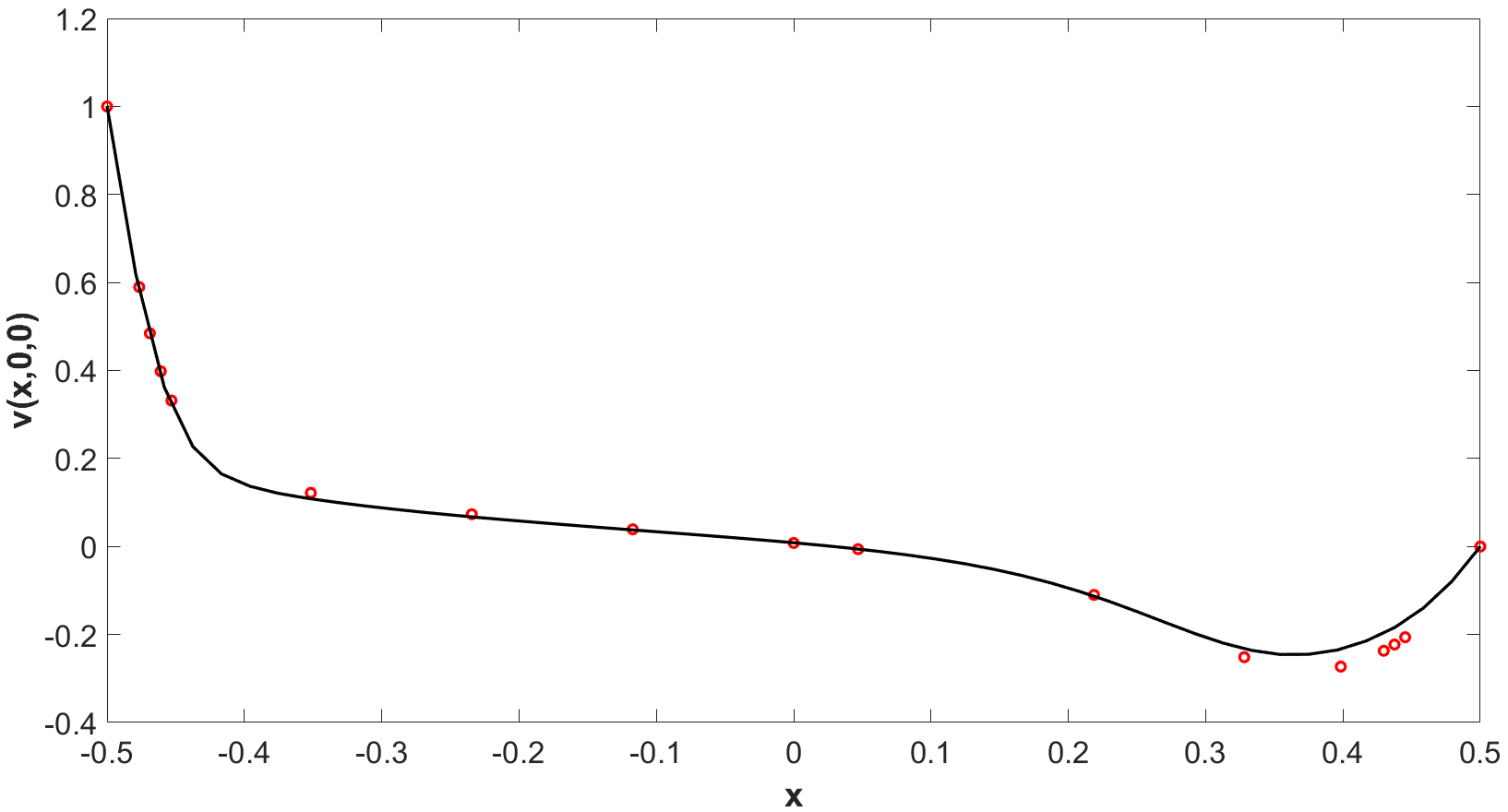} a)
  	\includegraphics[width=0.45\textwidth]{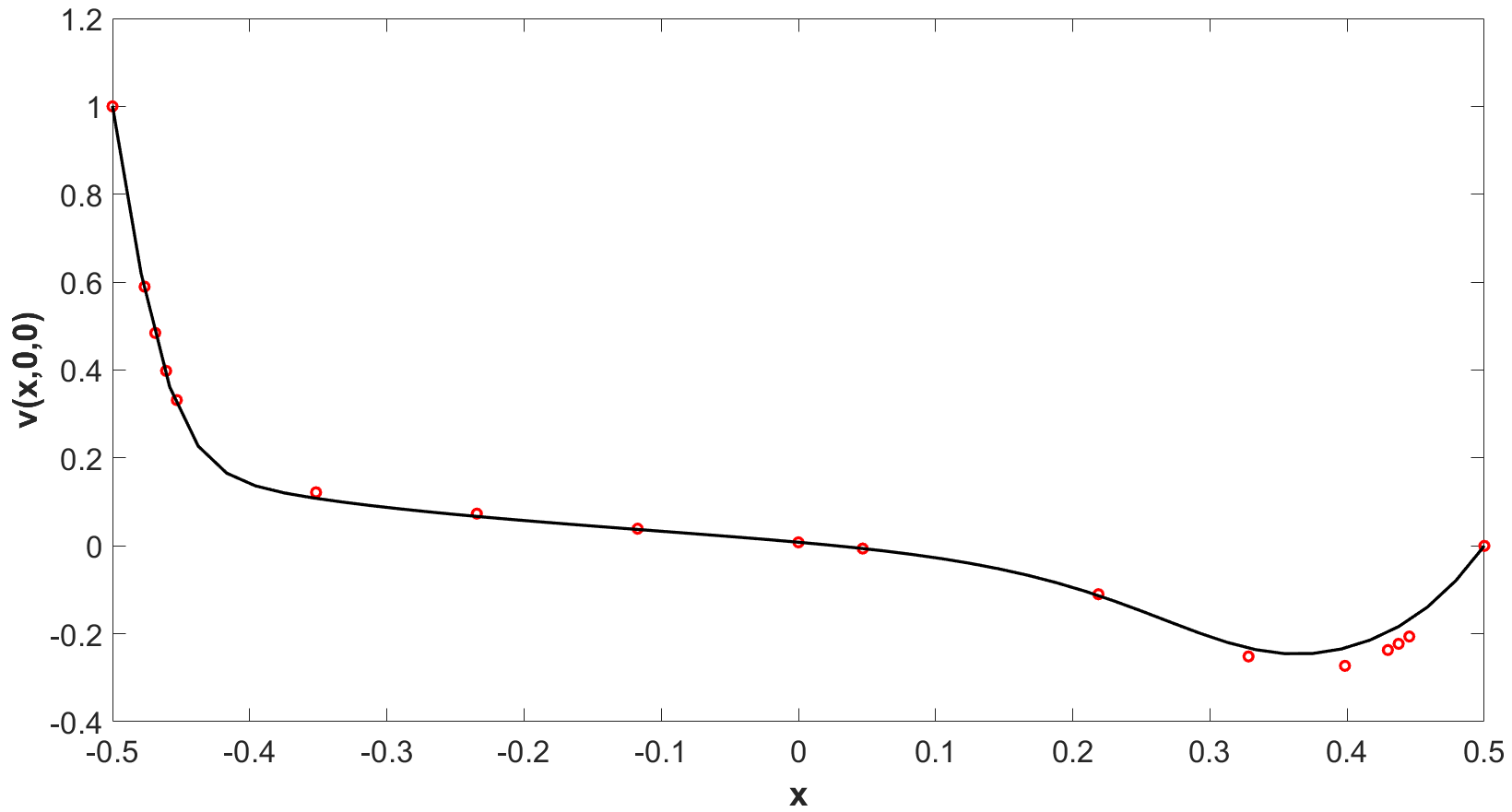} b) 
  	\caption{3D lid driven cavity benchmark at \(Re = 1000\), a) \(v\) velocity component values along the $x $ axis for adaptive mesh simulation, b) \(v\) velocity component values along the $x $ axis for fixed grid simulation. The continuous line denotes the numerical solution and the dots the reference solution values from \cite{albensoeder:2005}.}
  	\label{fig:3dcav_v_x_adapt}
\end{figure}

\begin{figure}[!h]
  	\includegraphics[width=0.45\textwidth]{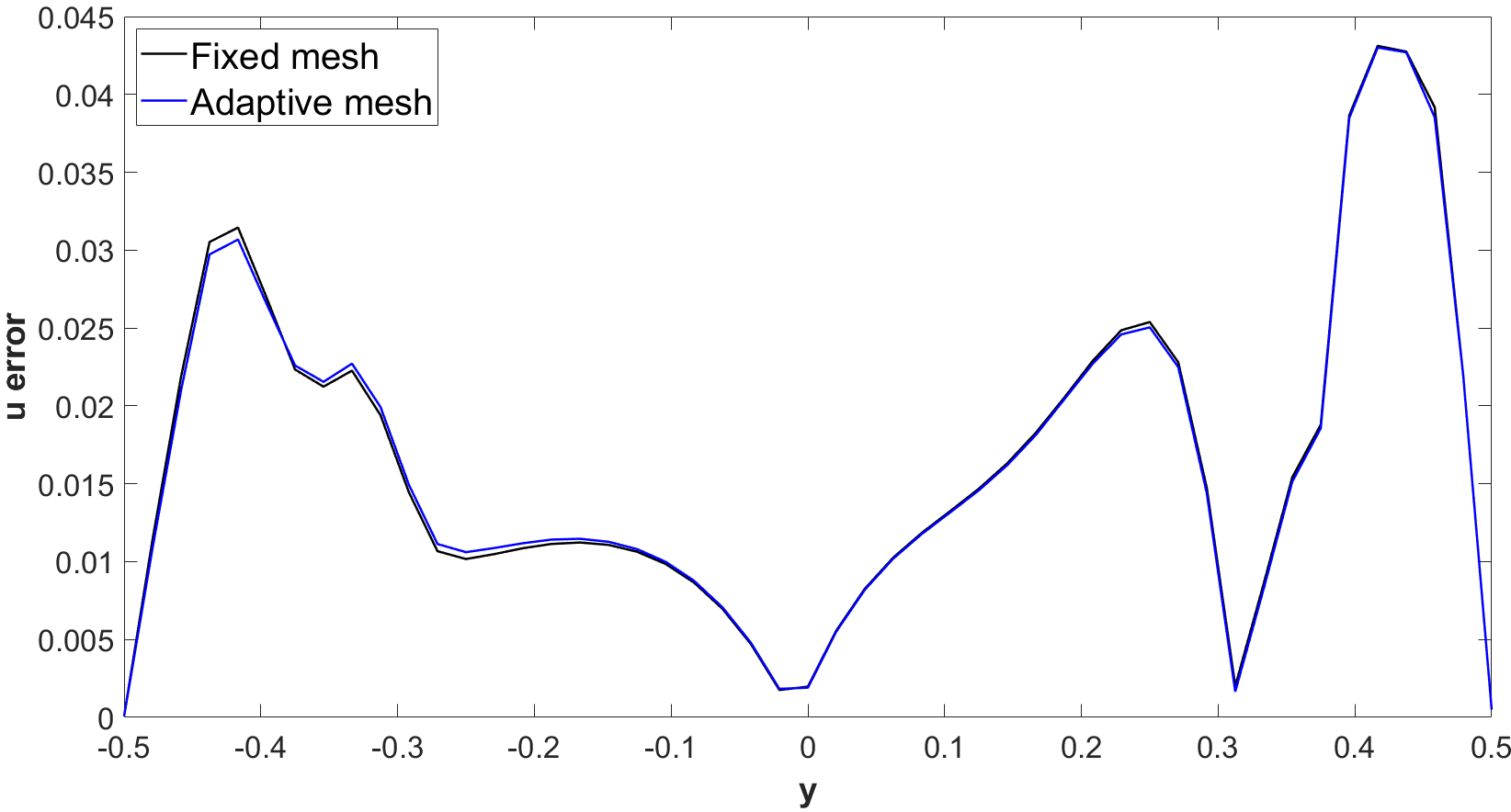} a)
  	\includegraphics[width=0.45\textwidth]{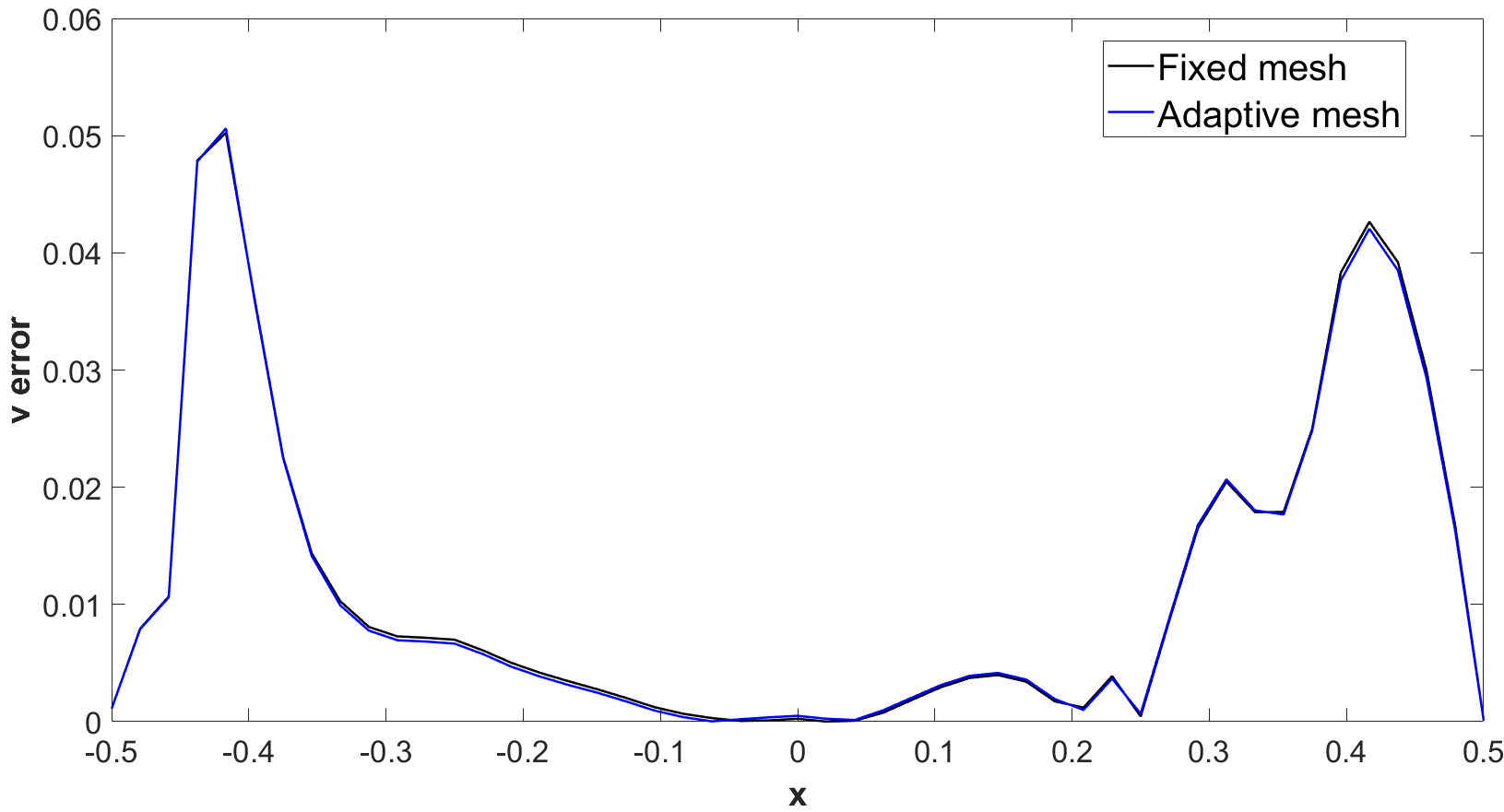} b)
  	\caption{3D lid driven cavity benchmark at \(Re = 1000\) , a) \(u\) velocity component component for the plane \(x=0,z=0\) with reference solution values from \cite{albensoeder:2005} interpolated, b) \(v\) velocity component comparison  for the plane \(y=0,z=0\) with reference solution values from \cite{albensoeder:2005} interpolated. The continuous black line denotes the result with fixed mesh, the blue one denotes the results with adaptive mesh refinement.}
  	\label{fig:velocity_error_comparison_3D_lid_driven}
\end{figure}

We have also exploited again the mesh adaptivity tool provided by \textit{deal.II} with the same refinement indicator introduced for the two-dimensional test. In particular, we started from a coarse mesh with \(N_e = 6\) elements along each direction and again we performed the refinement procedure on at most 10\% of the elements with the largest indicator value every 1000 time steps, while coarsening on at most 30\% of the elements with the smallest indicator values; moreover we have checked every 50 time steps if the refinement procedure had to be performed in advance in case the maximum difference between the velocities at two consecutive time steps was greater then \(10^{-2}\). The minimum element diameter allowed was \(\mathcal{H} = \frac{1}{48}\) in order to obtain \(C \approx 1\). In Figure \ref{fig:3dcav_u_y_adapt} and in Figure \ref{fig:3dcav_v_x_adapt} we report again the results for the \(v\) velocity component values along the $x$ axis and the \(u\) component of the velocity along the $y$ axis, respectively, compared with the results obtained using a fixed grid with \(N_e = 48\) elements along each direction. One can notice very good agreement between the two simulations, while the computational time required to perform the adaptive simulation is about half of that required by the fixed grid simulation. Moreover,  we have compared in Figure \ref{fig:velocity_error_comparison_3D_lid_driven} the errors of the two components for the velocity for the fixed and adaptive mesh, respectively. It is clear that, in spite of the different computational time, no significant differences arise.
  
\FloatBarrier

\begin{figure}[!h]
	\centering
	\includegraphics[width=0.8\textwidth]{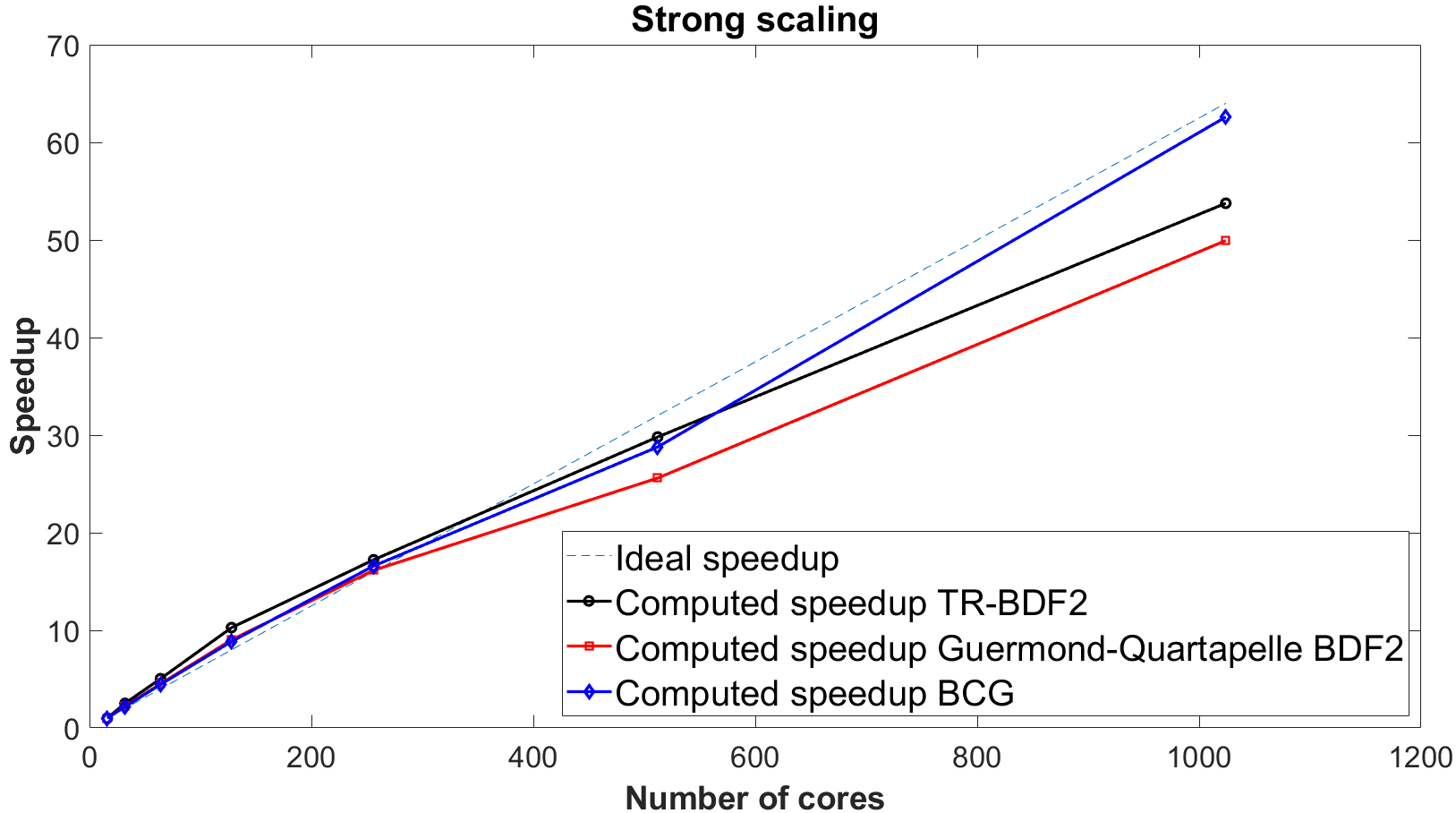} 
	\caption{3D lid driven cavity benchmark at \(Re = 1000\), strong scaling for the proposed method. The speedup is computed with respect to the time required with 16 cores.}
	\label{fig:strong_scaling}
\end{figure}

The size of the configuration employed for this test (we have used 15925248 dofs for the velocity and 1572864 dofs for the pressure) makes this benchmark a good candidate for a parallel scaling test. More specifically, we have performed a strong scaling analysis executing the same simulation up to time \(t = 0.6\) using from 16 up to 1024 2xCPU x86 Intel Xeon Platinum 8276-8276L @ 2.4Ghz cores of the HPC infrastructure GALILEO100 at the Italian supercomputing center CINECA. The results, reported in Figure \ref{fig:strong_scaling}, show a very good linear scaling, and even superlinear due to cache effects, up to 256 cores, while for a higher number of cores parallel performance is less optimal. A degradation of the performance for higher numbers of cores is observed, which we believe is mainly due to the fact that, given the size of the problem we were able to run, for these numbers of cores the amount of degrees of freedom owned by each core becomes very small so that the time needed by each core for computation   is dominated by the time needed for communication. Indeed, using 1024 cores, the number of unknowns per core is only 15552 for the velocity and 1536 for the pressure. Similar results are obtained for the Bell-Colella-Glaz and the Guermond-Quartapelle BDF2 methods already discussed in the previous Section. The apparent better behaviour of the Bell-Colella-Glaz projection scheme for a larger number of processors is probably due to the fact that this method, in view of the presence of nonlinear iterations, is the slowest one as reported in Table \ref{tab:scability_times} and therefore the least affected by communication costs.

\FloatBarrier

\begin{figure}[!h]
	\centering
	\includegraphics[width=0.8\textwidth]{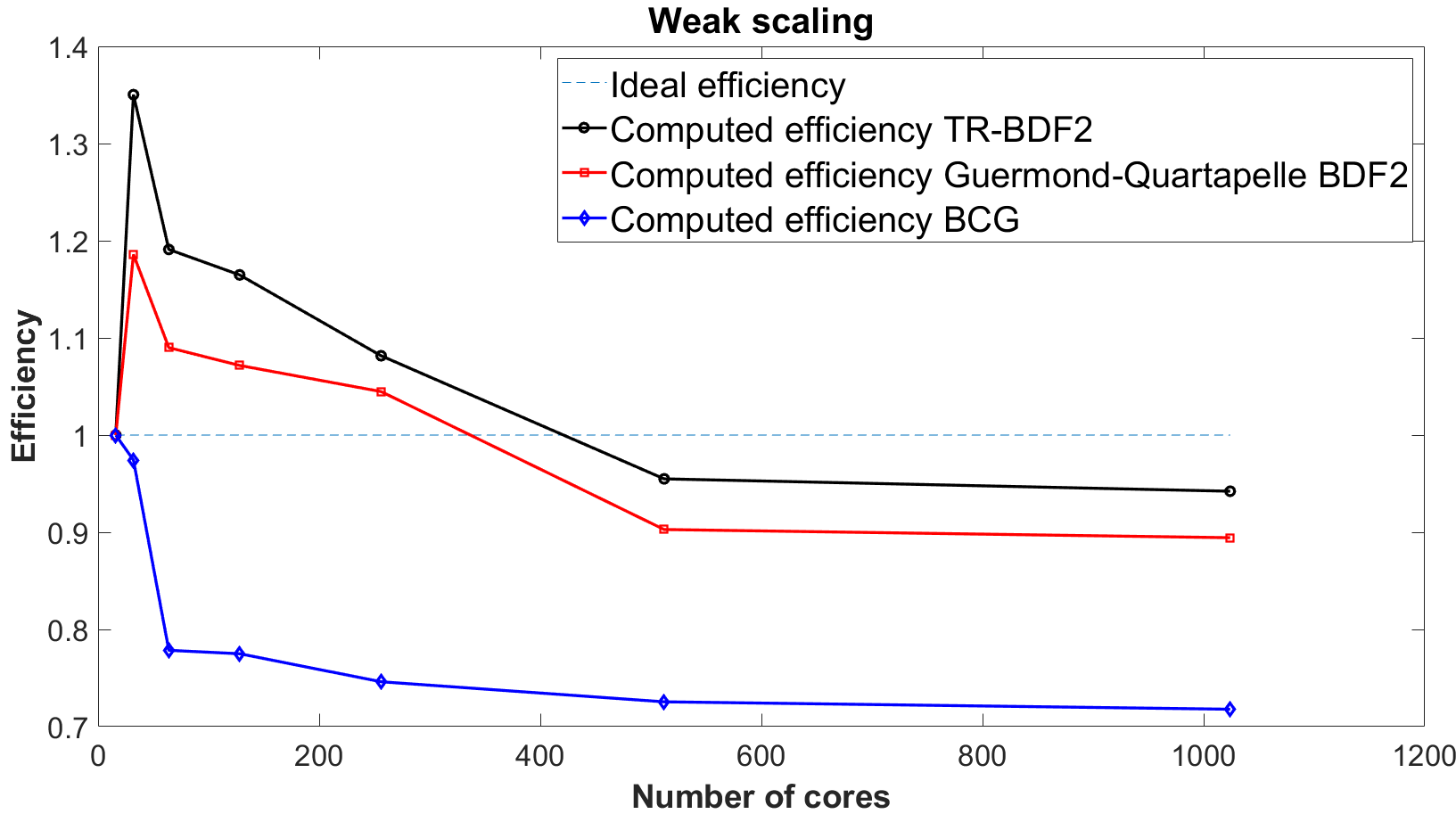}
	\caption{3D lid driven cavity benchmark at \(Re = 1000\), weak scaling for the proposed method. The efficiency is computed with respect to the time required with 16 cores.}
	\label{fig:weak_scaling}
\end{figure}

A weak scaling analysis has been performed using 124416 dofs per core for the velocity and 12288 dofs per core for the pressure. Figure \ref{fig:weak_scaling} shows the results for the three schemes described. One can easily notice that a good parallel efficiency is maintained up to 1024 cores. The overperformance of the TR-BDF2 scheme up to 256 cores can be due to a number of factors, such as the topology of the communication network in the specific architecture employed or the handling of communications between  different  groups of cores. 
This is also confirmed by the behaviour of the Bell-Colella-Glaz method which, since it requires the solution of more linear systems, is less dependent on these factors.

\FloatBarrier

\subsection{Flow past a cylinder}
\label{ssec:cylinder}

\begin{figure}[!h]
	\includegraphics[width=0.9\textwidth]{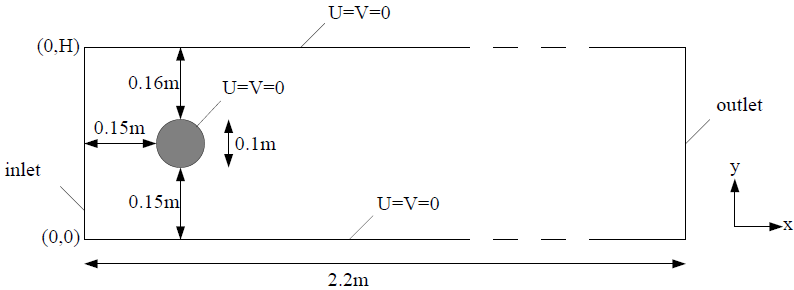} 
	\caption{Flow past a cylinder benchmark, geometry and boundary conditions (image from \cite{schafer:1996}).}
	\label{fig:cylinder_geometry}
\end{figure}

In this section, we consider another classical benchmark for the incompressible Navier-Stokes equations, namely the flow past a cylinder. We use the configuration described in \cite{schafer:1996}, that we summarize here for the reader's convenience. More in detail, the employed geometry and boundary conditions are reported in Figure \ref{fig:cylinder_geometry}, where \(H = \SI{0.41}{\meter}\) makes the domain non-symmetric and allows the vortex shedding in the wake of the cylinder. The inflow condition is

\[\mathbf{u}(0,y) = \left(\begin{array}{c}
									  4 U_m \frac{y \left(H - y\right)}{H^2} \\
									  0
									 \end{array}\right),\]

with \(U_{m} = \SI{1.5}{\meter\per\second}\). As explained in \cite{schafer:1996}, we consider as reference quantities the inflow velocity mean value \(U = \SI{1.0}{\meter\per\second}\), the diameter of the cylinder equal to \(L = \SI{0.1}{\meter}\) and \(\nu = \SI{0.001}{\square\meter\per\second}\), which yields \(Re = 100\).

\FloatBarrier

\begin{figure}[!h]
	\begin{center}
		\includegraphics[width=0.9\textwidth]{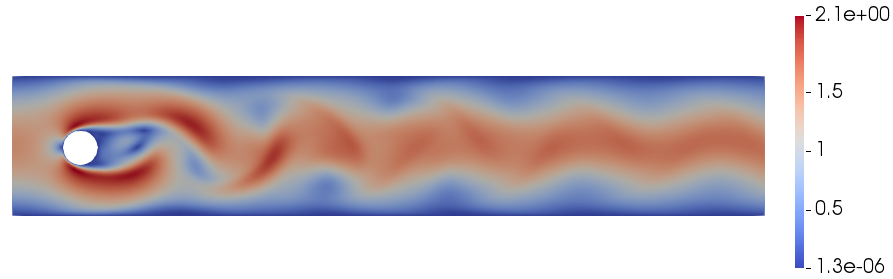}
	\end{center}
	\caption{Flow past a cylinder benchmark, contour plot of velocity magnitude at \(t = T\).}
	\label{fig:cylinder_wake}
\end{figure}

\begin{figure}[!h]
	\includegraphics[width=0.45\textwidth]{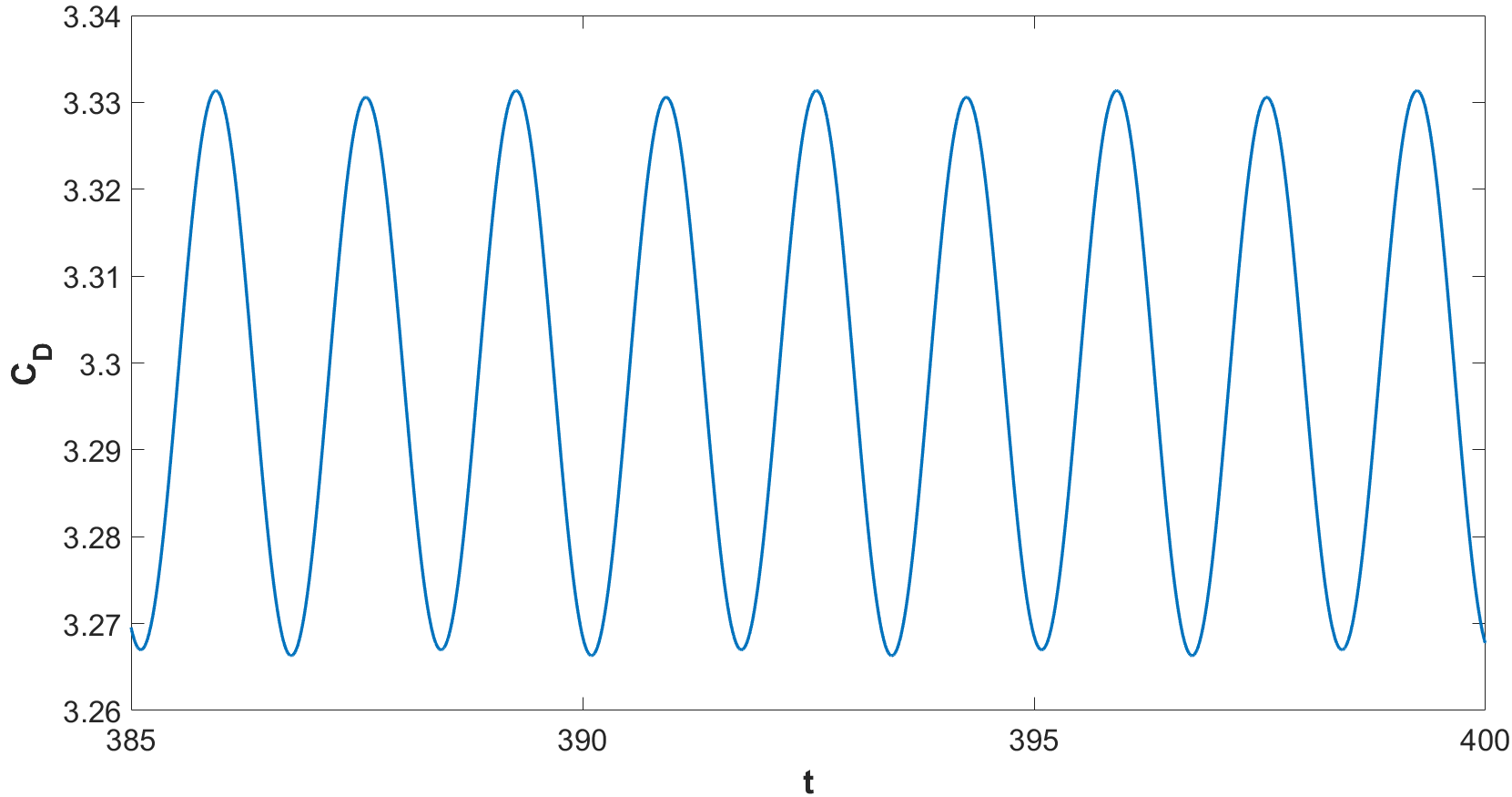} a)
	\includegraphics[width=0.45\textwidth]{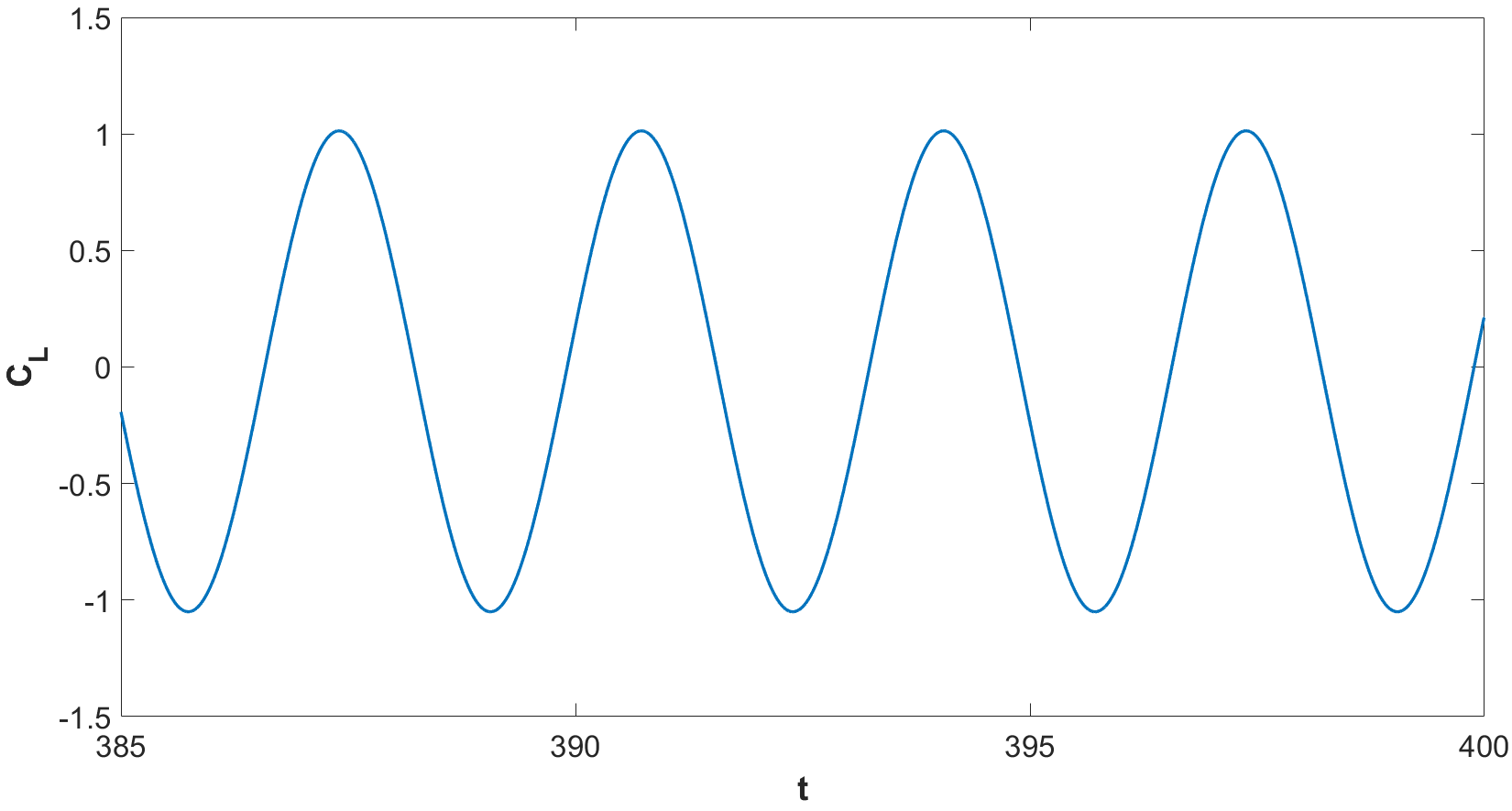} b)
	\caption{Flow past a cylinder benchmark, a) drag coefficient, b) lift coefficient.}
	\label{fig:cylinder_coefficients}
\end{figure}

We compute the drag and lift coefficients, defined as in \cite{schafer:1996}. Other reference values are the pressure drop \(\Delta p(t) = p(0.15, 0.2, t) - p(0.25, 0.2, t)\) and the Strouhal number \(St = \frac{Df}{U}\), where \(f\) is the frequency of separation computed as a function of the lift coefficient \(C_L\). The final time is \(T = 400\), which corresponds to a dimensional time of \(\SI{40}{\second}\), since the reference time value is \(\frac{L}{U} = \SI{0.1}{\second}\), and allows to obtain a fully developed wake. The grid is composed by 23552 elements and the time step \(\Delta t = 5 \cdot 10^{-3}\) is such that the maximum Courant number is around 1. Figure \ref{fig:cylinder_wake} shows the contour plot of the velocity magnitude at \(t = T\) and one can easily notice the formation of the vortices in the wake of the cylinder. Figure \ref{fig:cylinder_coefficients} reports the evolution of the lift and drag coefficients from \(t = 385\) to \(t = T\); it can be observed that the expected periodic behaviour is retrieved. The maximum drag coefficient and the pressure drop are \(3.33\) and \(2.60\), respectively, which are slightly larger values than  the intervals \(\left[3.22, 3.24\right]\) and \(\left[2.46, 2.50\right]\) proposed in \cite{schafer:1996}, even though they are in the overall range of the solutions proposed in the literature. The maximum lift coefficient is \(1.01\), which lies in the interval \(\left[0.99, 1.01\right]\) present in \cite{schafer:1996}, while the Strouhal number is equal to \(0.3\), which is again in the interval \(\left[0.295, 0.305\right]\) reported in \cite{schafer:1996}.

\FloatBarrier

\begin{figure}[!h]
	\includegraphics[width=0.9\textwidth]{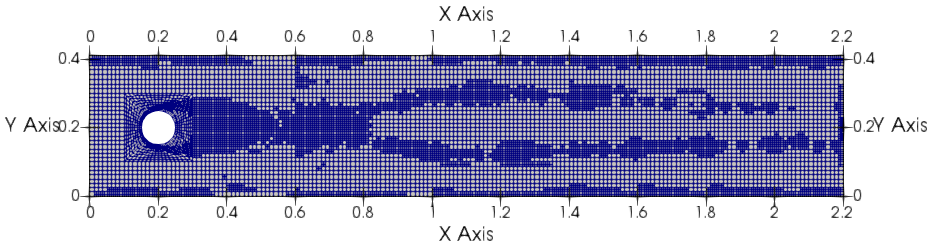} a)
	\includegraphics[width=0.45\textwidth]{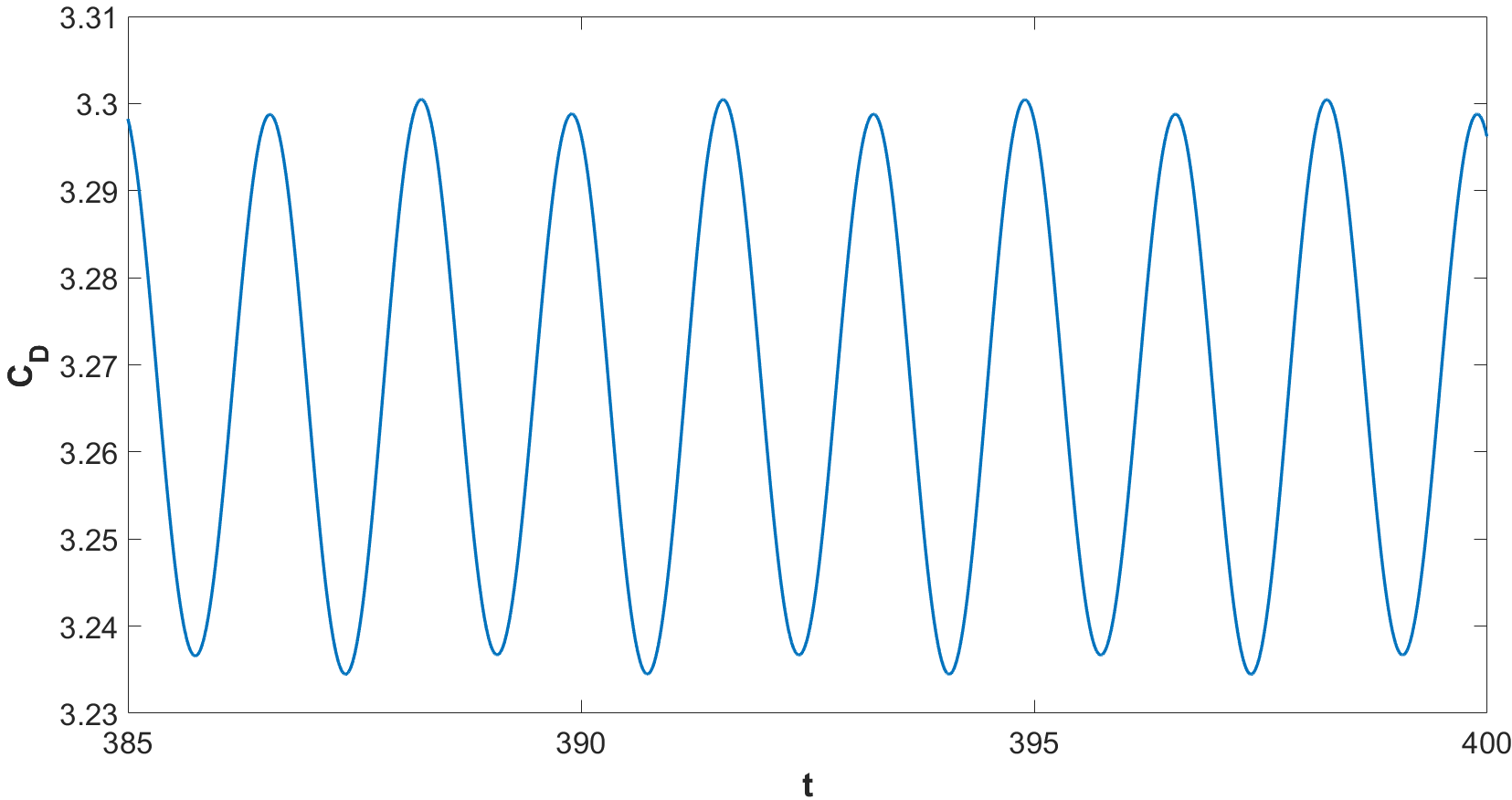} b)
	\includegraphics[width=0.45\textwidth]{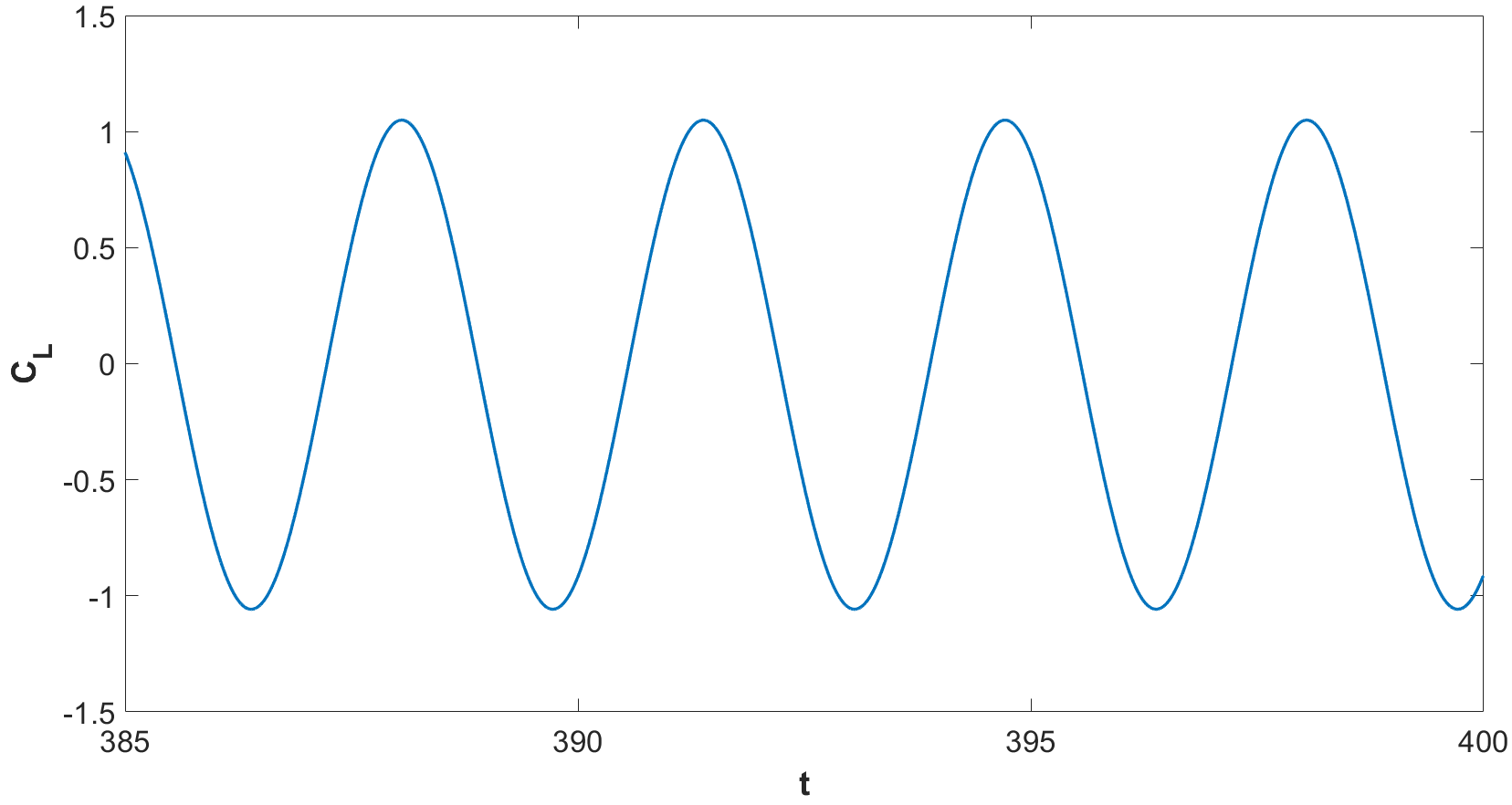} c)
	\caption{Flow past a cylinder benchmark, adaptive simulation, a) grid at \(t = T\), b) drag coefficient, c) lift coefficient.}
	\label{fig:cylinder_adaptive}
\end{figure}

The same test has been repeated using adaptive mesh refinement with the same criterion described in Section \ref{ssec:2dcavity}. The initial mesh is composed by 5558 elements and we allowed up to two local refinements, whereas the maximum element diameter is kept equal to the one of the initial grid. The same remeshing procedure described in \ref{ssec:2dcavity} was applied every 5000 time steps. Figure \ref{fig:cylinder_adaptive} reports the final mesh obtained and the values of the drag and lift coefficients. One can easily notice that more resolution is added in the wake of the cylinder and on its boundary and that the behaviour of the two coefficients is analogous to that in the  uniform mesh case. The final mesh consists of 11630 elements and a reduction of  computational  time of about 50\% is achieved with respect to the uniform mesh case.
 
\FloatBarrier  
 
\subsection{Complex geometry}
\label{ssec:compgeometry}

The matrix-free approach present in the \textit{deal.II} library makes the proposed solver attractive also for industrial applications that involve a large number of degrees of freedom. For this purpose we have tested the solver on the complex geometry of an heat exchanger of industrial interest \cite{burner:2018}. 
More specifically, a four channels module of a designed checkerboard pattern heat exchanger has been considered, with the goal of simulating its pure fluid-dynamic behaviour (i.e., in absence of heat exchanges) between the inlet and the outlet. The channel is long \(\SI{0.5}{\meter}\). 

\begin{figure}[!h]
	\includegraphics[width=0.8\textwidth]{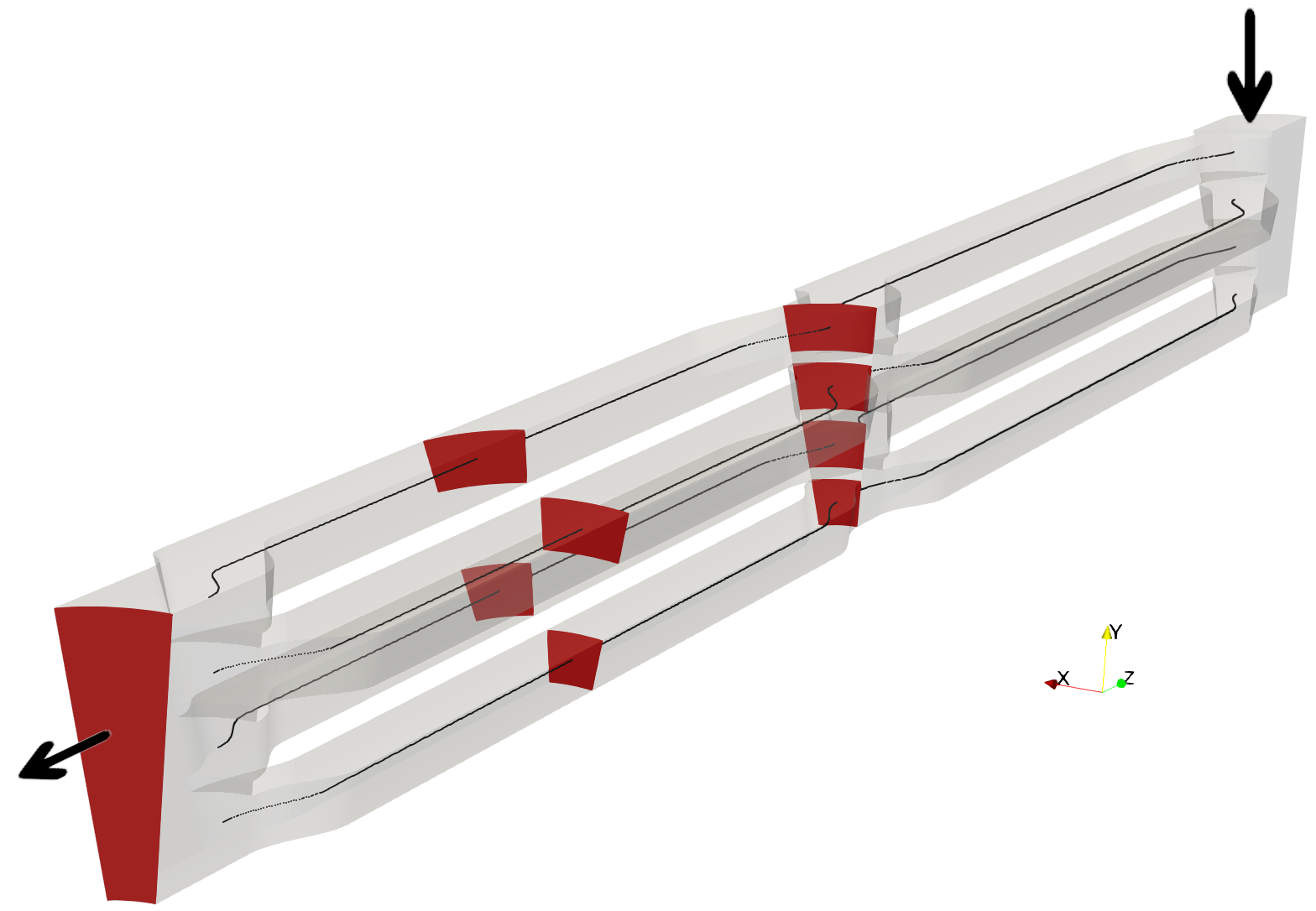} 
	\caption{Picture of the considered geometry} 
	\label{fig:compgeometry}
\end{figure} 
\FloatBarrier

\begin{figure}[!h]
 	\includegraphics[width=0.33\textwidth]{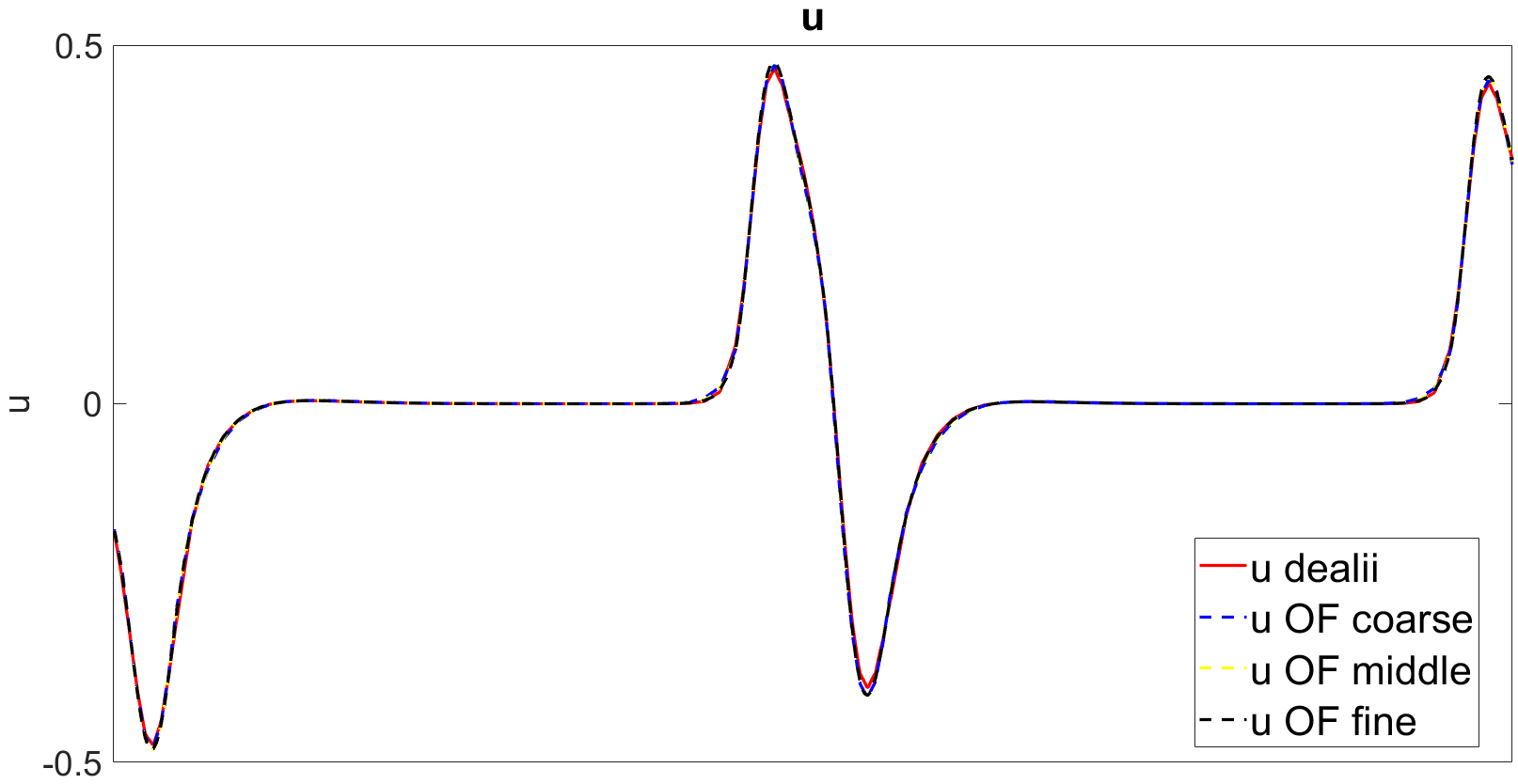} 
 	\includegraphics[width=0.33\textwidth]{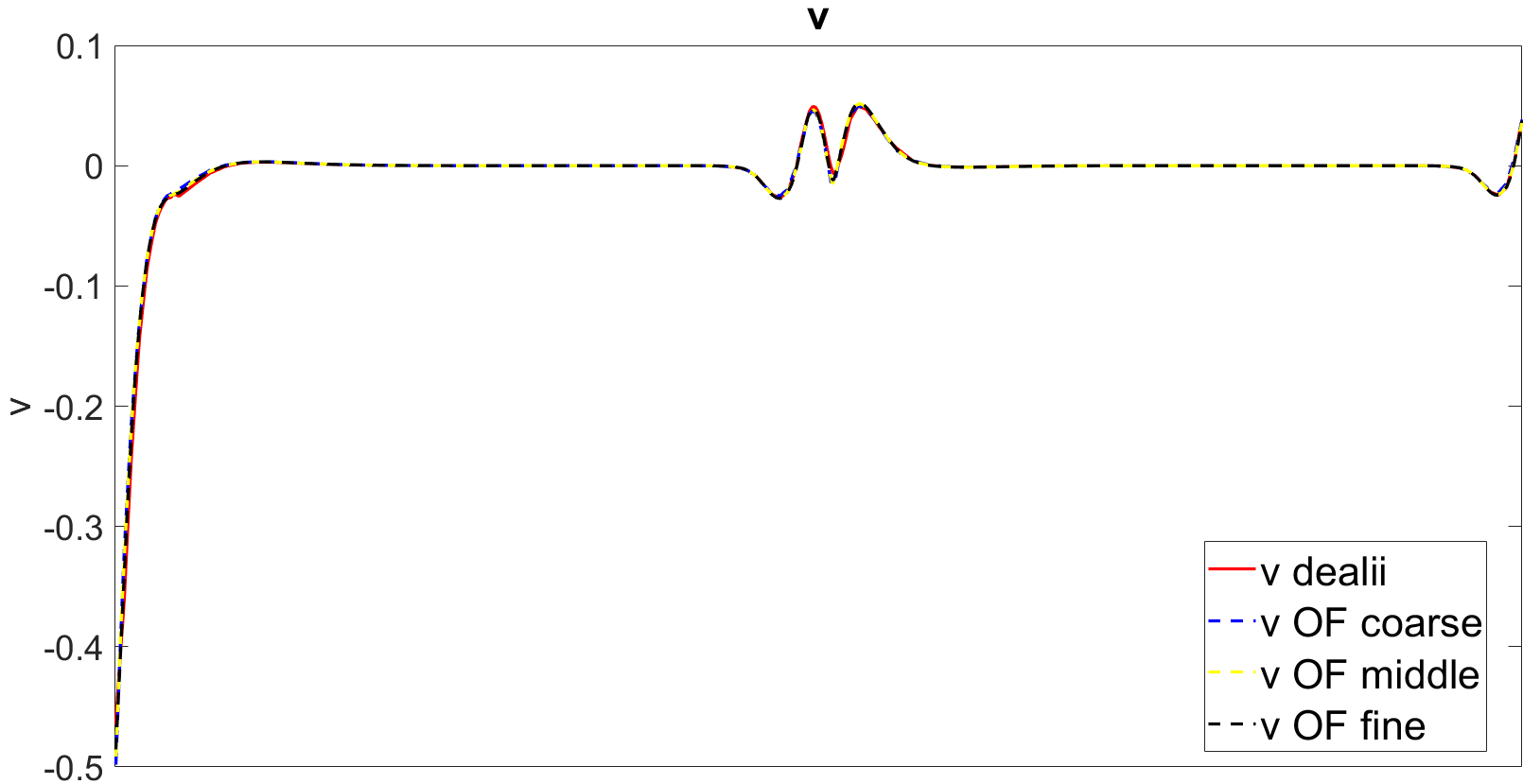} 
 	\includegraphics[width=0.33\textwidth]{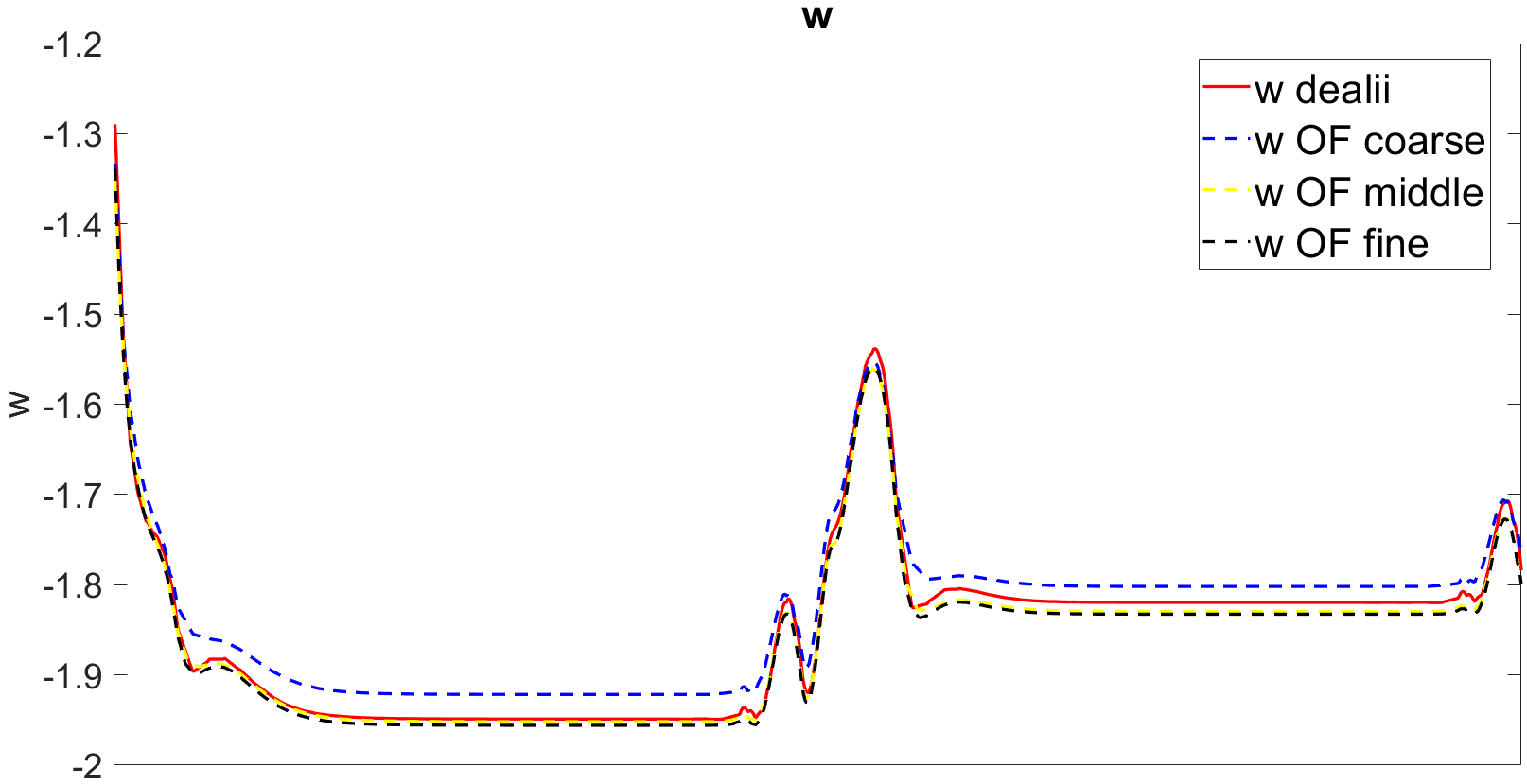} 
	\includegraphics[width=0.33\textwidth]{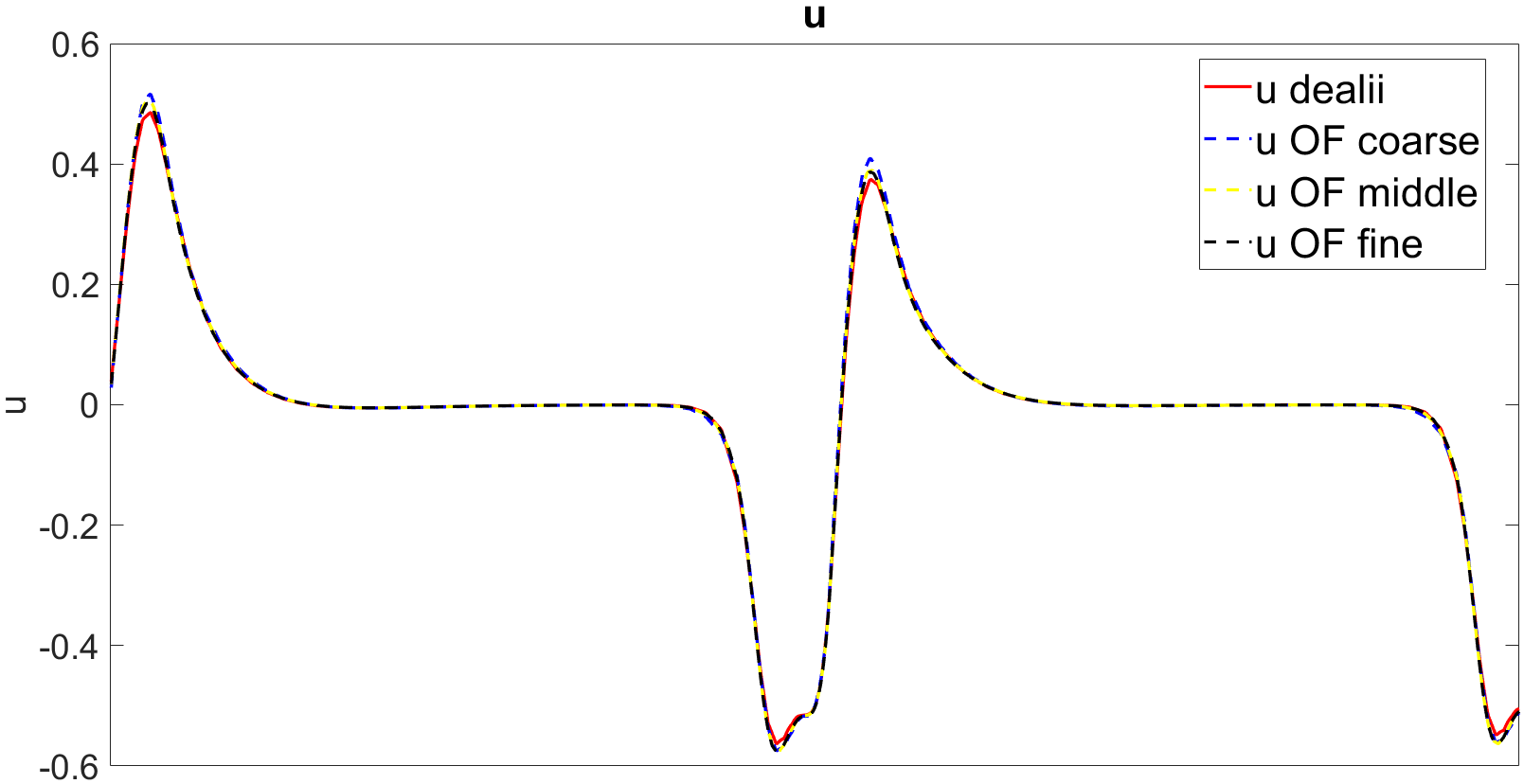} 
	\includegraphics[width=0.33\textwidth]{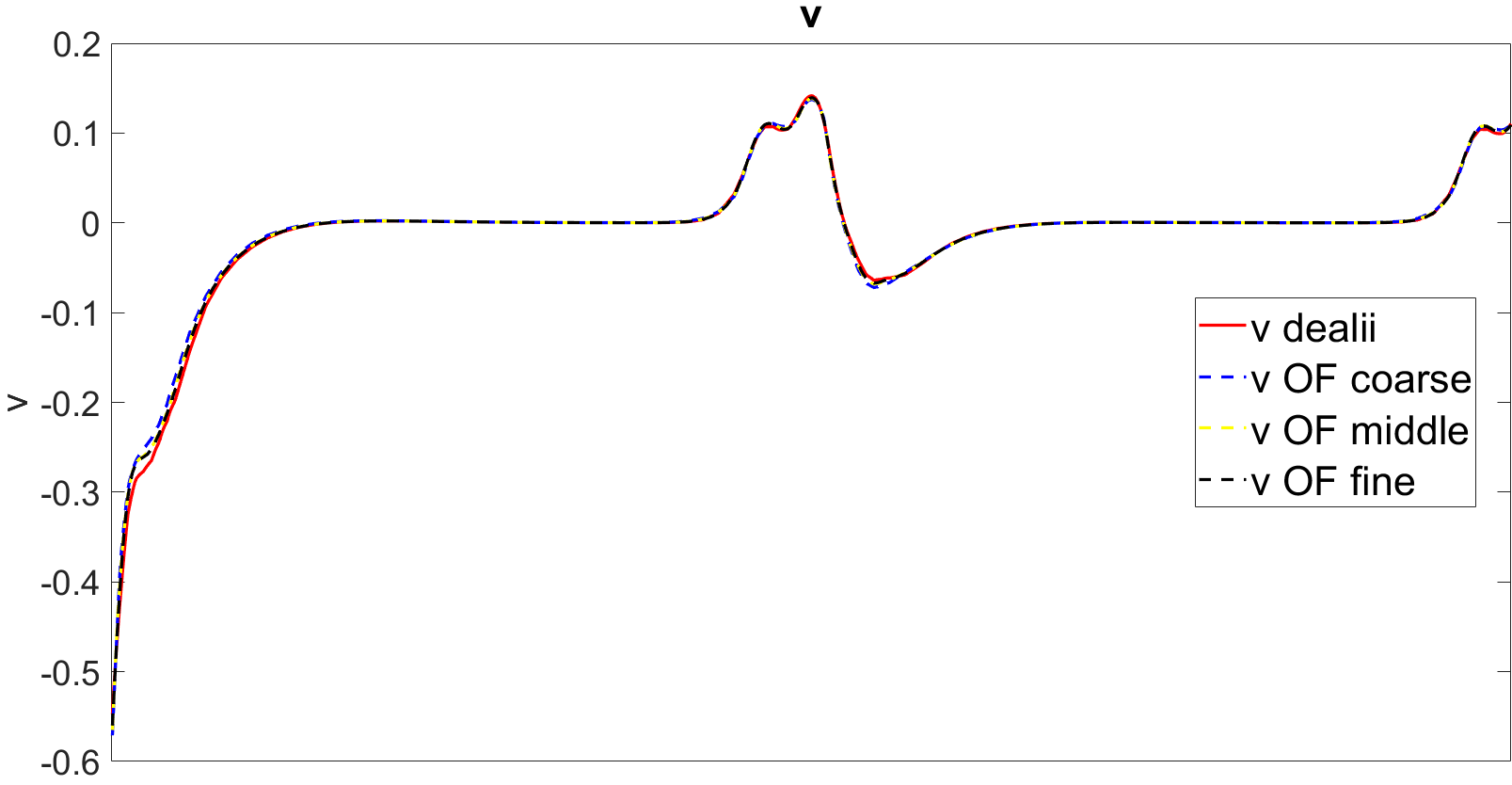} 
	\includegraphics[width=0.33\textwidth]{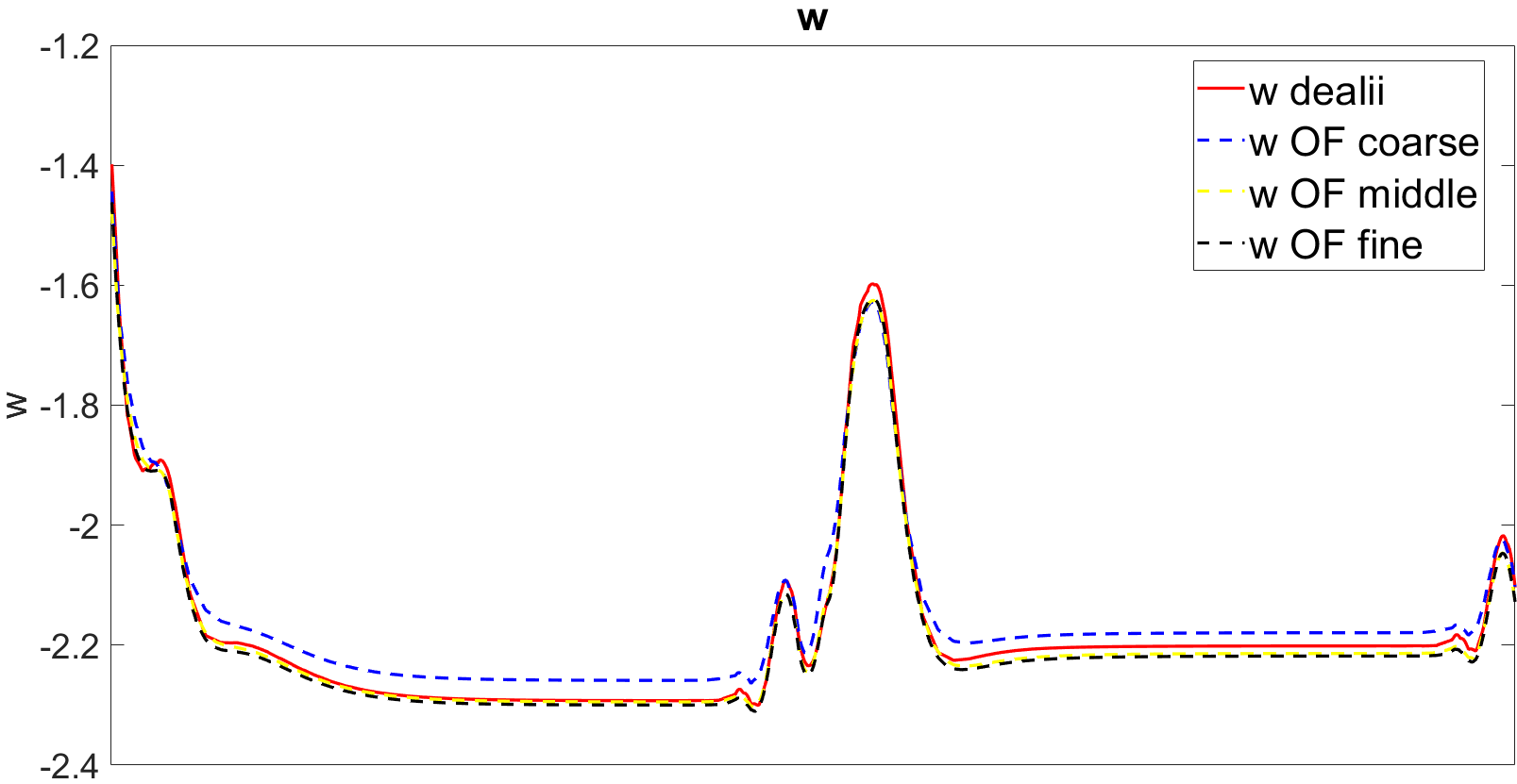} 
	\includegraphics[width=0.33\textwidth]{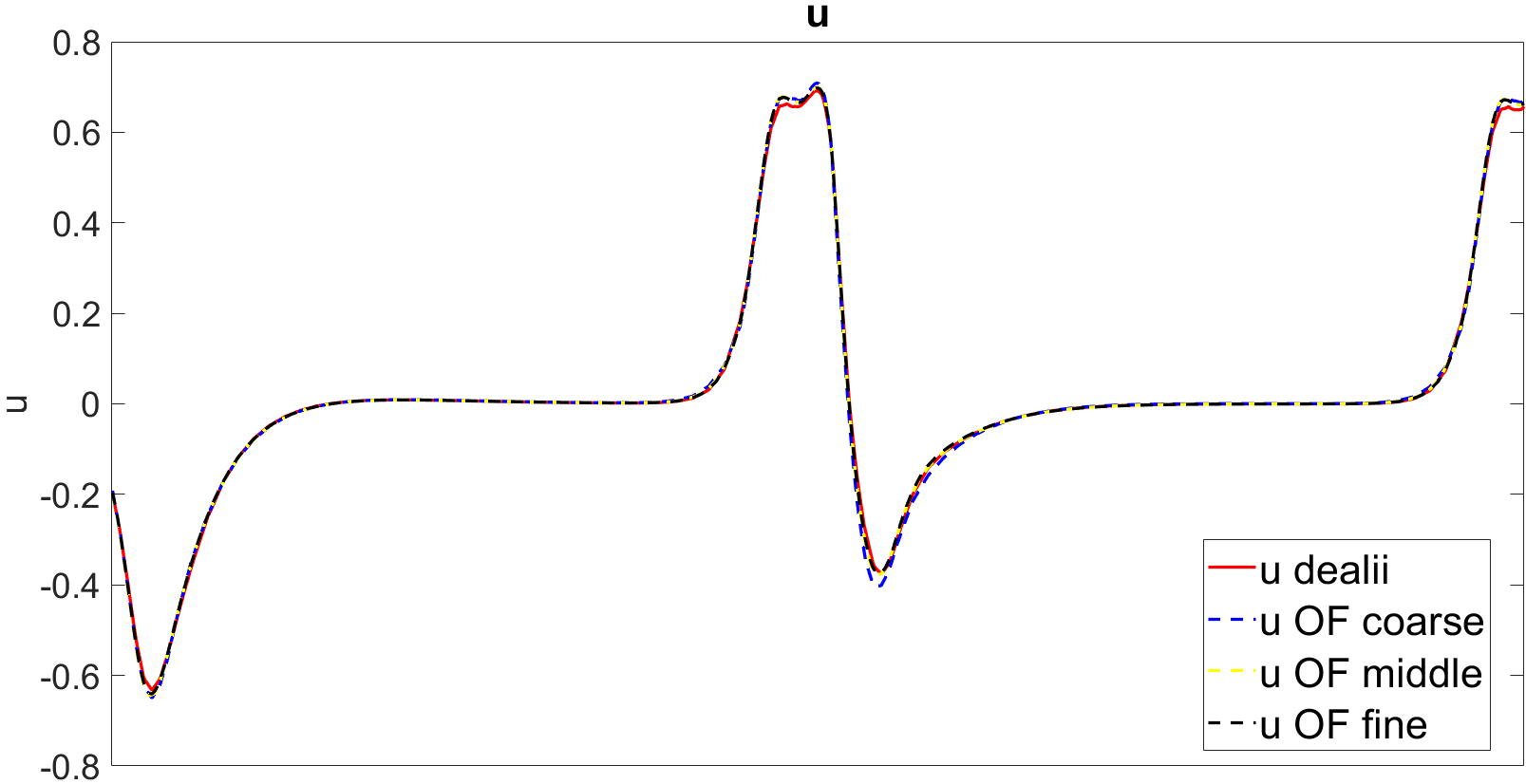} 
	\includegraphics[width=0.33\textwidth]{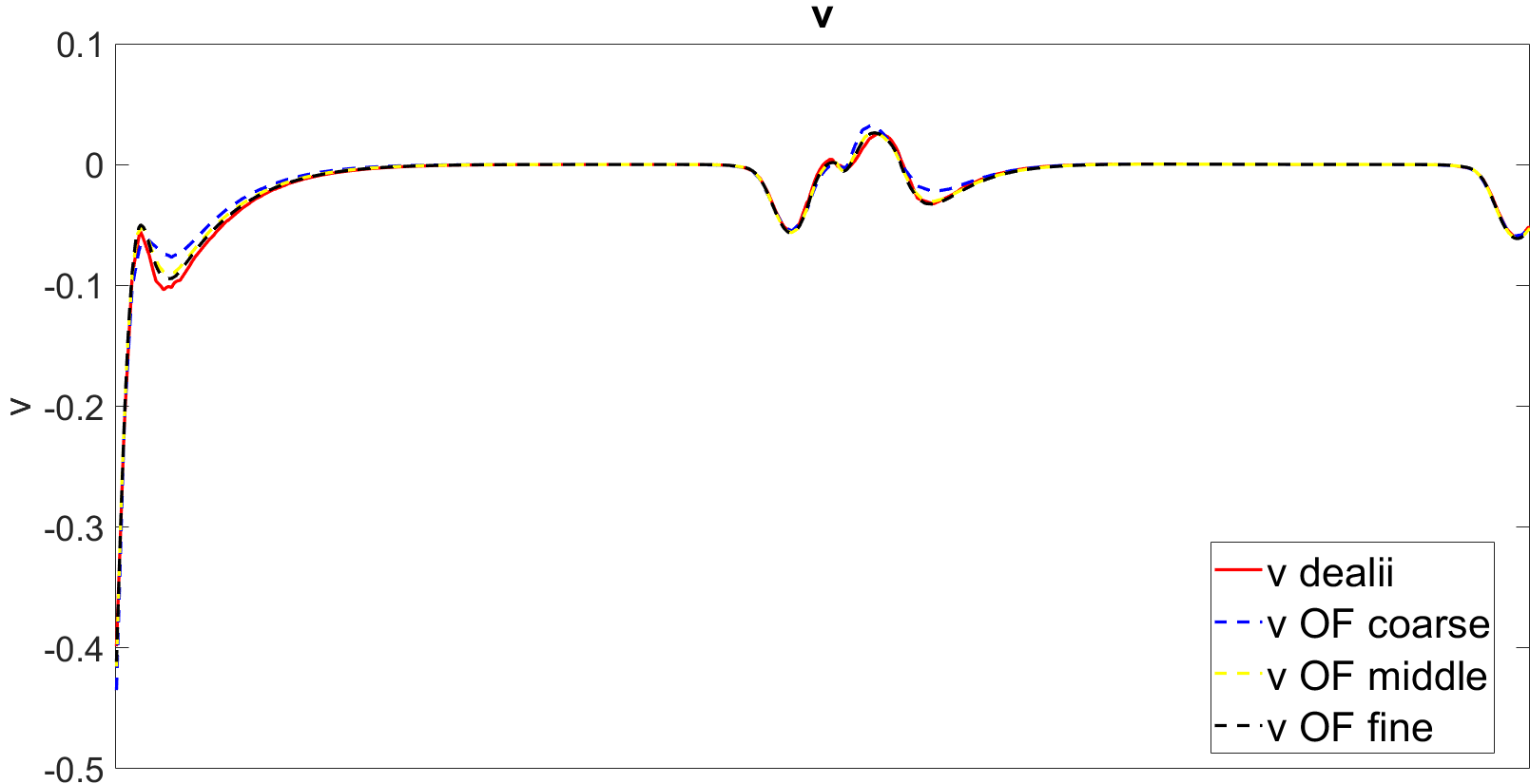} 
	\includegraphics[width=0.33\textwidth]{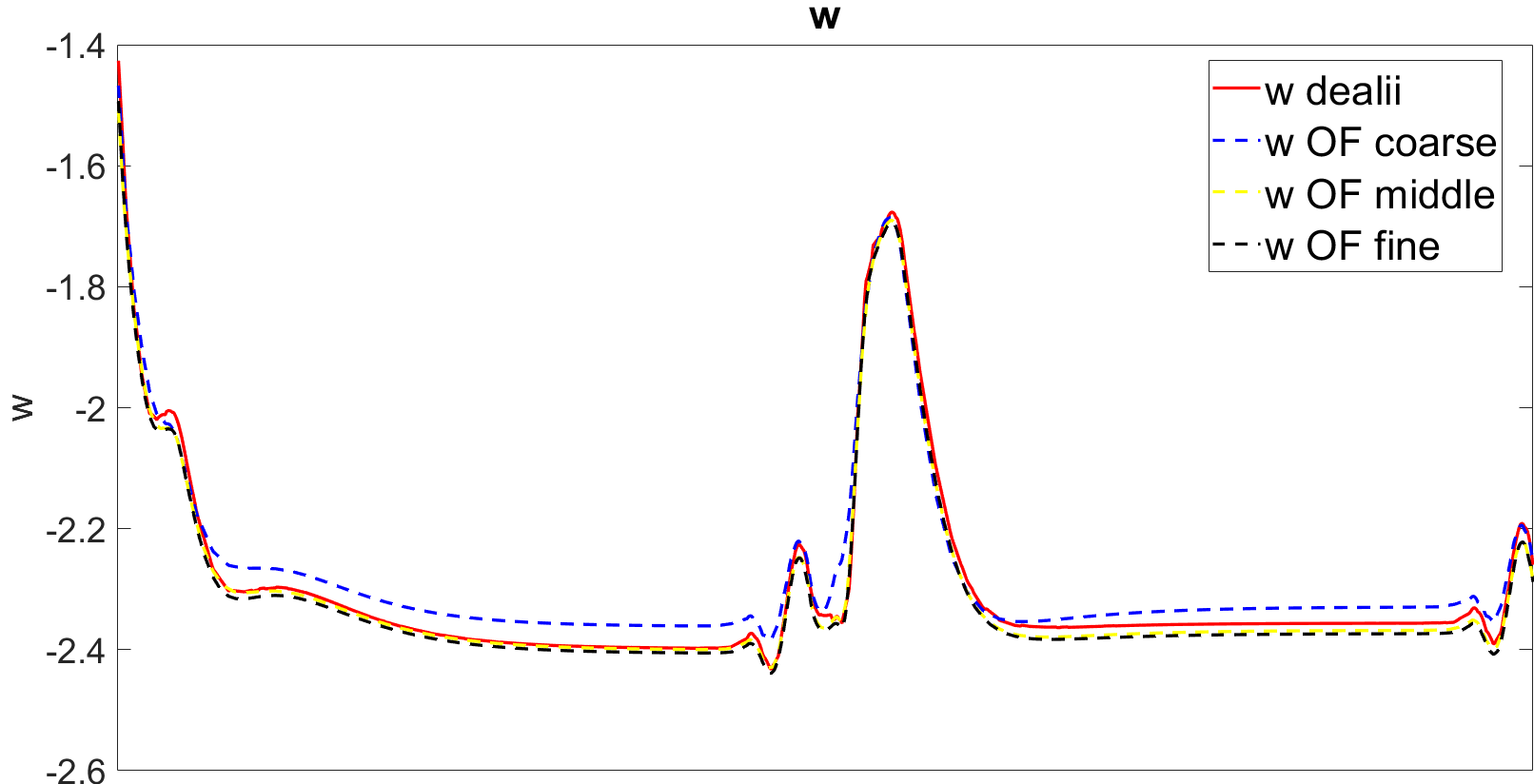} 
	\includegraphics[width=0.33\textwidth]{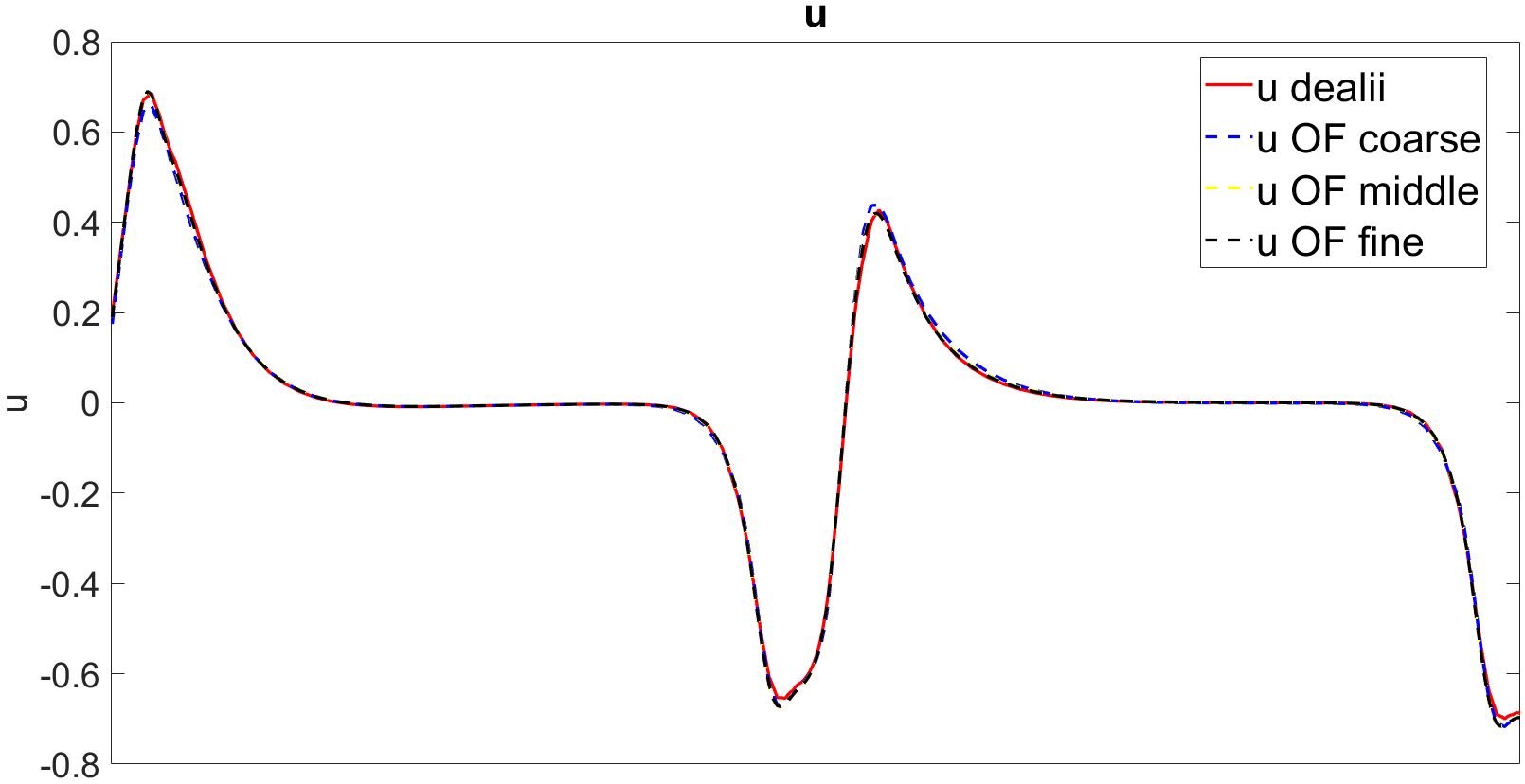} 
	\includegraphics[width=0.33\textwidth]{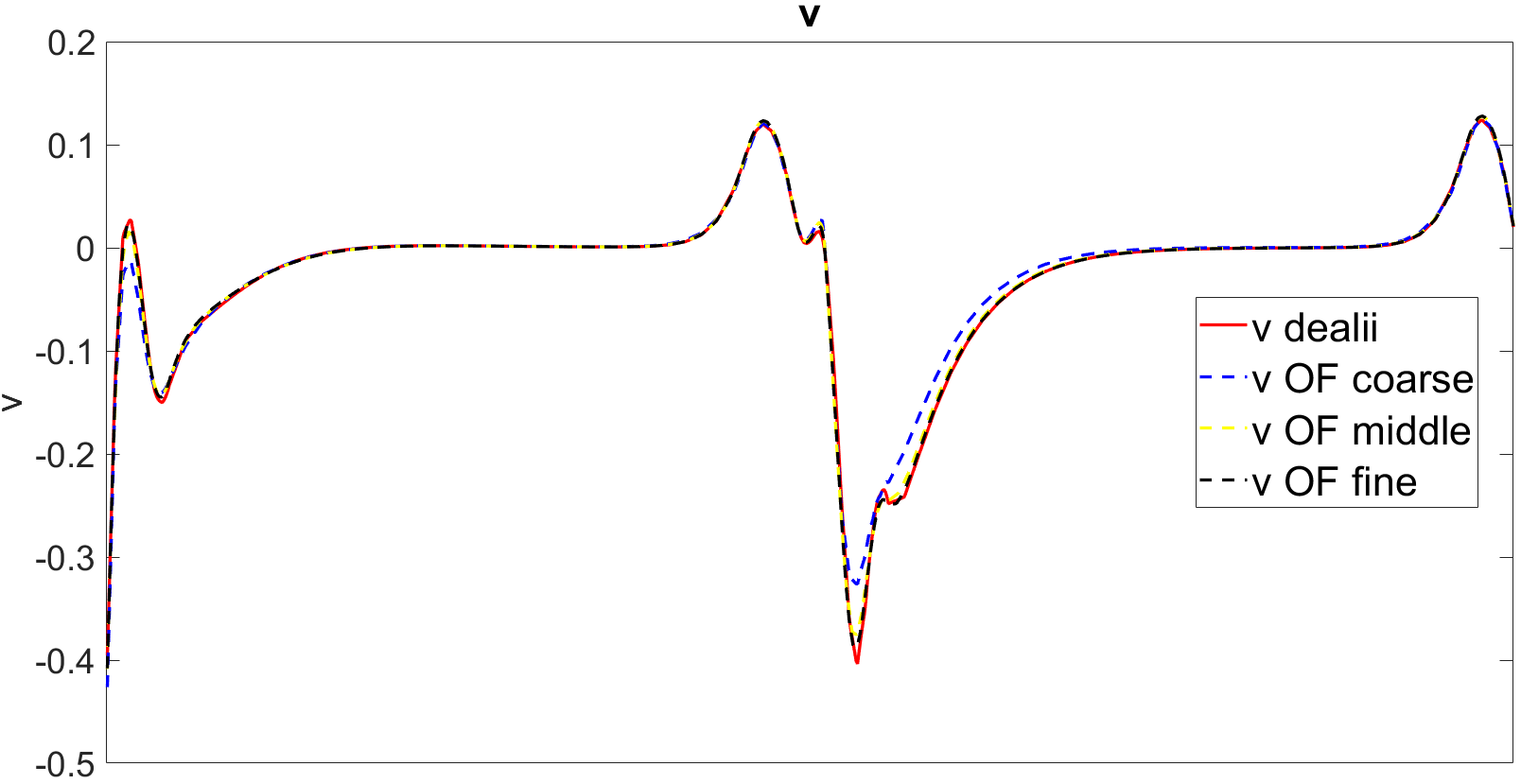} 
	\includegraphics[width=0.33\textwidth]{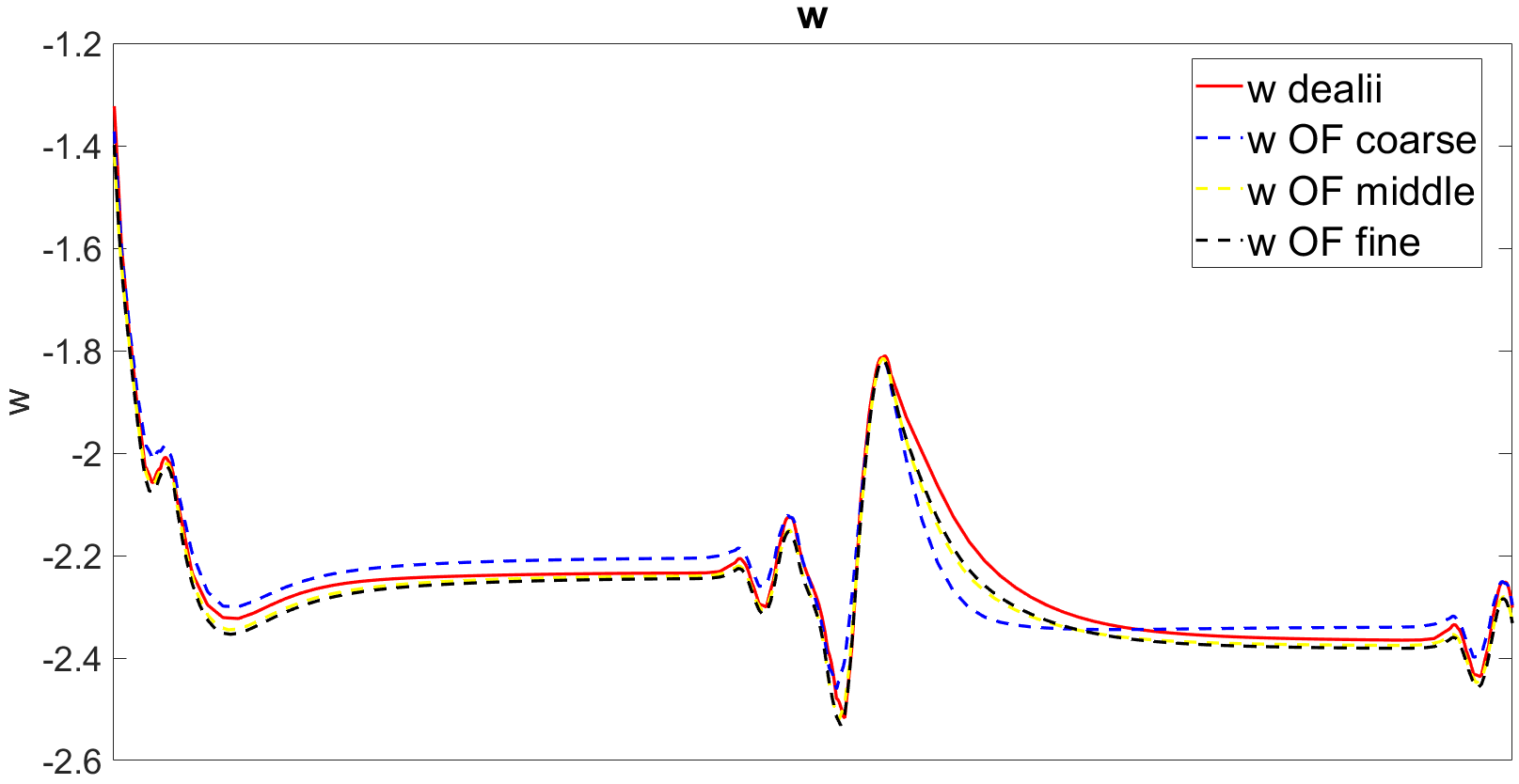} 
	\caption{Comparison between \textit{deal.II} and \textit{OpenFOAM} from line \(A\) (top) to line \(D\) (bottom) for \(u\) component of the velocity, the one along \(x\)-axis (left), \(v\) component of the velocity, the one along \(y\)-axis (center) and \(w\) component of the velocity, the one along \(z\)-axis (right)}
 	\label{fig:compgeometry_4lines}
\end{figure}

We consider \(Re = 5000\), assuming unitary inflow velocity considering the channel length as reference length. We set \(c = 300 \ m/s,\) which is of the order of magnitude of the speed of sound in air. The mesh consists of 129696 $\mathbf{Q}_2-Q_1$ elements, which yields \(10505376\) degrees of freedom for the discrete velocity variables and \(1037568\) for the discrete pressure variables. In order to verify the results of the simulation at steady state, various simulations with an OpenFoam steady state solver have been performed. More in detail, three meshes with different resolutions have been used with the OpenFoam solver.  The coarsest is the one previously described, an intermediate resolution one consists of 1382120 elements while the finest is composed by 2108119 elements. A comparison between the results obtained on each mesh is reported in Figure \ref{fig:compgeometry_4lines} for the midlines of the four channels depicted in Figure \ref{fig:compgeometry}. For the sake of simplicity, the channels are denote by \(A,B,C\) and \(D\) from bottom to top, respectively. 

\FloatBarrier
 
\begin{figure}[!h]
	\includegraphics[width=0.4\textwidth]{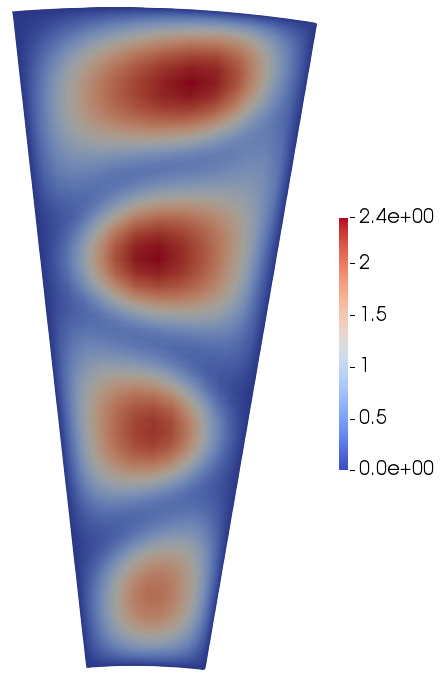} a) 
	\includegraphics[width=0.4\textwidth]{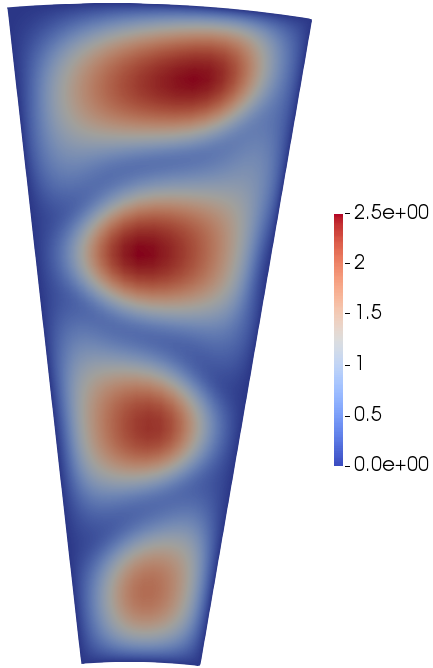} b) 
	\includegraphics[width=0.4\textwidth]{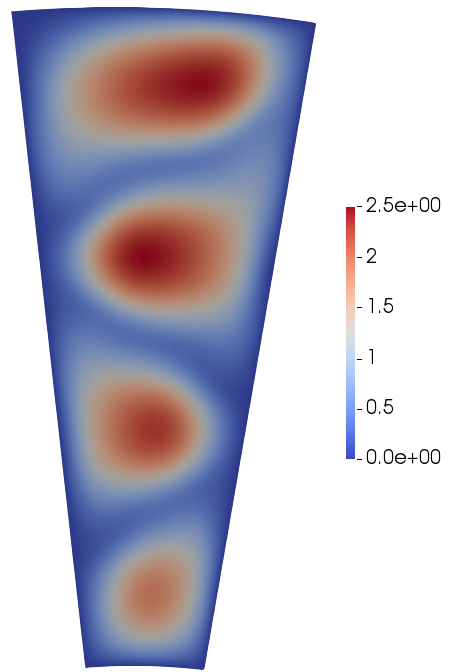} c)
	\hspace{2.5cm}
	\includegraphics[width=0.4\textwidth]{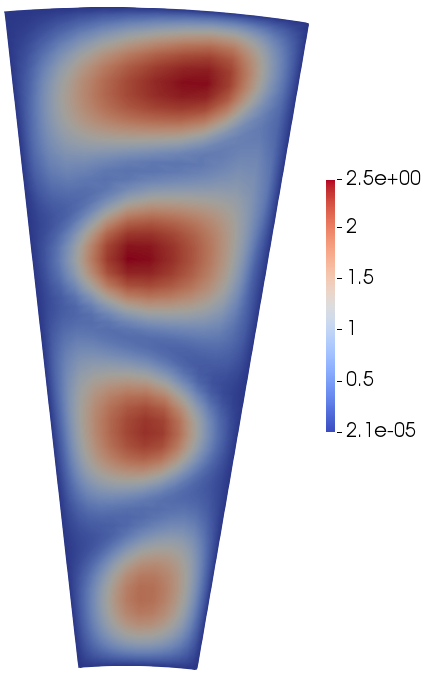} d)  
	\caption{Comparison between \textit{deal.II} and \textit{OpenFOAM} on the middle section, a) \textit{OpenFOAM} on coarse mesh, b) \textit{OpenFOAM} on middle mesh, c) \textit{OpenFOAM} on fine mesh, d) \textit{deal.II}}
	\label{fig:compgeometry_z022}
\end{figure}

\begin{figure}[!h]
	\includegraphics[width=0.45\textwidth]{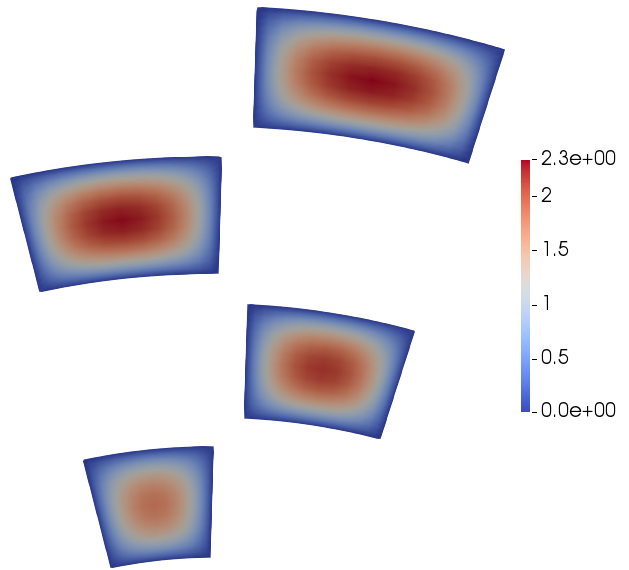} a) 
	\includegraphics[width=0.45\textwidth]{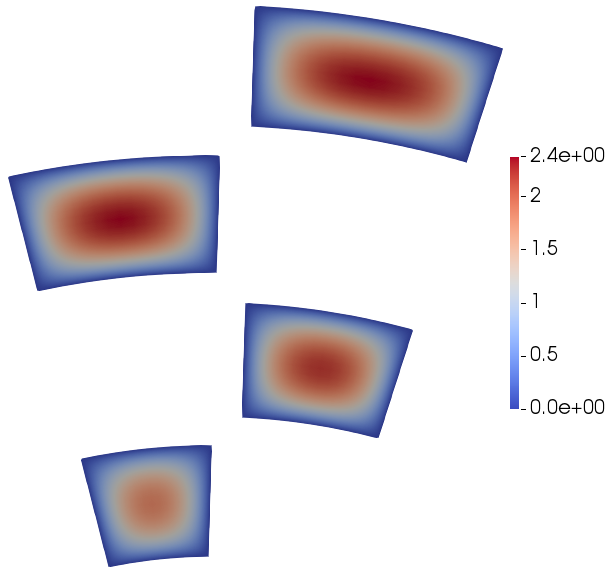} b)
	 
	\vspace{2cm}
	\includegraphics[width=0.45\textwidth]{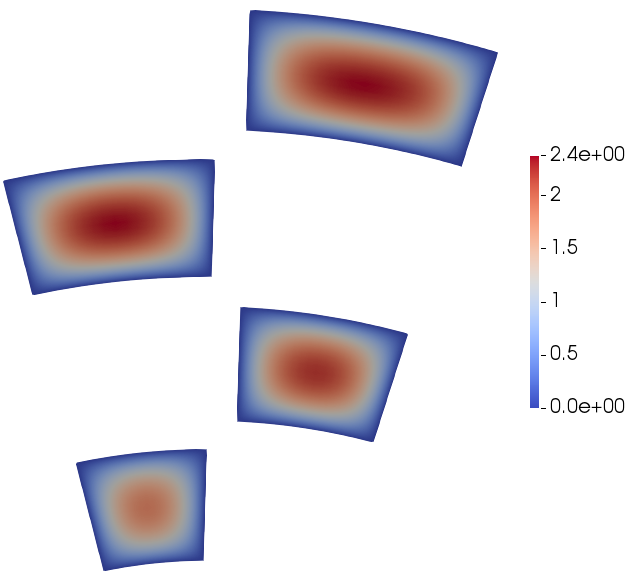} c) \hspace{0.5cm}
	\includegraphics[width=0.45\textwidth]{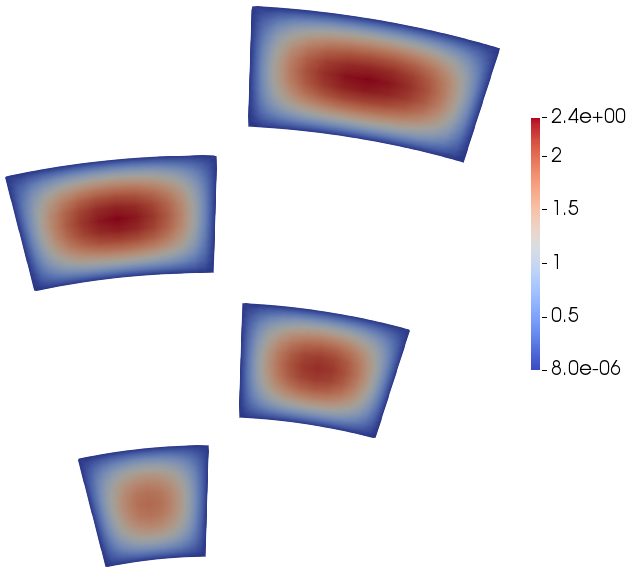} d)  
	\caption{Comparison between \textit{deal.II} and \textit{OpenFOAM} at three-quarters section, a) \textit{OpenFOAM} on coarse mesh, b) \textit{OpenFOAM} on middle mesh, c) \textit{OpenFOAM} on fine mesh, d) \textit{deal.II}}
	\label{fig:compgeometry_z010}
\end{figure}

\begin{figure}[!h]
	\includegraphics[width=0.4\textwidth]{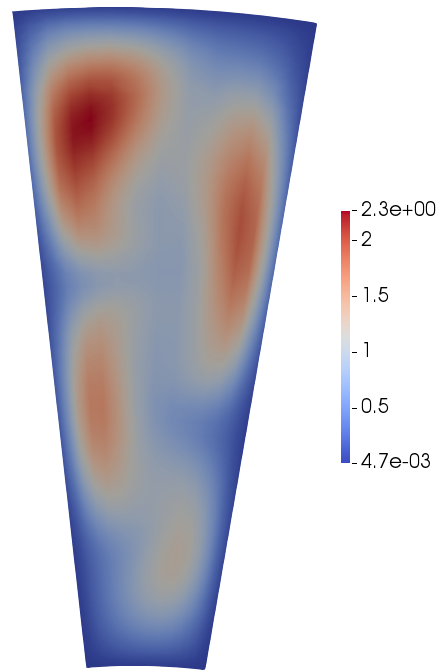} a) 
	\includegraphics[width=0.4\textwidth]{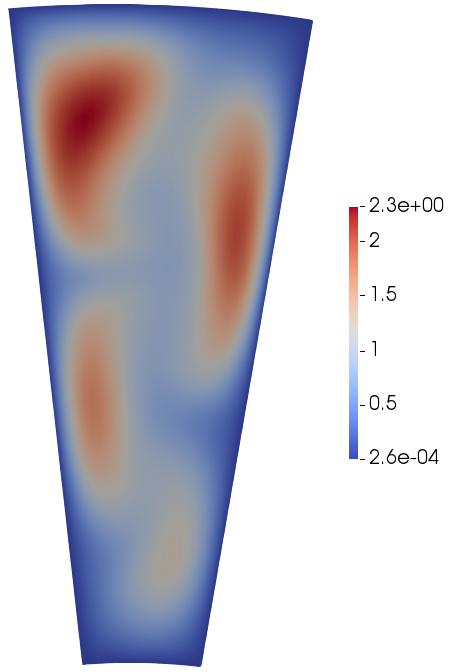} b) 
	\includegraphics[width=0.4\textwidth]{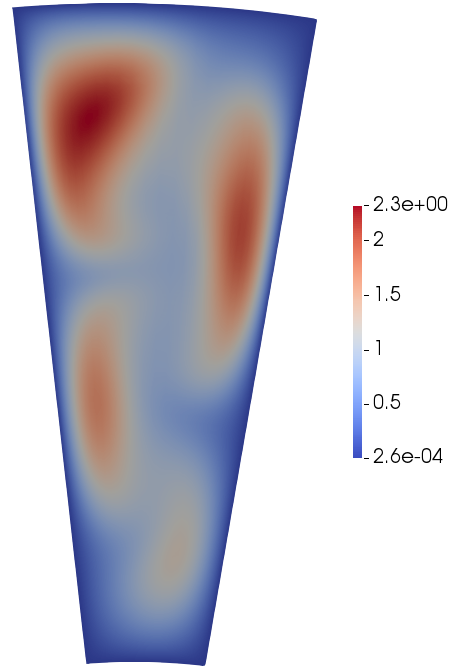} c)
	\hspace{2.5cm}
	\includegraphics[width=0.4\textwidth]{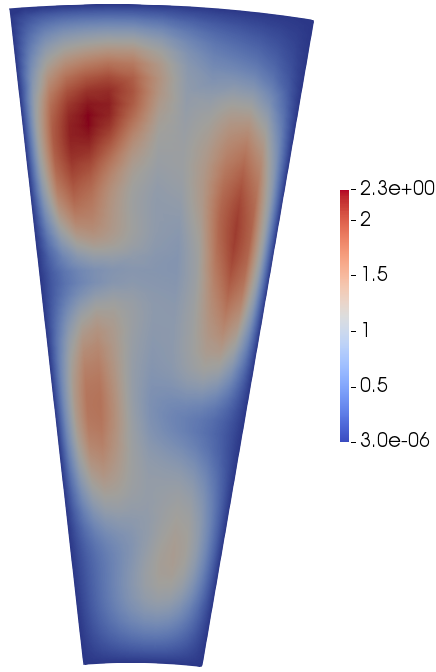} d)  
	\caption{Comparison between \textit{deal.II} and \textit{OpenFOAM} on the outlet, a) \textit{OpenFOAM} on coarse mesh, b) \textit{OpenFOAM} on middle mesh, c) \textit{OpenFOAM} on fine mesh, d) \textit{deal.II}}
	\label{fig:compgeometry_zoutlet}
\end{figure}

It can be observed that a good quantitative agreement between the two solvers has been obtained, taking into account the different features. Moreover, the solution computed with the DG approach is more similar to the results obtained with the OpenFoam solver on the finest meshes, as evident especially for the axial component \(w\). This is further confirmed by the pressure drop computed for the four lines and reported in Table \ref{tab:exchanger_pressure_drop}. Analogous considerations hold for the sections reported in Figure \ref{fig:compgeometry} where we have compared the contour of the velocity magnitude on the middle of the domain, at three-quarters of the domain and on the outlet in Figure \ref{fig:compgeometry_z022}, Figure \ref{fig:compgeometry_z010} and Figure \ref{fig:compgeometry_zoutlet}, respectively.

\FloatBarrier

\FloatBarrier

\FloatBarrier

\section{Conclusions and future perspectives}
\label{sec:conclu}

Building on the experience of \cite{dellarocca:2018}, we have proposed an accurate, efficient and robust projection method, based on the TR-BDF2 method. While time discretizations of the incompressible Navier-Stokes equations based on accurate implicit solvers have been proposed in many other papers, the specific combination of techniques presented in this work appears to be optimal under several viewpoints for the development of a second order adaptive flow solver.
    
The proposed fully implicit method has been implemented using discontinuous finite elements  in the framework of the numerical library \textit{deal.II}, with the aim of building a reliable, flexible and easily accessible tool for industrial applications that can ultimately be competitive with more conventional finite volume techniques. We have shown that the method has superior accuracy and efficiency with respect to some well known alternative schemes on a number of classical benchmarks.

In future work, besides application of the proposed approach to significant industrial applications and extensions to fully compressible and multiphase flow, an interesting development will be represented by the integration of more sophisticated \textit{a posteriori} error estimation techniques \cite{georgoulis:2009, hartmann:2006, oden:1994} to obtain optimal adaptive approaches. Furthermore, the multirate time integration version of the TR-BDF2 method \cite{bonaventura:2020a} could also be integrated in the discretization approach, so as to obtain a fully space-time adaptive technique based on a robust and unconditionally stable method.

\section*{Acknowledgements}
We would like to thank the two anonymous reviewers for their useful and productive comments, which helped us improving the initial version of this paper.
L.B. would also like to thank Roberto Ferretti, Elisabetta Carlini and Macarena G\'omez M\'armol for several useful discussions on numerical methods for incompressible flows. The parallel computations for the scaling test were performed at CINECA thanks to the computational resources made available through the SIDICoNS - HP10CLPLXI ISCRA C project.

\nocite{*}
\bibliography{DG_NS_IJNMF.bib}

\begin{thebibliography}{10}
\providecommand \doibase [0]{http://dx.doi.org/}%

\bibitem{quartapelle:2013}
Quartapelle L. {\it Numerical solution of the incompressible {N}avier-{S}tokes
  equations}.
\newblock Birkh{\"a}user .
\newblock 2013.

\bibitem{quarteroni:2008}
Quarteroni A, Valli A. {\it Numerical approximation of partial differential
  equations}. 23.
\newblock Springer Science \& Business Media .
\newblock 2008.

\bibitem{chorin:1968}
Chorin A. Numerical solution of the {N}avier-{S}tokes equations. {\it
  Mathematics of Computation} 1968\string; 22\string: 745--762.

\bibitem{temam:1969}
Temam R. Sur l'approximation de la solution des {\'e}quations de
  {N}avier-{S}tokes par la m{\'e}thode des pas fractionnaires (II). {\it
  Archive for Rational Mechanics and Analysis} 1969\string; 33\string:
  377--385.

\bibitem{guermond:2006}
Guermond J, Minev P, Shen J. An overview of projection methods for
  incompressible flows. {\it Computer methods in applied mechanics and
  engineering} 2006\string; 195\string: 6011--6045.

\bibitem{fletcher:1997}
Fletcher C. {\it Computational techniques for fluid dynamics, Volume 1:
  Fundamental and general techniques}.
\newblock Springer Verlag .
\newblock 1997.

\bibitem{chen:2014}
Chen G, Xiong Q, Morris P, Paterson E, Sergeev A, Wang Y. Open{FOAM} for
  computational fluid dynamics. {\it Notices of the AMS} 2014\string;
  61\string: 354--363.

\bibitem{jasak:2007}
Jasak H, Jemcov A, Tukovic Z. Open{FOAM: A C}++ library for complex physics
  simulations. In:  {\it International Workshop on Coupled Methods in Numerical
  Dynamics,}. 1000. IUC Dubrovnik Croatia. ; 2007\string: 1--20.

\bibitem{weller:1998}
Weller H, Tabor G, Jasak H, Fureby C. A tensorial approach to computational
  continuum mechanics using object-oriented techniques. {\it Computers in
  Physics} 1998\string; 12\string: 620--631.

\bibitem{dellarocca:2018}
Della~Rocca A. {\it Large-Eddy Simulations of Turbulent Reacting Flows with
  Industrial Applications}. PhD thesis. Politecnico di Milano,  2018.

\bibitem{giraldo:2020}
Giraldo F. {\it An Introduction to Element-Based {G}alerkin Methods on
  Tensor-Product Bases}.
\newblock Springer {N}ature .
\newblock 2020.

\bibitem{karniadakis:2005}
Karniadakis G, Sherwin S. {\it Spectral $hp-${E}lement {M}ethods for
  {C}omputational {F}luid {D}ynamics}.
\newblock {O}xford {U}niversity {P}ress .
\newblock 2005.

\bibitem{bassi:2005}
Bassi F, Crivellini A, Rebay S, Savini M. Discontinuous Galerkin solution of
  the {R}eynolds-averaged {N}avier--{S}tokes and k--$\omega$ turbulence model
  equations. {\it Computers \& Fluids} 2005\string; 34\string: 507--540.

\bibitem{fehn:2019}
Fehn N, Kronbichler M, Lehrenfeld C, Lube G, Schroeder P. High-order {DG}
  solvers for under-resolved turbulent incompressible flows: A comparison of
  $L^2$ and $H(div)$ methods. {\it International Journal of Numerical Methods
  in Fluids} 2019\string; 91\string: 533--556.

\bibitem{fehn:2017}
Fehn N, Wall W, Kronbichler M. On the stability of projection methods for the
  incompressible {N}avier--{S}tokes equations based on high-order discontinuous
  {G}alerkin discretizations. {\it Journal of Computational Physics}
  2017\string; 351\string: 392--421.

\bibitem{fehn:2018}
Fehn N, Wall W, Kronbichler M. Robust and efficient discontinuous {G}alerkin
  methods for under-resolved turbulent incompressible flows. {\it Journal of
  Computational Physics} 2018\string; 372\string: 667--693.

\bibitem{giorgiani:2014}
Giorgiani G, Fern{\'a}ndez-M{\'e}ndez S, Huerta A. Hybridizable discontinuous
  {G}alerkin with degree adaptivity for the incompressible {N}avier--{S}tokes
  equations. {\it Computers \& Fluids} 2014\string; 98\string: 196--208.

\bibitem{schotzau:2003}
Sch{\"o}tzau D, Schwab C, Toselli A. Stabilized {D}{G}{F}{E}{M} for
  incompressible flows. {\it Mathematical Models and Methods in Applied
  Sciences} 2003\string; 13\string: 1413--1436.

\bibitem{tumolo:2015}
Tumolo G, Bonaventura L. A semi-implicit, semi-{L}agrangian discontinuous
  {G}alerkin framework for adaptive numerical weather prediction: {SISL-D}G
  Framework for Adaptive {NWP}. {\it Quarterly Journal of the Royal
  Meteorological Society} 2015\string; 141\string: 2582--2601.

\bibitem{tumolo:2013}
Tumolo G, Bonaventura L, Restelli M. A semi-implicit, semi-{L}agrangian,
  $p-$adaptive discontinuous {G}alerkin method for the shallow water equations.
  {\it Journal of Computational Physics} 2013\string; 232\string: 46--67.

\bibitem{bank:1985}
Bank R, Coughran W, Fichtner W, Grosse E, Rose D, Smith R. {T}ransient
  {S}imulation of {S}ilicon {D}evices and {C}ircuits. {\it IEEE {T}ransactions
  on {E}lectron {D}evices.} 1985\string; 32\string: 1992-2007.

\bibitem{hosea:1996}
Hosea M, Shampine L. Analysis and implementation of {TR}-{BDF}2. {\it Applied
  Numerical Mathematics} 1996\string; 20\string: 21--37.

\bibitem{bangerth:2007}
Bangerth W, Hartmann R, Kanschat G. deal {II}: a general-purpose
  object-oriented finite element library. {\it ACM Transactions on Mathematical
  Software (TOMS)} 2007\string; 33\string: 24--51.

\bibitem{bassi:2007}
Bassi F, Crivellini A, Di~Pietro D, Rebay S. An implicit high-order
  discontinuous {G}alerkin method for steady and unsteady incompressible flows.
  {\it Computers \& Fluids} 2007\string; 36\string: 1529--1546.

\bibitem{bassi:2015}
Bassi F, Botti L, Colombo A, Ghidoni A, Massa F. Linearly implicit
  {R}osenbrock-type {R}unge--{K}utta schemes applied to the {D}iscontinuous
  {G}alerkin solution of compressible and incompressible unsteady flows. {\it
  Computers \& Fluids} 2015\string; 118.

\bibitem{rhebergen:2013}
Rhebergen S, Cockburn B, Van Der~Vegt J. A space--time discontinuous Galerkin
  method for the incompressible {N}avier-{S}tokes equations. {\it Journal of
  Computational Physics} 2013\string; 233\string: 339--358.

\bibitem{tavelli:2014}
Tavelli M, Dumbser M. {A staggered semi-implicit discontinuous Galerkin method
  for the two dimensional incompressible Navier--Stokes equations}. {\it
  Applied Mathematics and Computation} 2014\string; 248\string: 70--92.

\bibitem{tavelli:2016}
Tavelli M, Dumbser M. A staggered space--time discontinuous {G}alerkin method
  for the three-dimensional incompressible {N}avier--{S}tokes equations on
  unstructured tetrahedral meshes. {\it Journal of Computational Physics}
  2016\string; 319\string: 294--323.

\bibitem{chorin:1967}
Chorin A. A numerical method for solving incompressible viscous flow problems.
  {\it Journal of Computational Physics} 1967\string; 2\string: 12--26.

\bibitem{nithiarasu:2003}
Nithiarasu P. An efficient artificial compressibility ({AC}) scheme based on
  the characteristic based split ({CBS}) method for incompressible flows. {\it
  International Journal for Numerical Methods in Engineering} 2003\string;
  56\string: 1815--1845.

\bibitem{rogers:1990}
Rogers S, Kwak D. Upwind differencing scheme for the time-accurate
  incompressible {N}avier-{S}tokes equations. {\it AIAA journal} 1990\string;
  28\string: 253--262.

\bibitem{ekaterinaris:2004}
Ekaterinaris J. High-order accurate numerical solutions of incompressible flows
  with the artificial compressibility method. {\it International Journal of
  Numerical Methods in Fluids} 2004\string; 45\string: 1187--1207.

\bibitem{merkle:1987}
Merkle C. Time-accurate unsteady incompressible flow algorithms based on
  artificial compressibility. In:  {\it 8th Computational Fluid Dynamics
  Conference,}; 1987\string: 1137.

\bibitem{rahman:2008}
Rahman M, Siikonen T. An artificial compressibility method for viscous
  incompressible and low {M}ach number flows. {\it International Journal for
  Numerical Methods in Engineering} 2008\string; 75\string: 1320--1340.

\bibitem{casulli:1984}
Casulli V, Greenspan D. Pressure method for the numerical solution of
  transient, compressible fluid flows. {\it International Journal for Numerical
  Methods in Fluids} 1984\string; 4\string: 1001--1012.

\bibitem{dumbser:2016}
Dumbser M, Casulli V. A conservative, weakly nonlinear semi-implicit finite
  volume scheme for the compressible {N}avier-{S}tokes equations with general
  equation of state. {\it Applied Mathematics and Computation} 2016\string;
  272\string: 479--497.

\bibitem{giraldo:2013}
Giraldo F, Kelly J, Constantinescu E. Implicit-Explicit Formulations Of A
  Three-Dimensional Nonhydrostatic Unified Model Of The Atmosphere
  ({N}{U}{M}{A}). {\it SIAM Journal of Scientific Computing} 2013\string;
  35\string: 1162--1194.

\bibitem{bonaventura:2018b}
Bonaventura L, Fern{\'a}ndez-Nieto E, Garres-D{\'\i}az J, Narbona-Reina G.
  Multilayer shallow water models with locally variable number of layers and
  semi-implicit time discretization. {\it Journal of Computational Physics}
  2018\string; 364\string: 209--234.

\bibitem{garres:2021}
Garres-D{\'\i}az J, Bonaventura L. Flexible and efficient discretizations of
  multilayer models with variable density. {\it Applied Mathematics and
  Computation} 2021\string; 402\string: 126097.

\bibitem{bonaventura:2017}
Bonaventura L, Della~Rocca A. Unconditionally Strong Stability Preserving
  Extensions of the {TR-BDF2} Method. {\it Journal of Scientific Computing}
  2017\string; 70\string: 859--895.

\bibitem{kennedy:2016}
Kennedy C, Carpenter M. Diagonally Implicit {R}unge-{K}utta Methods for
  Ordinary Differential Equations, a Review. Tech. Rep. TM-2016-219173, NASA;
  2016.

\bibitem{bonaventura:2021}
Bonaventura L, G\'omez~Marmol M. The {TR-BDF} method for second order problems
  in structural mechanics. {\it Computers \& Mathematics with Applications}
  2021\string; 95\string: 13--26.

\bibitem{bell:1989}
Bell J, Colella P, Glaz H. A second-order projection method for the
  incompressible {N}avier-{S}tokes equations. {\it Journal of Computational
  Physics} 1989\string; 85\string: 257--283.

\bibitem{guermond:1998}
Guermond J, Quartapelle L. On incremental projection methods. {\it
  International Conference on Navier-Stokes Equations: Theory and Numerical
  Methods} 1998\string; 338\string: 277--288.

\bibitem{toselli:2002}
Toselli A. $H-p$ {D}iscontinuous {G}alerkin Approximations for the {S}tokes
  Problem. {\it Mathematical Models and Methods in Applied Sciences}
  2002\string; 12\string: 1565--1597.

\bibitem{brezzi:2002}
Arnold D, Brezzi F, Cockburn B, Marini L. Unified analysis of {D}iscontinuous
  {G}alerkin methods for elliptic problems. {\it SIAM Journal of Numerical
  Analysis} 2002\string; 39\string: 1749--1779.

\bibitem{arnold:1982}
Arnold D. An Interior Penalty Finite Element Method with Discontinuous
  Elements. {\it SIAM Journal of Numerical Analysis} 1982\string; 19\string:
  742--760.

\bibitem{green:1937}
Green A, Taylor G. Mechanism of the production of small eddies from large ones.
  {\it Proceedings of The Royal Society A: Mathematical, Physical and
  Engineering Sciences} 1937\string; 158\string: 499-521.

\bibitem{galloway:1987}
Galloway D, Frisch U. A note on the stability of a family of space-periodic
  {B}eltrami flows. {\it Journal of Fluid Mechanics} 1987\string; 158\string:
  557--564.

\bibitem{issa:1986}
Issa R, Ahmadi-Befrui B, Beshay K, Gosman A. Solution of the implicitly
  discretised reacting flow equations by operator-splitting. {\it Journal of
  Computational Physics} 1986\string; 62\string: 388-410.

\bibitem{auteri:2002}
Auteri F, Parolini N, Quartapelle L. Numerical Investigation on the Stability
  of Singular Driven Cavity Flow. {\it Journal of Computational Physics}
  2002\string; 183\string: 1-25.

\bibitem{botella:1998}
Botella O, Peyret R. Benchmark Spectral Results On The Lid-Driven Cavity Flow.
  {\it Computers and Fluids} 1998\string; 27\string: 421--433.

\bibitem{bruneau:2006}
Bruneau C, Saad M. The 2{D} lid-drivenvcavity flow revisited. {\it Computers
  and Fluids} 2006\string; 35\string: 326--348.

\bibitem{albensoeder:2005}
Albensoeder S, Kuhlmann H. Accurate three-dimensional lid-driven cavity flow.
  {\it Journal of Computational Physics} 2005\string; 206\string: 536--558.

\bibitem{jiang:1994}
Jiang B, Lin T, Povinelli L. Large-scale computation of incompressible viscous
  flow by least-squares finite element method. {\it Computer Methods in Applied
  Mechanics and Engineering} 1994\string; 114(3)\string: 213-231.

\bibitem{schafer:1996}
Sch{\"a}fer M, Turek S, Durst F, Krause E, Rannacher R. {\it Benchmark
  Computations of Laminar Flow Around a Cylinder}\string: 547--566; Wiesbaden:
  Vieweg+Teubner Verlag .
\newblock 1996

\bibitem{burner:2018}
{RFCS (Research Fund for Coal and Steel) research project ``Burner 4.0,
  Development of a new burner concept: Industry 4.0 technologies applied to the
  best available combustion system for the Steel Industry'', project ID
  847237}. tech. rep.,  2018.

\bibitem{georgoulis:2009}
Georgoulis E, Hall E, Houston P. {Discontinuous {G}alerkin methods on
  $hp$-anisotropic meshes II: A posteriori error analysis and adaptivity}. {\it
  Applied Numerical Mathematics} 2009\string; 59\string: 2179--2194.

\bibitem{hartmann:2006}
Hartmann R, Houston P. {Symmetric interior penalty DG methods for the
  compressible Navier--Stokes equations II: Goal--oriented a posteriori error
  estimation}. {\it International Journal of Numerical Analysis \& Modeling}
  2006\string; 3\string: 141--162.

\bibitem{oden:1994}
Oden T, Wu W, Ainsworth M. An a posteriori error estimate for finite element
  approximations of the {Navier-Stokes} equations. {\it Computer Methods in
  Applied Mechanics and Engineering} 1994\string; 111\string: 185--202.

\bibitem{bonaventura:2020a}
Bonaventura L, Casella F, Delpopolo~Carciopolo L, Ranade A. A self adjusting
  multirate algorithm for robust time discretization of partial differential
  equations. {\it Computers and Mathematics with Applications} 2020\string;
  79\string: 2086--2098.

\bibitem{bassi:2006}
Bassi F, Crivellini A, Di~Pietro D, Rebay S. An artificial compressibility flux
  for the discontinuous {G}alerkin solution of the incompressible
  {N}avier--{S}tokes equations. {\it Journal of Computational Physics}
  2006\string; 218\string: 794--815.

\bibitem{chorin:1993}
Chorin A, Marsden J. {\it A mathematical introduction to fluid mechanics. Third
  edition}.
\newblock Springer .
\newblock 1993.

\bibitem{ham:2004}
Ham F, Iaccarino G. Energy conservation in collocated discretization schemes on
  unstructured meshes. tech. rep., Center for Turbulence Research, NASA Ames -
  Stanford University;  2004.

\bibitem{quarteroni:1994}
Quarteroni A, Valli A. {\it Numerical approximation of partial differential
  equations}.
\newblock Springer Verlag .
\newblock 1994.

\bibitem{shen:1995}
Shen J. On error estimates of the penalty method for unsteady {N}avier-{S}tokes
  equations. {\it SIAM Journal on Numerical Analysis} 1995\string; 32\string:
  386--403.

\bibitem{shen:1997}
Shen J. Pseudo-compressibility methods for the unsteady incompressible
  {N}avier-{S}tokes equations. In:  {\it Proceedings of the 1994 Beijing
  symposium on nonlinear evolution equations and infinite dynamical
  systems,}ZhongShan University Press. ; 1997\string: 68--78.

\end{thebibliography}

\begin{table}[h!]
	\centering
	\begin{tabular}{|c|c|c|c|c|c|c|}
		\hline
		$\Delta t $   & $N_{el} $ & $\mu$ & \(H^1\) rel. error \(\mathbf{u}\) & \(H^1\)  rate \(\mathbf{u}\) & \(L^2\) rel. error \(\mathbf{u}\) & \(L^2\) rate \(\mathbf{u}\) \\
		\hline
		0.64    &   8 & 0.04 &          1.5          & &                 0.38 & \\
		\hline
		0.32   &   16 & 0.08 &       0.65   &     1.22       &    0.095     & 2.01 \\
		\hline
		0.16    &  32 & 0.17 &    0.12   &   2.45      &      0.016  &    2.58 \\
		\hline
		0.08  &    64 & 0.33 &      0.023 &       2.38 &           0.0031  &    2.37 \\
		\hline
	\end{tabular}
	\caption{Convergence test for the Green-Taylor vortex benchmark computed at  \(C = 1.63\) with  $\mathbf{Q}_2-Q_1$ elements, relative errors for the velocity in $H^1$ and $L^2$ norms.}
	\label{tab:greentaylor_q2q1_u}
\end{table}

\begin{table}[h!]
	\centering
	\begin{tabular}{|c|c|c|c|c|}
		\hline
		$\Delta t $   & $N_{el} $ & $\mu$ & \(L^2\) rel. error \(p\) & \(L^2\)  rate \(p\)  \\
		\hline
		0.64    &   8 & 0.04 &          0.43          & \\
		\hline
		0.32   &   16 & 0.08 &       0.14   &     1.60       \\
		\hline
		0.16    &  32 & 0.17 &    0.04   &   1.72      \\
		\hline
		0.08  &    64 & 0.33 &      0.011 &       1.91  \\
		\hline
	\end{tabular}
	\caption{Convergence test for the Green-Taylor vortex benchmark computed at  \(C = 1.63\) with  $\mathbf{Q}_2-Q_1$ elements, relative errors for the pressure in $L^2$ norm.}
	\label{tab:greentaylor_q2q1_p}
\end{table}

\begin{table}[h!]
	\centering
	\begin{tabular}{|c|c|c|c|c|c|c|}
		\hline
		$\Delta t $   & $N_{el} $ & $\mu$ & \(H^1\) rel. error \(\mathbf{u}\) & \(H^1\)  rate \(\mathbf{u}\) & \(L^2\) rel. error \(\mathbf{u}\) & \(L^2\) rate \(\mathbf{u}\) \\
		\hline
		0.43    &   8 & 0.06  &    0.28                & &                             0.062 & \\
		\hline
		0.21   &   16 & 0.12  &    0.033     &     3.12       &   0.0068        & 3.18 \\
		\hline
		0.11    &  32 & 0.25  &   0.0044     &   2.88      &      0.00044  &    3.93 \\
		\hline
		0.053  &    64 &  0.50 &   0.00059    &       2.92 &       0.000031       &    3.85 \\
		\hline
	\end{tabular}
	\caption{Convergence test for the Green-Taylor vortex benchmark computed at  \(C = 1.63\) with  $\mathbf{Q}_3-Q_2$ elements, relative errors for the velocity in $H^1$ and $L^2$ norms.}
	\label{tab:greentaylor_q3q2_u}
\end{table}

\begin{table}[h!]
	\centering
	\begin{tabular}{|c|c|c|c|c|}
		\hline
		$\Delta t $   & $N_{el} $ & $\mu$ & \(L^2\) rel. error \(p\) & \(L^2\)  rate \(p\) \\
		\hline
		0.43    &   8 & 0.06  &    0.087                & \\
		\hline
		0.21   &   16 & 0.12  &    0.011     &     2.93        \\
		\hline
		0.11    &  32 & 0.25  &   0.00075     &   3.92      \\
		\hline
		0.053  &    64 &  0.50 &   0.000029    &       4.72\\
		\hline
	\end{tabular}
	\caption{Convergence test for the Green-Taylor vortex benchmark computed at  \(C = 1.63\) with  $\mathbf{Q}_3-Q_2$ elements, relative errors for the pressure in $L^2$ norm.}
	\label{tab:greentaylor_q3q2_p}
\end{table}

\begin{table}[h!]
	\centering
	\begin{tabular}{|c|c|c|c|c|c|c|}
		\hline
		$\Delta t $   & $N_{el} $ & $\mu$ & \(H^1\) rel. error \(\mathbf{u}\) & \(H^1\)  rate \(\mathbf{u}\) & \(L^2\) rel. error \(\mathbf{u}\)& \(L^2\) rate \(\mathbf{u}\) \\
		\hline
		0.52    &   8 & 0.05 &          1.2          & &                 0.32 & \\
		\hline
		0.22   &   16 & 0.12  &       0.55   &     1.13       &    0.081     & 1.96 \\
		\hline
		0.11    &  32 & 0.25 &    0.12   &   2.16     &      0.012  &    2.77 \\
		\hline
		0.052  &    64 & 0.50  &      0.021 &       2.51 &           0.0023 &    2.41 \\
		\hline
	\end{tabular}
	\caption{Convergence test for the Green-Taylor vortex benchmark computed on a distorted mesh at  \(C = 1.63\) with  $\mathbf{Q}_2-Q_1$ elements, relative errors for the velocity in $H^1$ and $L^2$ norms.}
	\label{tab:greentaylor_q2q1_dist_u}
\end{table}

\begin{table}[h!]
	\centering
	\begin{tabular}{|c|c|c|c|c|}
		\hline
		$\Delta t $   & $N_{el} $ & $\mu$ & \(L^2\) rel. error \(p\)& \(L^2\) rate \(p\) \\
		\hline
		0.52    &   8 & 0.05 &          0.43          &  \\
		\hline
		0.22   &   16 & 0.12  &       0.077   &     2.02       \\
		\hline
		0.11    &  32 & 0.25 &    0.024   &   1.68      \\
		\hline
		0.052  &    64 & 0.50  &      0.0064 &       1.91 \\
		\hline
	\end{tabular}
	\caption{Convergence test for the Green-Taylor vortex benchmark computed on a distorted mesh at  \(C = 1.63\) with  $\mathbf{Q}_2-Q_1$ elements, relative errors for the pressure in $L^2$ norm.}
	\label{tab:greentaylor_q2q1_dist_p}
\end{table}

\begin{table}[h!]
	\centering
	\begin{tabular}{|c|c|c|c|c|c|c|}
		\hline
		$\Delta t $   & $N_{el} $ & $\mu$ & \(H^1\) rel. error \(\mathbf{u}\) & \(H^1\)  rate \(\mathbf{u}\) & \(L^2\) rel. error \(\mathbf{u}\) & \(L^2\) rate \(\mathbf{u}\) \\
		\hline
		0.64    &   8 & 0.04 &          1.49          & &                 0.48 & \\
		\hline
		0.32   &   16 & 0.08 &       0.73   &     1.04       &    0.13     & 1.84 \\
		\hline
		0.16    &  32 & 0.17 &    0.14   &   2.34      &     0.0311  &    2.11 \\
		\hline
		0.08  &    64 & 0.33 &      0.03 &       2.25 &           0.0070  &    2.16 \\
		\hline
	\end{tabular}
	\caption{Convergence test for the Green-Taylor vortex benchmark computed at  \(C = 1.63\) with  $\mathbf{Q}_2-Q_1$ elements and the projection method of \cite{bell:1989}, relative errors for the velocity in $H^1$ and $L^2$ norms.}
	\label{tab:greentaylor_q2q1_bell_u}
\end{table}

\begin{table}[h!]
	\centering
	\begin{tabular}{|c|c|c|c|c|}
		\hline
		$\Delta t $   & $N_{el} $ & $\mu$ & \(L^2\) rel. error \(p\) & \(L^2\) rate \(p\) \\
		\hline
		0.64    &   8 & 0.04 &          0.41          & \\
		\hline
		0.32   &   16 & 0.08 &       0.12   &     1.75        \\
		\hline
		0.16    &  32 & 0.17 &    0.037   &   1.75       \\
		\hline
		0.08  &    64 & 0.33 &      0.0087 &       2.08  \\
		\hline
	\end{tabular}
	\caption{Convergence test for the Green-Taylor vortex benchmark computed at  \(C = 1.63\) with  $\mathbf{Q}_2-Q_1$ elements and the projection method of \cite{bell:1989}, relative errors for the pressure in $L^2$ norm.}
	\label{tab:greentaylor_q2q1_bell_p}
\end{table}

\begin{table}[h!]
	\centering
	\begin{tabular}{|c|c|c|c|c|c|c|}
		\hline
		$\Delta t $   & $N_{el} $ & $\mu$ & \(H^1\) rel. error \(\mathbf{u}\) & \(H^1\)  rate \(\mathbf{u}\) & \(L^2\) rel. error \(\mathbf{u}\) & \(L^2\) rate \(\mathbf{u}\)\\
		\hline
		0.64    &   8 & 0.04 &          0.89          & &                 0.28 & \\
		\hline
		0.32   &   16 & 0.08 &       0.40   &     1.15       &    0.09      & 1.69 \\
		\hline
		0.16    &  32 & 0.17 &    0.085   &   2.24      &      0.023  &    1.92 \\
		\hline
		0.08  &    64 & 0.33 &      0.029 &       1.57 &           0.0077 &    1.60 \\
		\hline
	\end{tabular}
	\caption{Convergence test for the Green-Taylor vortex benchmark computed at  \(C = 1.63\) with  $\mathbf{Q}_2-Q_1$ elements and the projection method of \cite{guermond:1998}, relative errors for the velocity in $H^1$ and $L^2$ norms.}
	\label{tab:greentaylor_q2q1_guermond_u}
\end{table}

\begin{table}[h!]
	\centering
	\begin{tabular}{|c|c|c|c|c|}
		\hline
		$\Delta t $   & $N_{el} $ & $\mu$ & \(L^2\) rel. error \(p\) & \(L^2\) rate \(p\)\\
		\hline
		0.64    &   8 & 0.04 &          0.41          &  \\
		\hline
		0.32   &   16 & 0.08 &       0.10   &     2.06        \\
		\hline
		0.16    &  32 & 0.17 &    0.026   &     1.91 \\
		\hline
		0.08  &    64 & 0.33 &      0.0067   &  1.97 \\
		\hline
	\end{tabular}
	\caption{Convergence test for the Green-Taylor vortex benchmark computed at  \(C = 1.63\) with  $\mathbf{Q}_2-Q_1$ elements and the projection method of \cite{guermond:1998}, relative errors for the pressure in $L^2$ norm.}
	\label{tab:greentaylor_q2q1_guermond_p}
\end{table}

\begin{table}[h!]
	\centering
	\begin{tabular}{|c|c|c|c|c|c|c|c|}
		\hline
		$\Delta t $   & $N_{el} $ & $\mu$  & \(H^1\) rel. error \(\mathbf{u}\) & \(H^1\)  rate \(\mathbf{u}\) & \(L^2\) rel. error \(\mathbf{u}\) & \(L^2\) rate \(\mathbf{u}\) \\
		\hline
		1.18    &   8 & 0.08 &          1.33          & &                 0.39 & \\
		\hline
		0.59   &   16 & 0.15 &       0.63   &     1.07       &    0.11      & 1.79 \\
		\hline
		0.29    &  32 & 0.31 &    0.12   &   2.35      &      0.028  &    2.02 \\
		\hline
		0.15  &    64 & 0.61 &      0.028 &       2.17 &           0.0059 &    2.23 \\
		\hline
	\end{tabular}
	\caption{Convergence test for the Green-Taylor vortex benchmark computed at  \(C = 3\) with  $\mathbf{Q}_2-Q_1$ elements, relative errors for the velocity in $H^1$ and $L^2$ norms.}
	\label{tab:greentaylor_q2q1_CFL_3_u}
\end{table}

\begin{table}[h!]
	\centering
	\begin{tabular}{|c|c|c|c|c|c|}
		\hline
		$\Delta t $   & $N_{el} $ & $\mu$  & \(L^2\) rel. error \(p\) & \(L^2\) rate \(p\) \\
		\hline
		1.18    &   8 & 0.08 &          0.49          & \\
		\hline
		0.59   &   16 & 0.15 &       0.13   &     1.87        \\
		\hline
		0.29    &  32 & 0.31 &    0.04   &   1.60       \\
		\hline
		0.15  &    64 & 0.61 &      0.013 &       1.75  \\
		\hline
	\end{tabular}
	\caption{Convergence test for the Green-Taylor vortex benchmark computed at  \(C = 3\) with  $\mathbf{Q}_2-Q_1$ elements, relative errors for the pressure in $L^2$ norm.}
	\label{tab:greentaylor_q2q1_CFL_3_p}
\end{table}

\begin{table}[h!]
	\centering
	\begin{tabular}{|c|c|c|c|c|c|c|}
		\hline
		$\Delta t $   & $N_{el} $ & $\mu$ & \(H^1\) rel. error \(\mathbf{u}\) & \(H^1\) rate \(\mathbf{u}\) & \(L^2\) rel. error \(\mathbf{u}\) & \(L^2\) rate \(\mathbf{u}\)\\
		\hline
		0.32    &   8 & 2.08 &          0.019           & &                 0.0078 & \\
		\hline
		0.16   &   16 & 4.15 &       0.0054    &     1.85      &    0.0022     & 1.86 \\
		\hline
		0.08    &  32 & 8.30  &    0.0014    &   1.98      &      0.00056  &    1.99 \\
		\hline
		0.04  &    64 & 16.60  &      0.00036 &       1.91 &           0.00017 &    1.75 \\
		\hline
	\end{tabular}
	\caption{Convergence test for the ABC flow benchmark computed at  \(C = 1.63\) with  $\mathbf{Q}_2-Q_1$ elements, relative errors for the velocity in $H^1$ and $L^2$ norms.}
	\label{tab:abc_q2q1_u}
\end{table}

\begin{table}[h!]
	\centering
	\begin{tabular}{|c|c|c|c|c|}
		\hline
		$\Delta t $   & $N_{el} $ & $\mu$ & \(L^2\) rel. error \(p\) & \(L^2\) rate \(p\)\\
		\hline
		0.32    &   8 & 2.08 &          1.0           &  \\
		\hline
		0.16   &   16 & 4.15 &       0.13    &     2.93       \\
		\hline
		0.08    &  32 & 8.30  &    0.039    &   1.74       \\
		\hline
		0.04  &    64 & 16.60  &      0.011 &       1.79 \\
		\hline
	\end{tabular}
	\caption{Convergence test for the ABC flow benchmark computed at  \(C = 1.63\) with  $\mathbf{Q}_2-Q_1$ elements, relative errors for the pressure in $L^2$ norm.}
	\label{tab:abc_q2q1_p}
\end{table}

\begin{table}[h!]
	\centering
	\begin{tabular}{|c|c|c|c|c|c|c|}
		\hline
		$\Delta t $   & $N_{el} $ & $\mu$ & \(H^1\) rel. error \(\mathbf{u}\) & \(H^1\) rate \(\mathbf{u}\) & \(L^2\) rel. error \(\mathbf{u}\) & \(L^2\) rate \(\mathbf{u}\) \\
		\hline
		0.21    &   8 & 3.11  &    0.0036               & &                             0.0019 & \\
		\hline
		0.11   &   16 & 6.23  &    0.0010     &     1.80       &   0.0068        & 2.05 \\
		\hline
		0.053    &  32 & 12.45  &   0.00037     &   1.5      &     0.00014  &    1.68 \\
		\hline
	\end{tabular}
	\caption{Convergence test for the ABC flow benchmark computed at  \(C = 1.63\) with  $\mathbf{Q}_3-Q_2$ elements, relative errors for the velocity in $H^1$ and $L^2$ norms.}
	\label{tab:abc_q3q2_u}
\end{table}

\begin{table}[h!]
	\centering
	\begin{tabular}{|c|c|c|c|c|}
		\hline
		$\Delta t $   & $N_{el} $ & $\mu$ & \(L^2\) rel. error \(p\) & \(L^2\) rate \(p\) \\
		\hline
		0.21    &   8 & 3.11  &    0.25               &  \\
		\hline
		0.11   &   16 & 6.23  &    0.033     &     2.93        \\
		\hline
		0.053    &  32 & 12.45  &   0.0097     &   1.72       \\
		\hline
	\end{tabular}
	\caption{Convergence test for the ABC flow benchmark computed at  \(C = 1.63\) with  $\mathbf{Q}_3-Q_2$ elements, relative errors for the pressure in $L^2$ norm.}
	\label{tab:abc_q3q2_p}
\end{table}

\begin{table}[h!]
	\centering
	\begin{tabular}{|c|c|c|c|c|c|c|}
		\hline
		$\Delta t $   & $N_{el} $ & $C$ & \(H^1\) rel. error \(\mathbf{u}\)& \(H^1\)  rate \(\mathbf{u}\)& \(L^2\) rel. error \(\mathbf{u}\)& \(L^2\) rate \(\mathbf{u}\)\\
		\hline
		0.32    &   8 & 1.57 &          0.019           & &                  0.0071 & \\
		\hline
		0.08   &   16 & 0.79 &       0.0045    &     2.05      &    0.0013     & 2.5 \\
		\hline
		0.02    &  32 & 0.39  &    0.0012    &   1.97      &      0.00031  &    2.02 \\
		\hline
		0.005  &    64 & 0.20  &      0.00029 &       1.98 &           0.000053  &    2.54 \\
		\hline
	\end{tabular}
	\caption{Convergence test for the ABC flow benchmark computed at  \(\mu = 2\) with  $\mathbf{Q}_2-Q_1$ elements, relative errors for the velocity in $H^1$ and $L^2$ norms.}
	\label{tab:abc_q2q1_muc_u}
\end{table}

\begin{table}[h!]
	\centering
	\begin{tabular}{|c|c|c|c|c|}
		\hline
		$\Delta t $   & $N_{el} $ & $C$ & \(L^2\) rel. error \(p\)& \(L^2\) rate \(p\)\\
		\hline
		0.32    &   8 & 1.57 &          1.0           & \\
		\hline
		0.08   &   16 & 0.79 &       0.16    &     2.66      \\
		\hline
		0.02    &  32 & 0.39  &    0.042    &   1.93       \\
		\hline
		0.005  &    64 & 0.20  &      0.011 &       1.94  \\
		\hline
	\end{tabular}
	\caption{Convergence test for the ABC flow benchmark computed at  \(\mu = 2\) with  $\mathbf{Q}_2-Q_1$ elements, relative errors for the pressure in $L^2$ norm.}
	\label{tab:abc_q2q1_muc_p}
\end{table}

\begin{table}[h!]
	\centering
	\begin{tabular}{|c|c|c|c|c|c|c|}
		\hline
		$\Delta t $   & $N_{el} $ & $C$ & \(H^1\) rel. error \(\mathbf{u}\)& \(H^1\)  rate \(\mathbf{u}\)& \(L^2\) rel. error \(\mathbf{u}\) & \(L^2\) rate \(\mathbf{u}\)\\
		\hline
		0.14    &   8 & 1.05  &    0.0025               & &                             0.00089 & \\
		\hline
		0.034   &   16 & 0.52  &    0.00024     &     2.70       &   0.00011        & 3.08 \\
		\hline
		0.0086    &  32 & 0.26  &   0.000071     &   1.78      &      0.000018  &    2.51 \\
		\hline
	\end{tabular}
	\caption{Convergence test for the ABC flow benchmark computed at  \(\mu = 2\) with  $\mathbf{Q}_3-Q_2$ elements, relative errors for the velocity in $H^1$ and $L^2$ norms.}
	\label{tab:abc_q3q2_muc_u}
\end{table}

\begin{table}[h!]
	\centering
	\begin{tabular}{|c|c|c|c|c|}
		\hline
		$\Delta t $   & $N_{el} $ & $C$ & \(L^2\) rel. error \(p\) & \(L^2\) rate \(p\)\\
		\hline
		0.14    &   8 & 1.05  &    0.20               & \\
		\hline
		0.034   &   16 & 0.52  &    0.027     &     2.89        \\
		\hline
		0.0086    &  32 & 0.26  &   0.0043     &   2.65       \\
		\hline
	\end{tabular}
	\caption{Convergence test for the ABC flow benchmark computed at  \(\mu = 2\) with  $\mathbf{Q}_3-Q_2$ elements, relative errors for the pressure in $L^2$ norm.}
	\label{tab:abc_q3q2_muc_p}
\end{table}

\begin{table}[h!]
	\centering
	\begin{tabular}{|c|c|c|c|c|}
		\hline
		& OpenFOAM (coarse)  & OpenFOAM (middle) & OpenFOAM (fine) & deal.II \\
		\hline
		line \(A\)    & 18.6801 & 18.8735 &  18.9112  &  18.7687 \\
		\hline
		line \(B\)   &   18.5748 & 18.7453 &      18.7973    &     18.6780      \\
		\hline
		line \(C\)    &  18.2494 & 18.4158  &    18.4706    &   18.3596       \\
		\hline
		line \(D\)  &    17.0799 & 17.2534  &      17.3135 &       17.2452  \\
		\hline
	\end{tabular}
	\caption{Pressure drop along the four midlines of the channels for the different simulations.}
	\label{tab:exchanger_pressure_drop}
\end{table}

\begin{table}[h!]
	\centering
	\begin{tabular}{|c|c|c|c|}
		\hline
		Number of cores & Wallclock time TR-BDF2  & Wallclock time Bell-Colella-Glaz & Wallclock time Guermond-Qaurtapelle BDF2  \\
		\hline
		16    & \(1.86 \cdot 10^{3}\) & \(3.80 \cdot 10^{3}\) &  \(8.09 \cdot 10^{2}\)   \\
		\hline
		32    & \(7.39 \cdot 10^{2}\) & \(1.75 \cdot 10^{3}\) &  \(3.64 \cdot 10^{2}\)   \\
		\hline
		64   &   \(3.69 \cdot 10^{2}\) &  \(8.54 \cdot 10^{2}\) &  \(1.79 \cdot 10^{2}\)    \\
		\hline
		128    &  \(1.81 \cdot 10^{2}\) & \(4.29 \cdot 10^{2}\)  &    \(8.96 \cdot 10^{1}\)     \\
		\hline
		256  &    \(1.08 \cdot 10^{2}\) & \(2.29 \cdot 10^{2}\)  &      \(5.01 \cdot 10^{1}\)   \\
		\hline
		512  &    \(6.24 \cdot 10^{1}\) & \(1.32 \cdot 10^{2}\)  &      \(3.16 \cdot 10^{1}\)   \\
		\hline
		1024  &    \(3.46 \cdot 10^{1}\) & \(6.07 \cdot 10^{1}\)  &      \(1.62 \cdot 10^{1}\)   \\
		\hline
	\end{tabular}
	\caption{Wallclock times in seconds of the different simulations performed for the strong scaling analysis.}
	\label{tab:scability_times}
\end{table}

\end{document}